\documentclass[a4paper,10pt]{article}

\usepackage[utf8]{inputenc}
\usepackage[T1]{fontenc} 
\usepackage[english]{babel}
\usepackage{lmodern}
\usepackage{pifont}
\usepackage[plain]{fullpage}
\usepackage{authblk} 
\usepackage{amsmath,amsfonts,amssymb,dsfont,bm,mathtools,amsthm}
\usepackage{graphicx}
\graphicspath{{Figures/}}
\usepackage{xcolor}
\usepackage{csquotes}
\usepackage{natbib}
\usepackage[normalem]{ulem}
\usepackage{accents}
\usepackage{enumitem}   

\usepackage[hyperfootnotes=false]{hyperref}
\hypersetup{
  colorlinks,
  citecolor=blue,
  linkcolor=red,
  urlcolor=gray}

\DeclareMathOperator{\MSPE}{MSPE}

\DeclareMathOperator{\Tr}{Tr}
\DeclareMathOperator{\Span}{Span}
\DeclareMathOperator{\Id}{I}
\DeclareMathOperator{\Ind}{Ind}
\DeclareMathOperator{\HS}{HS}
\DeclareMathOperator{\supp}{supp}
\DeclareMathOperator*{\argmin}{arg\,min}

\DeclareMathOperator{\KL}{\mathbf{KL}}
\DeclareMathOperator{\pen}{pen}
\DeclareMathOperator{\EMSPE}{EMSPE}

\newcommand{\ubar}[1]{\underaccent{\bar}{#1}}
\def\RR{\textsf{R}\/}
\newtheorem{remark}{Remark}

\newtheorem{proposition}{Proposition}
\newtheorem{corollary}{Corollary}
\newtheorem{lemma}{Lemma}
\newtheorem{theorem}{Theorem}

\begin{document}

\title{Adaptive nonparametric estimation in the functional linear model with functional output}

\author[1]{Gaëlle Chagny\footnote{gaelle.chagny@univ-rouen.fr}}
\author[1]{Anouar Meynaoui\footnote{anouar.meynaoui@gmail.com}}
\author[2]{Angelina Roche\footnote{roche@ceremade.dauphine.fr}}
\affil[1]{LMRS, UMR CNRS 6085, Université de Rouen Normandie}
\affil[2]{CEREMADE, UMR CNRS 7534, Université Paris-Dauphine}

\date{\today}
\maketitle

\begin{abstract}
In this paper, we consider a functional linear regression model, where both the covariate and the response variable are functional random variables. We address the problem of optimal nonparametric estimation of the conditional expectation operator in this model. A collection of projection estimators over finite dimensional subspaces is first introduce. We provide a non-asymptotic bias-variance decomposition for the Mean Square Prediction error in the case where these subspaces are generated by the (empirical) PCA functional basis. The automatic trade-off is realized thanks to a model selection device which selects the best projection dimensions: the penalized contrast estimator satisfies an oracle-type inequality and is thus optimal in an adaptive point of view. These upper-bounds allow us to derive convergence rates over ellipsoidal smoothness spaces. The rates are shown to be optimal in the minimax sense: they match with a lower bound of the minimax risk, which is also proved. Finally, we conduct a numerical study, over simulated data and over two real-data sets. 
\end{abstract}

\section{Introduction}
\label{Introduction}
Functional data analysis \citep{RS05,FV06,FR11} has attracted a growing interest from the past decades. 
In this context, regression models involving functional data as covariate are of particular interest. The case where the variable to predict is a real variable, called \emph{functional linear model with scalar output} or simply \emph{functional linear model} has been widely studied (see e.g. \citealt{cai2006prediction,cardot2007clt,LH07,hilgert2013minimax,cai2012minimax}) and is now well understood. In particular, the minimax rates for the estimation of the slope function in this model have been computed by \citet{cardot2010thresholding} and adaptive estimators have been built \citep{comte2010adaptive,comte2012adaptive,BR15,brunel2016non}. On the contrary, the case of the \emph{functional linear model with functional output}, where the variable to predict is also a functional variable has been less studied. This paper is dedicated to minimax adaptive estimation in this framework.

We assume here that we observe a sample $\{(X_i,Y_i), i=1,\hdots,n\}$, $n\in\mathbb N\backslash\{0\}$ of independent copies of a couple of functional data $(X,Y)$. For simplicity, we assume that both $X$ and $Y$ are random variables in the same functional space $\mathbb H=L^2([0,1])$, the space of square integrable functions on the interval $[0,1]$, equipped with its usual scalar product $\langle \cdot,\cdot\rangle$ defined by $\langle f,g\rangle=\int_0^1f(t)g(t)dt$, $f,g\in\mathbb H$ and norm $\|\cdot\|=\sqrt{\langle\cdot,\cdot\rangle}$.  
The link between the functional variable of interest $Y\in \mathbb H$ and the functional covariate $X\in\mathbb H$ is linear: there exists an operator $S\in\mathcal L(\mathbb H)$, the space of continuous linear operators on $\mathbb H$, such that 
\begin{equation}\label{eq:model}
Y=SX+\varepsilon,
\end{equation}
where $\varepsilon\in\mathbb H$ stands for an (unobserved) noise. The functional variables $X$ and $\varepsilon$ are supposed to be both centered, and independent. The noise $\varepsilon$  satisfied $\sigma_{\varepsilon}^2=\mathbb E\|\varepsilon\|^2<\infty$. 
The \textit{slope operator $S$} is an integral operator and we denote by $\mathcal S\in L^2([0,1]^2)$ its kernel:
$$\begin{array}{ccl}\displaystyle  S :\mathbb H&\longrightarrow &\mathbb H\\
  \displaystyle                             f  &\longmapsto     & Sf\, :\, t\in [0,1]\,\mapsto\, Sf(t)=\displaystyle\int_0^1\mathcal S(s,t)f(s)ds.\end{array}$$
The aim is to estimate the unknown operator $S$ (or its kernel $\mathcal S$)  from the sample $(X_i,Y_i)_{i\in\{1,\ldots,n\}}$.
 
It seems that the first article about estimation in this model is the one of \citet{CFF02}. In the fixed design case, they propose a histogram estimator and prove its consistency under strong assumption on the design matrix. A wavelets estimator has been considered by \citet{AOV08} and a splines estimator by \citet{APDeRS10}. The majority of the literature focus on estimators by projection onto the basis of principal components of the covariate $X$ \citep{CMW04,YMW05}. 
The interest of functional Principal Components Analysis (PCA in the sequel) may be seen in the Karhunen-Lo\`eve decomposition of $X$ that is to say the writing of $X$ as a series (with convergence in $\mathbb H$) 
\begin{equation}\label{eq:KL}
X = \sum_{j\geq 1}\sqrt{\lambda_j}\xi_j\varphi_j,
\end{equation}
where $(\lambda_j)_{j\geq 1}$ is a non-increasing sequence of non-negative real numbers, $(\xi_j)_{j\geq 1}$ is a sequence of standardized random variables (the \emph{principal components scores}) and $(\varphi_j)_{j\geq 1}$ is an orthonormal basis of $\mathbb H$ (the \emph{principal components basis}). It can be proved that, for a given dimension $D$, the space $\rm{span}\{\varphi_1,\hdots,\varphi_D\}$ is the best approximation space for $X$ in the sense of the $L^2$-loss i.e. 
\[
\rm{span}\{\varphi_1,\hdots,\varphi_D\}={\arg\min}\left\{\mathbb E\left[\|X-{\rm proj}_S(X)\|^2\right], S\text{ lin. sub. of } \mathbb H, \dim(S)=D\right\}
\]
where, for a linear subspace $S$ of $\mathbb H$, ${\rm proj}_S(X)$ is the orthogonal projection of $X$ into $S$, see \citet[Chapter 8]{FR11} or \citet[Theorem 7.2.8]{hsing_theoretical_2015}. We also refer to \citet{dauxois_asymptotic_1982,MR15}  for other reviews on PCA for functional data. 
A procedure to estimate the $\varphi_j$'s is described in Section~\ref{sec:estimPCA}.  

To our knowledge, few articles investigate the theoretical properties of slope operator or kernel estimators in Model \eqref{eq:model}. \citet{crambes2013asymptotics} study an estimator of the slope operator $S$ by projection onto the principal components basis. They provide a bias-variance decomposition of the mean squared prediction risk and compute optimal rates of convergence: such type of results can be stated only under some smoothness assumptions on the target operator $S$ (as usual in nonparametric estimation)  but also under assumptions of the process $X$, through the rate of decay of the eigenvalues (the $\lambda_j$'s in \eqref{eq:KL}) of the associated covariance operator. \citet{crambes2013asymptotics} also derive weak convergence properties of their estimator. In their procedure, the smoothness indices of the target operator and of the covariate $X$ (the decreasing rate of the covariance operator eigenvalues for example) are required to choose the projection dimension that permits to achieve the optimal rate. Thus, the method is not adaptive.
More recently, \citet{imaizumi2018pca} study two estimation procedures called \emph{simple} and \emph{double} truncature. The simple truncature estimator corresponds to the one of \citet{crambes2013asymptotics}. They obtain lower and asymptotic upper-bounds on the estimation risk of the slope kernel $\mathcal S$. As in \citet{crambes2013asymptotics}, the procedure is not adaptive, and the results are valid only when the decay rate of the eigenvalues of the covariance operator is a polynomial one.

In the present work, we propose a procedure which leads to an optimal estimate for the slope operator in Model \eqref{eq:model}, both from the minimax and the adaptive estimation point of view, for the mean squared prediction error. We first introduce a collection of projection estimators, by minimizing a least-squares contrast function over subspaces of $\mathbb H$ spanned by the first elements of the PCA basis, corresponding to the \textit{double truncature} procedure of \cite{imaizumi2018pca}. We focus on the mean squared prediction error and compute a non-asymptotic upper-bound in Theorem \ref{thm-UB-Emp-Risk}. This bound exhibits a bias-variance decomposition allowing us to derive rates of convergence, under some regularity assumption on the operator $S\Gamma^{1/2}$ (see Corollary \ref{corr-upper-bound-MSPE}).
We then show that these bounds match with lower bounds that we also proved in Theorem \ref{thm-lower-bound-MSPE} (see also Corollary \ref{coro:minimax_rates}). One of the other main original contributions of this paper is to propose an entirely data-driven procedure to automatically select the best projection dimensions. The method relies on classical model selection tools \citep{Mas2007}, and takes advantage of the definition of the estimates as minimized-contrast estimators. A penalized version of the contrast function permits to derive data driven estimate, which satisfies an oracle-type inequality, and achieves the optimal minimax convergence rates. Our selection rule does not depend on smoothness parameters of $X$ and $S$.  The procedure is then adaptive and minimax optimal. 

The paper is organized as follows. Our least-squares estimators are constructed in Section \ref{Estimation}. Upper and lower-bounds for the risk we choose are established in Section \ref{Risque&Optimality}, after the description of the main hypotheses. Section \ref{sec:adaptation} is devoted to adaptive estimation: the penalization strategy is described, and the oracle inequality as well as adaptive convergence rates are stated.  Numerical results illustrate the theoretical
properties in Section \ref{sec:simus}. We first calibrate our estimator and study its performances on simulated data in Section \ref{Sdata}. Then we apply our procedure on two real-data sets problems in Section \ref{sec:Rdata}: the prediction of the appliances electricity consumption of a day given the ones of the day before \citep{CFD17}, and the prediction of the evolution of the electricity prices from the wind power in-feed \citep{liebl2013modeling}. Finally, the proofs are gathered in Section \ref{sec:proofs}.

\section{Estimation method}
\label{Estimation}
 
\subsection{Notations}\label{sec:notations}

We introduce here some notations which will be used all along this document.  We denote by $\mathcal{L}_2(\mathbb H)$ the subspace of Hilbert-Schmidt operators on $\mathbb H$ equipped with its usual Hilbert-Schmidt norm defined for any operator $T\in\mathcal{L}_2(\mathbb H)$ as follows 
$$\|T\|_{\HS}=\left(\sum_{j=1}^\infty\|T \phi_j \|^2\right)^{1/2},$$
where $(\phi_j)_{j\geq 1}$ is a Hilbertian basis of $\mathbb H$. Note that  the Hilbert-Schmidt norm is independent of the Hilbertian basis choice. It is also worth mentioning that an integral operator is Hilbert-Schmidt if and only if the associated kernel is square integrable. This means that by assumptions, our target operator $S$ is Hilbert-Schmidt. We also need to define two operators that play a key role in the estimation procedure, namely the covariance and cross-covariance operators. To do so, we first define the tensor product between two elements $a$ and $b$ of $\mathbb H$ as
$$\begin{array}{lccl} b \otimes a : &\mathbb H &\longrightarrow &\mathbb H\\
                               &f  &\longmapsto     & \langle a,f\rangle b.\end{array}$$
The covariance operator of $X$, denoted $\Gamma$ is the operator defined by
$$ \begin{array}{cccl}\Gamma :&\mathbb H&\longrightarrow &\mathbb H\\
                         & \displaystyle      f  &\longmapsto     & \mathbb E[X\otimes X(f)]=\mathbb E[\langle X,f\rangle X].\end{array}.$$
Note that the covariance operator is a natural extension of the covariance matrix, in the infinite dimensional framework. The $(\lambda_j,\varphi_j)_j$ involved in \eqref{eq:KL} are the eigenelements of $\Gamma$. We also introduce the        cross-covariance operator $\Delta$ of $(X,Y)$ given by
$$ \begin{array}{cccl}\Delta :&\mathbb H&\longrightarrow &\mathbb H\\
                          &\displaystyle      f  &\longmapsto     & \mathbb E[Y\otimes X(f)]=\mathbb E[\langle X,f\rangle Y].\end{array}$$
Empirical counterparts of $\Gamma$ and $\Delta$, respectively denoted by $\Gamma_n$ and $\Delta_n$ will be useful in the definition of our estimators. These operators are naturally defined on $\mathbb{H}$ by
$$\Gamma_n = \frac{1}{n} \sum_{i=1}^n X_i \otimes X_i \quad \mbox{and} \quad \Delta_n = \frac{1}{n} \sum_{i=1}^n Y_i \otimes X_i.$$

In order to study the estimator behaviors, we use an optimality risk called the Mean Square Prediction Error (MSPE). This criterion is also used in \citet{crambes2013asymptotics,cardot2010thresholding, crambes2009smoothing} or \cite{brunel2016non}. The MSPE of a given estimator $\widehat S$ of $S$ is defined~as
$$\MSPE (\widehat S_n) = \mathbb{E} \Vert \widehat S_n(X_{n+1}) - S (X_{n+1}) \Vert^2,$$ where $X_{n+1}$ is a new observation of $X$, independent of $(X_i, \varepsilon_i), i = 1, \ldots, n$. This risk can also be written
\begin{equation}\label{eq:MSPE}
\MSPE(\widehat S_n)= \mathbb E\left[\Vert \widehat Y_{n+1}-\mathbb E\left[Y_{n+1}|X_{n+1}\right]\Vert^2\left| (X_i,Y_i)_{i=1,\ldots,n}\right.\right],
\end{equation}
where  $Y_{n+1}=SX_{n+1}+\varepsilon_{n+1}$, $\widehat Y_{n+1}=\widehat S_nX_{n+1}$, and $\mathbb{E}[\cdot|Z]$ is the conditional expectation given a variable $Z$. It is also linked with the Hilbert-Schmidt norm as follows, 
\begin{align}\label{eq:MSPE-HS}
\mathbb{E} \Vert \widehat S_n (X_{n+1}) - S (X_{n+1}) \Vert^2 &= \mathbb{E} \Vert (\widehat S_n - S) \Gamma^{1/2} \Vert_{\HS}^2,
\end{align}
see Lemma \ref{lemma-ps-Xn+1-HS} in the proof (Section \ref{sec:proof_bv}).

For two sequences $(a_j)_{j\geq 1}$ and $(b_j)_{j\geq 1}$ of real numbers, we write $a_j\asymp b_j$ if there exists $c\geq 1$ such that $c^{-1}a_j\leq b_j\leq ca_j$.

\subsection{Least-squares estimation}   
\label{LSE}

\subsubsection{Minimum contrast estimation}                      

The main goal of statistical estimation is to build an estimate that leads to a small risk.  Following the model selection device introduced by \cite{birge1998minimum}, we minimize an empirical counterpart of the risk, called the contrast function, over  finite dimensional subspaces of $\mathcal{L}_2(\mathbb H)$, to build projection-type estimators. Let $(\phi_j)_{j\geq 1}$ be an orthonormal basis of $\mathbb H=L^2([0,1])$. We introduce a collection of finite linear subspaces of $\mathcal{L}_2(\mathbb H)$, called the models and denoted $V_{m_1,m_2}$ for given $m_1,m_2$ in $\mathbb{N} \backslash \{0\}$. These models are defined as
$$V_{m_1,m_2}=\Span\{\phi_k\otimes \phi_j,\;1\leq j\leq m_1,\; 1\leq k\leq m_2\}.$$ 
Note that $V_{m_1,m_2}$ only contains integral operators. Subsequently, for any $ T\in \mathcal{L}_2(\mathbb H)$, let
$$\gamma_n(T) = \frac{1}{n}\sum_{i=1}^n\left\|Y_i-T(X_i)\right\|^2.$$
The operator $\gamma_n: \mathcal{L}_2(\mathbb H)\rightarrow\mathbb R$ is defined in the spirit of other regression contrast introduced for example by \cite{baraud2002model} and \cite{brunel2016non} and stands for an empirical version of the risk \eqref{eq:MSPE}. Thus, we set
\begin{equation}
\label{eq:estim_MC}
\widehat{S}_{m_1,m_2}\in\arg\min_{ T\in V_{m_1,m_2}}\gamma_n(T).
\end{equation}
To compute $\widehat{S}_{m_1,m_2}$, we introduce the matrices $A$ and $Y_{\phi}$ given by 
$$A=(\langle \Gamma_n \phi_j, \phi_k\rangle)_{j,k\in\{1,\ldots,m_1\}} \quad \mbox{and} \quad Y_{\phi}=(\langle \Delta_n \phi_j,\phi_k \rangle)_{j\in\{1,\ldots,m_1\}, k\in\{1,\ldots,m_2\}}.$$
The following result, which proof can be found in Section  \ref{sec:proof_prop_MC}, gives a condition for the existence of $\widehat{S}_{m_1,m_2}$.
\begin{proposition}\label{pro:estim_MC}
If the matrix $A$ is invertible, then $\widehat S_{m_1,m_2}$ in \eqref{eq:estim_MC} is uniquely defined, and 
$$\widehat{S}_{m_1,m_2}=\sum_{j=1}^{m_1}\sum_{k=1}^{m_2}\widehat b_{j,k}\phi_k\otimes\phi_j,$$
with $\widehat b=(\widehat b_{j,k})_{j\in\{1,\ldots,m_1\},k\in\{1,\ldots,m_2\}}$ defined by
$\widehat b = A^{-1}Y_{\phi}$.
\end{proposition}

\begin{remark} 
Since $S$ is an integral operator with kernel $\mathcal S$, we could also have defined a contrast function over the space of the  kernel functions: for any $\mathcal T\in L^2([0,1]^2)$, let  $$\gamma_n'(\mathcal T) = n^{-1}\sum_{i=1}^n\|Y_i-\int_{[0,1]}\mathcal T(s,\cdot)X(s)ds\|^2.$$ If we denote by $V'_{m_1,m_2}=\Span\{(t,s)\,\mapsto\,\phi_j(s)\phi_k(t),\;1\leq j\leq m_1,\;1\leq k\leq m_2\}$, we can set $\widehat{\mathcal S}_{m_1,m_2}\in\arg\min_{\mathcal T\in V'_{m_1,m_2}}\gamma_n'(\mathcal T)$. The estimator $\widehat{\mathcal S}_{m_1,m_2}$ is also uniquely defined under the assumptions of Proposition \ref{pro:estim_MC}, and for any $f\in \mathbb H$,
 \begin{equation}\label{eq:estim_MC_op}
\widehat{S}_{m_1,m_2}f=\int_{0}^1\widehat{\mathcal S}_{m_1,m_2}(s,\cdot)f(s)ds.
\end{equation}
\end{remark}
\begin{remark}
Defined with \eqref{eq:estim_MC}, $\widehat{S}_{m_1,m_2}$ estimates the orthogonal projection $\Pi_{m_1,m_2}^{op}S$
 of the operator $S$ onto the operator space $V_{m_1,m_2}$.  This projection operator can be written, for any $T\in \mathcal{L}_2(\mathbb H)$, 
 \begin{equation}\label{eq:egalite_proj}\Pi_{m_1,m_2}^{op}T=\Pi_{m_1}T\Pi_{m_2},
 \end{equation}
  where, for any $m\in\mathbb N\backslash\{0\}$, $\Pi_{m}$ is the projection operator on the subspace $\Span\{\phi_k,\; k\in\{1,\ldots,m\}\}$. The proof of \eqref{eq:egalite_proj} can be found in Section \ref{sec:proof_egalite_proj}.
\end{remark}

\subsubsection{Specific choice of the projection spaces : principal component basis}
\label{sec:estimPCA}

In the rest of this article, we focus on the basis of principal components. Recall that, by definition, the empirical covariance operator $\Gamma_{n}$ is self-adjoint. Moreover, since it is a finite-rank operator, it is also a compact operator. Then $\Gamma_n$ is diagonalisable in a Hilbertian basis, denoted by $(\widehat\varphi_j)_{j\geq 1}$. We also denote by $(\widehat\lambda_j)_{j\geq 1}$ its eigenelements, which are sorted in a decreasing order. The $(\widehat\varphi_j)_{j\geq 1}$ is called the empirical PCA basis of $X$. Notice that the operator $\Gamma_{n}$ is not invertible, since it has finite rank at most equal to  $n$. This means that the eigenvalues $(\widehat\lambda_j)_{j \geq 1}$ are zero from a given rank. Let us introduce its pseudo-inverse, $\Gamma_{n,m_1}^{\dag} $, defined for an index $m_1\in\mathbb N\backslash\{0\}$ by
$$\Gamma_{n,m_1}^{\dag} =\sum_{j=1}^{m_1}\frac{1}{\widehat\lambda_j}\widehat\varphi_j\otimes\widehat\varphi_j,$$
for $m_1 \leq m_{\max}$, with $m_{\max}=\max_{m\geq 1}\{\widehat\lambda_m>0\}$ is the rank from which the eigenvalues are equal to zero, and $\Gamma_{n,m_1}^{\dag}=\Gamma_{n,m_{\max}}^{\dag}$ for $m_1>m_{\max}$.

\medskip

We obtain the following expressions for the least-squares estimators of the linear operator $S$ and its kernel $\mathcal S$ . The proof can be found in Section \ref{sec:proof_prop_MC_ACP}.

\begin{proposition}\label{pro:estim_MC_ACP} On the PCA basis, the least-squares estimator for the kernel $\mathcal S$ exists, and is uniquely defined by 
\begin{equation}\label{eq:MC_ACP_kernel}
\widehat{\mathcal S}_{m_1,m_2}(s,t)=\sum_{j=1}^{m_1}\sum_{k=1}^{m_2}\frac{1}{\widehat\lambda_j}\langle\Delta_n\widehat\varphi_j,\widehat\varphi_k\rangle\widehat\varphi_j(s)\widehat\varphi_k(t), \; (s,t)\in [0,1]^2.
\end{equation}
Moreover, the expression of the resulting estimator for the linear operator $S$ is 
\begin{equation}\label{eq:MC_ACP_op}
\widehat S_{m_1,m_2}=\Delta_n\Gamma_{n,m_1}^{\dag} +\sum_{\substack{1\leq j\leq m_1  \\ 1\leq k\leq m_2\\ j\neq k}}\frac{1}{\widehat\lambda_j}\langle\Delta_n\widehat\varphi_j,\widehat\varphi_k\rangle\widehat\varphi_j\otimes\widehat\varphi_k.
\end{equation}
\end{proposition}

Remark that our estimator can also be written $\widehat S_{m_1,m_2}= \widehat{\Pi}_{m_2}\Delta_n\Gamma_{n,m_1}^{\dag} $, where  $\widehat\Pi_{m_2}$ is the projection operator onto   the finite dimensional subspace $\Span\{\widehat\varphi_k, k=1,\ldots,m_2\}$.  Thus, $\widehat S_{m_1,m_2}$ can be compared to the estimator of \cite{crambes2013asymptotics} which writes $\widehat{S}_{m_1}^{CM}=\Delta_n\Gamma_{n,m_1}^{\dag} $. Our choice is based on the fact that the initial regression problem comes down to estimate the kernel $\mathcal{S} \in L^2 ([0,1]^2)$ of the operator $S$,  which brings out two projection dimensions. 
Our estimate is thus the same as the  estimator with \enquote{double truncation} of \cite{imaizumi2018pca} (see their equation (7) p.19), even if they do not introduce it as a minimum of contrast estimator. The definition of  $\widehat S_{m_1,m_2}$ as an operator that minimizes a contrast allows us to derive non-asymptotic upper-bounds for the prediction error, and to propose a data-driven way to select the best projection dimensions.  

\section{Upper and lower bounds of the estimation risk }
\label{Risque&Optimality}

In this section, we provide sharp upper bounds for the estimation risk of the estimator $\widehat S_{m_1,m_2}$, for any but fixed $(m_1,m_2)\in(\mathbb N\backslash\{0\})^2$, after stating the main hypotheses. We also establish a lower bound for the prediction risk, to ensure that the collection of estimates is reasonable.

\subsection{Assumptions}
\label{Assumptions}

Classically, we need to make some assumptions for optimal theoretical results. We distinguish different types of assumptions: regularity assumptions on $S\Gamma^{1/2}$, regularity assumptions on $X$, moment assumptions on $X$ and moment assumptions on the noise $\varepsilon$.  

\medskip
\textit{Assumption on $S\Gamma^{1/2}$:} We consider by analogy with \cite{brunel2016non}, the regularity space of type ellipsoid, defined for all positive $\alpha, \beta, R$ by  
\begin{equation*}
\mathcal{W}_{\alpha, \beta}^R = \left\{ T \in \mathcal{L}_2(\mathbb H) , \sum_{j = 1}^{+\infty} \sum_{r = 1}^{+\infty} \eta_{\alpha} (j) \psi_{\beta} (r) \langle T (\varphi_j) , \varphi_r \rangle^2 \leq R^2 \right\}, 
\end{equation*}
where $\alpha, \beta > 0$ and for all $\gamma > 0$, the functions $\eta_{\gamma}$ is defined such that
$$\eta_{\gamma}(j)\asymp j^{\gamma}\mbox{ or }\eta_{\gamma}(j)\asymp \exp(j^{\gamma}),$$
and the same for $\psi_{\gamma}$. We speak about the \enquote{polynomial case} or the \enquote{exponential case} in the sequel. These regularity spaces are generalization of the ellipsoid sets in the finite dimensional framework. Moreover, the regularity parameters $\alpha$ and $\beta$ are respectively the convergence rates towards 0 of the operator components in both directions.  

\smallskip

$\boldsymbol{\mathcal{A}_{1}}:$ We assume that $S \Gamma^{1/2}$ belongs to $\mathcal{W}_{\alpha, \beta}^R$ for some positive regularity parameters $\alpha, \beta, R$. 

\smallskip

Asssumption $\boldsymbol{\mathcal{A}_{1}}$ is a smoothness assumption on the operator $S$ we want to recover. In nonparametric estimation, it is classical, and permits to control the bias term of the risk and to derive convergence rates (see \textit{e.g.} \cite{tsybakov2008introduction}). The kind of smoothness ball (ellipsoid space) we choose is also classical for projection type estimation (again, see \cite{tsybakov2008introduction}, but also  \cite{barron1999risk}, or \cite{brunel2016non}).

The specificity of our assumption  is that it is a joint regularity assumption both on  $S$ and on the covariate $X$. This technical choice is related to the choice of the mean squared prediction error we consider in this work: this risk is linked to the Hilbert-Schmidt norm of $S\Gamma^{1/2}$, see \eqref{eq:MSPE-HS}, it is thus natural that the smoothness assumption refers to this operator.  This was also the case in the paper of \cite{crambes2013asymptotics} and was pointed out by \cite{hilgert2013minimax}. In a similar but different way same discussions appear in \cite{comte2012adaptive}: the link between the smoothness of $S$ and $X$ appears in a \enquote{compatibility} assumption on the norms, see their section 3.1. If we replace the prediction risk with a quadratic risk, the most appropriate choice would be to impose a regularity assumption only on $S$, as done in \cite{imaizumi2018pca}.

\medskip

\textit{Assumptions on $X$:} Requiring a regularity on $X$ is tantamount to making assumptions on $\Gamma$ and its eigenvalues $(\lambda_j)_{j \geq 1}$. In particular, we consider that the eigenvalues are all distincts. 
In addition, we assume the following, 

\smallskip
$\boldsymbol{\mathcal{A}_{2}}:$ For all $j$ in $\mathbb{N}\backslash\{0\}$, we have $\lambda_j \psi_{\beta} (j) \geq 1$. 

\smallskip
The previous assumption ensures a separability condition on the eigenvalues of $\Gamma$. Indeed, considering that the model estimation is based on the estimation of the eigenfunctions of $\Gamma_n$, a separability condition on the eigenvalues of $\Gamma$, ensuring that they are not too close to each others is needed. Separation conditions on the eigenvalues of the covariance operator are usual in functional PCA regression. A usual alternative is to consider assumptions on the gap between two consecutive eigenvalues, as in \cite{imaizumi2018pca} or \citet{MR2332269}. 

\smallskip
$\boldsymbol{\mathcal{A}_{3}}:$ There exists a convex positive function $x \mapsto \lambda(x)$ such that, for all $j$ in $\mathbb{N}\backslash\{0\}$ : $\lambda_j = \lambda(j)$.

\smallskip
$\boldsymbol{\mathcal{A}_{4}}:$ There exists a constant $\gamma >0$ for which the sequence $\left( j \lambda_j \max \{\ln^{1+\gamma}(j),1 \} \right)_{j \geq 1}$ decreases.

\smallskip
Assumptions $\boldsymbol{\mathcal{A}_{3}}$ and $\boldsymbol{\mathcal{A}_{4}}$ permits to obtain some decreasing rate of convergence of the eigenvalues of $\Gamma$ and are classical in obtaining the optimal convergence rate of estimation. Similar assumptions have also been made in \cite{crambes2013asymptotics} and \cite{brunel2016non}.  

\medskip
\textit{Moment assumptions on $X$:}  

\smallskip 
$\boldsymbol{\mathcal{A}_{5}}:$ There exists a constant $b > 0$ such that, for all $l$ in $\mathbb{N}\backslash\{0\}$, 
\begin{equation*}
\sup_{j \geq 1} \mathbb{E} \left[ \frac{\langle X , \varphi_j \rangle^{2l}}{\lambda_j^l} \right] \leq l! b^{l-1}.
\end{equation*}

\smallskip

Like Assumption $\boldsymbol{\mathcal{A}_{3}}$, one can also find Assumption $\boldsymbol{\mathcal{A}_{5}}$ in \cite{crambes2013asymptotics} and \cite{brunel2016non}. The control of the moment of the random variables $\langle X , \varphi_j \rangle$, $j\geq 1$ is required to apply Bernstein's  exponential inequality. 
  
\smallskip
$\boldsymbol{\mathcal{A}_{6}}:$ For all $j \neq k$, $\langle X, \varphi_j \rangle$ and $\langle X, \varphi_k \rangle$ are independent.

\smallskip

Both assumptions $\boldsymbol{\mathcal{A}_{5}}$ and $\boldsymbol{\mathcal{A}_{6}}$ are satisfied when $X$ is a Gaussian process (see \citealt[Section 1.4]{ash_topics_1975}). For general (non Gaussian) processes, we know that $\langle X, \varphi_j \rangle$ and $\langle X, \varphi_k \rangle$ are, at least, uncorrelated since $\mathbb E[\langle X, \varphi_j \rangle\langle X, \varphi_k \rangle]=\langle \Gamma\varphi_j,\varphi_k\rangle=\lambda_j\mathbf 1_{\{j=k\}}$.
\medskip 

\textit{Moment assumption on $\varepsilon$:}

\smallskip
$\boldsymbol{\mathcal{A}_{7}}:$ There exists $p > 6$ such that $\tau_p = \mathbb{E} \Vert \varepsilon \Vert^p < +\infty$.

\smallskip 
The moment assumption $\boldsymbol{\mathcal{A}_{7}}$ is classically needed to obtain adaptive estimators of the model $S$. Indeed, when no assumption is imposed on the fluctuation of noise, it is not possible to construct optimal estimator without knowing the regularity of the model. We can deduce easily, e.g. from \citet[Lemma 8]{CR14}, that it is verified if $\mathbb E[\|\varepsilon\|^2]<+\infty$ and $\mathbb E[\langle \varepsilon,f\rangle^p]<+\infty$, for all $f\in\mathbb H$. As particular cases, $\boldsymbol{\mathcal{A}_{7}}$ is verified as soon as the noise is Gaussian, or bounded.


\subsection{Upper bound of the Mean Square Prediction Error (MSPE)}
%
Theorem \ref{thm-upper-bound-MSPE} below gives a first sharp upper bound of the Mean Squared Prediction Error of the estimator $\widehat S_{m_1,m_2}$ with respect to the projection dimensions $m_1$ and $m_2$. This permits to have (up to a positive constant) an order of magnitude of the prediction error for each theoretical choice of $m_1$ and $m_2$. The proof of the results of the section can be found in Section \ref{sec:proof_bv}. 

\begin{theorem}
\label{thm-upper-bound-MSPE}
Under Assumptions $\boldsymbol{\mathcal{A}_{1}}$ to $\boldsymbol{\mathcal{A}_{6}}$, the Mean Squared Prediction Error of the estimator $\widehat S_{m_1,m_2}$ is upper bounded by
\begin{equation}
\label{eq-thm-upper-bound-MSPE}
\MSPE (\widehat S_{m_1,m_2}) \leq \sigma_{\varepsilon}^2 \frac{m_1}{n} + 3 \sum_{j = m_1 + 1}^{+\infty} \Vert S \Gamma^{1/2} (\varphi_j) \Vert^2 + 3 \sum_{j = 1}^{m_1} \Vert (\Id - \Pi_{m_2}) S \Gamma^{1/2} (\varphi_j) \Vert^2 + A_{n,m_1} + B_{n,m_1} + D_{n,m_2} + E_n,
\end{equation}
where, for a constant $C$ which does not depend on $n$, $m_1$, $m_2$,
\begin{align*}
A_{n,m_1} &= \sigma_{\varepsilon}^2 \frac{C m_1^2 \ln^2 (m_1)}{n^2}, \quad B_{n,m_1} = \frac{C m_1^2 \lambda_{m_1} \Vert S \Vert_{\HS}}{n}, \quad E_n = \frac{C}{n^2} \Vert S \Gamma^{1/2} \Vert_{\HS}^2 \\
D_{n,m_2} &= \frac{C m_2 \lambda_{m_2}}{n} + \frac{C m_2^2 \ln(m_2)}{n \psi_{\beta} (m_2)} + \frac{C m_2^3}{n \psi_{\beta} (\lfloor m_2/2 \rfloor)} + \frac{C \ln^4 (n)}{n^2} \left( \sum_{k = 1}^{m_2} \frac{k^2 \ln^2 (k)}{\sqrt{\psi_{\beta} (k)}} \right)^2.
\end{align*}
\end{theorem}

In Theorem \ref{thm-upper-bound-MSPE} appears a bias-variance trade-off. The first term in the right side of Equation \eqref{eq-thm-upper-bound-MSPE} is a variance term, which increases with $m_1$. The two following terms are bias terms: one is decreasing with $m_1$, the other one with $m_2$. Both are related to the smoothness of $S \Gamma^{1/2}$. Notice right now that it is not the smoothness of the target function $S$ that influences the result, but the one of $S \Gamma^{1/2}$. This is consistent with the choice of the risk, since the prediction error we study can also be written
$$\MSPE(\widehat S_{m_1,m_2})=\mathbb E\Vert (\widehat {S}_{m_1,m_2} - S) \Gamma^{1/2} \Vert_{\HS}^2,$$
see  Lemma \ref{lemma-ps-Xn+1-HS} below. The same phenomenom occurs for \cite{crambes2013asymptotics}. Compared to their result for the estimator $\widehat S_{m_1}^{CM}$ (see Theorem 2 p.2633 in \citealt{crambes2013asymptotics}), the first two terms of the bias-variance decomposition \eqref{eq-thm-upper-bound-MSPE} are the same, but we have an additional bias term (third term in the right-hand-side of \eqref{eq-thm-upper-bound-MSPE}), which depends on the index $m_2$. We prove in Corollary~\ref{corr-upper-bound-MSPE} below that the other terms are negligible.

\begin{corollary}
\label{corr-upper-bound-MSPE}
Assume that we are in the case where the function $\psi_{\beta}$ is polynomial with $\beta > 6$ or exponential. Assume also that there exists $\nu>0$ such that $\lambda_j\leq j^{-1-\nu}$, for any $j\geq 1$. Under the assumptions of Theorem \ref{thm-upper-bound-MSPE}, we have the following bound of the non-asymptotic  maximal prediction risk of $\widehat S_{m_1,m_2}$.
\begin{equation}\label{eq:compro_bv}
 \inf_{\substack{m_1, m_2 \in  \mathbb{N}\backslash\{0\}\\  m_1\leq n/\ln^2(n)}} \sup_{S \Gamma^{1/2} \in \mathcal{W}_{\alpha, \beta}^R} \MSPE (\widehat S_{m_1,m_2}) \leq \inf_{\substack{m_1 \in  \mathbb{N}\backslash\{0\}\\  m_1\leq n/\ln^2(n)}} \left\{ \sigma_{\varepsilon}^2 \frac{m_1}{n} +  \frac{3}{\eta_{\alpha} (m_1)} \right\} + \frac{c}{n},
\end{equation}
where  $c$ is a positive constant.
\end{corollary}

Some comments are needed at this point. The dimension parameter $m_2$ does not appear in the upper-bound. A similar phenomenon has been observed by \cite{imaizumi2018pca}. It is mainly due to the fact that the variance is completely independent of it. Hence, since the bias decreases to 0 when $m_2\to+\infty$, it is sufficient to choose $m_2$ sufficiently large so that the bias is negligible (remind that the estimator is well defined even in the case $m_2=+\infty$). 

Notice also that the additional assumption $\lambda_j\leq j^{-1-\nu}$ is very mild and useful only for technical purpose (we recall that, since the operator $\Gamma$ is trace-class, $\sum_{j\geq 1}\lambda_j<+\infty$). It is satisfied if the eigenvalues decrease at a polynomial or exponential rate. It can also be relaxed to allow us to choose $\lambda_j=(j\ln^{\mu}(j))^{-1}$ for some $\mu\geq 1$.

Corollary \ref{corr-upper-bound-MSPE} gives the sharpest possible upper-bound of the prediction risk for the estimator we define by projection onto the basis of principal components.  In the next section, we show that the upper-bound of Corollary~\ref{corr-upper-bound-MSPE} is optimal over the ellipsoidal regularity spaces we consider here.


%
%

\subsection{Lower bound of the minimax Mean Square Prediction risk}\label{sec:lower-bound}

In this section, we demonstrate that the upper-bound of the Mean Square Prediction risk obtained in Corollary \ref{corr-upper-bound-MSPE} is optimal in the minimax sense in a non-asymptotic framework. This result is stated in Theorem \ref{thm-lower-bound-MSPE} below. The demonstration of this result is based on a reduction scheme to a finite number of hypotheses, as explained in \cite{tsybakov2008introduction}. We apply the Kullback-Leibler version of Assouad's Lemma, and  the Cameron-Martin theorem \citep{lifshits2012lectures}. It permits to control the likelihood expectation between different possible data distributions in the finite model collection. 

\begin{theorem}
\label{thm-lower-bound-MSPE}
Let $\alpha>0$, $\beta>0$ and $R>0$, we have the following lower bound, for a constant $C>0$,
\begin{equation*}
\inf_{\widehat S_n} \sup_{S \Gamma^{1/2} \in \mathcal{W}_{\alpha, \beta}^R} \MSPE (\widehat S_{n}) \geq C \inf_{m_1 \in \mathbb{N}\backslash\{0\}} \left\{ \sigma_{\varepsilon}^2 \frac{m_1}{n} +  \frac{3}{\eta_{\alpha} (m_1)} \right\},
\end{equation*}
where the infimum is taken over all estimators $\widehat S_n$ calculated from a sample $\{(X_i,Y_i), i=1,\hdots,n\}$ following model~\eqref{eq:model}, under the assumption that the noise $\varepsilon$ is a Gaussian process.
\end{theorem}

This lower bound permits to derive the minimax explicit convergence rates in the  polynomial and exponential cases. 

\begin{corollary}\label{coro:minimax_rates}
Under the assumptions of Theorem \ref{thm-lower-bound-MSPE}, we compute the two following convergence decay for the minimax estimation risk, up to a constant $C>0$.
\begin{enumerate}
\item If $\eta_{\alpha}(j)\asymp j^\alpha$ (polynomial case) then,
\begin{equation*}
 \inf_{\widehat S_n}  \sup_{S \Gamma^{1/2} \in \mathcal{W}_{\alpha, \beta}^R} \MSPE (\widehat S_{n}) \geq Cn^{-\frac{\alpha}{\alpha+1}}.
\end{equation*}

\item If $\eta_{\alpha}(j)\asymp\exp(j^{\alpha})$ (exponential case) then,
\begin{equation*}
 \inf_{\widehat S_n}\sup_{S \Gamma^{1/2} \in \mathcal{W}_{\alpha, \beta}^R}\MSPE (\widehat S_{n})\geq C\frac{(\ln(n))^{1/{\alpha}}}{n}.
\end{equation*}
\end{enumerate}

\end{corollary}

From Corollaries~\ref{corr-upper-bound-MSPE} and~\ref{coro:minimax_rates}, we deduce that the projection estimators onto the PCA bases achieve the minimax rate for a suitable choice of the dimension $m_1$ and $m_2=+\infty$. For example, it can be deduced in  the polynomial case, that the optimal sharp upper-bound in Corollary~\ref{corr-upper-bound-MSPE}  is obtained for $m_1 = C n^{1/(1+\alpha)}$ and $m_2 \rightarrow \infty$, where $C$ is a universal positive constant.  This leads to an upper-bound  of order $n^{-\alpha/(\alpha+1)}$, which effectively matches with the lower bound of Corollary~\ref{coro:minimax_rates}. The latter estimation rate is known to be optimal in many other nonparametric estimation problems, see for example \cite{tsybakov2008introduction}. It can  also be remarked that, in both polynomial and exponential cases, the rates we get are very similar to the minimax rates obtained by \citet[Theorem 4]{brunel2016non} or \citet[Proposition 3.1]{cardot2010thresholding} in the functional linear model with scalar outputs.  

Although minimax optimal if the projection dimensions $m_1$ and $m_2$ are well choosen, the projection estimates  are not adaptive at this stage. Indeed, the optimal dimension $m_1$ depends on the regularity $\alpha$ of the operator $S\Gamma^{1/2}$, which is generally unknown. In the next section, we focus on the construction of an adaptive estimator of the model, meaning that is does not imply any knowledge of the unknown model regularity and achieves the optimal required estimation rate.

\section{Adaptive estimation}\label{sec:adaptation}

\subsection{Model selection}
The objective is to perform adaptive model selection, which does not depend on the unknown smoothness of the model $S$, but only on the available data.  As a reminder, for given projection dimensions $m_1$ and $m_2$, we estimate the operator $S$ by $\widehat{S}_{m_1, m_2} = \widehat{\Pi}_{m_2} \Delta_n \Gamma_{n,m_1}^{\dagger} $, where the operators $\widehat{\Pi}_{m_2}$, $\Delta_n$ and $\Gamma_{n,m_1}^{\dagger} $ are defined in Section \ref{LSE}.

The idea is to propose a procedure which automatically selects the best projection dimensions $m_1$ and $m_2$, that is the best estimator in the collection $(\widehat S_{m_1,m_2})_{m_1,m_2}$.  According to the result of Corollary \ref{corr-upper-bound-MSPE}, we choose $m_2 \rightarrow +\infty$ and we select $m_1$ in a collection $\mathcal{M}_n = \{1, \ldots, N_n \}$, where the size of the collection $N_n$ statisfied $N_n\leq \lfloor n/\ln^2(n)\rfloor$ where $\lfloor \cdot \rfloor$ is the floor function, associating to each $x$ in $\mathbb{R}$ the largest integer less or equal to $x$. Thus, the issue we consider now is the choice of an estimator in the collection $(\widehat S_{m_1,\infty})_{m_1\in \mathcal M_n}$, where $\widehat{S}_{m_1, \infty} = \Delta_n \Gamma_{n,m_1}^{\dagger} $ corresponds in fact to the estimator of \cite{crambes2013asymptotics}. The method we use is derived from the model selection tools developed by \cite{barron1999risk}, as in \cite{brunel2016non} or  \cite{comte2012adaptive}. A clear and detailed account is given in \cite{Mas2007}. We want to select the ''best'' estimator in the collection $(\widehat S_{m_1,\infty})_{m_1\in \mathcal M_n}$, that is the one which has the smaller risk. Since the risk is unknown in practice, the oracle $m_1^*=\arg\min_{m_1\in\mathcal M_n} \MSPE (\widehat S_{m_1,\infty})$ is also unknown, and the risk $\MSPE(\widehat S_{m_1,\infty})$ should be replaced by an empirical counterpart. Since the contrast function is an empirical version of the risk, the first idea is to choose $\arg\min_{m_1\in\mathcal M_n}\gamma_n (\widehat S_{m_1,\infty})$. However, since the contrast function decreases when $m_1$ grows, the choice of $\arg\min_{m_1\in\mathcal M_n}\gamma_n (\widehat S_{m_1,\infty})$ will lead to the selection of the largest index in the collection  $\mathcal{M}_n$. One of the main idea of model selection theory is to introduce a penalty to balance this decrease, usually of the order of the variance. The dimension parameter $m_1$ is choosen as the one which minimizes a penalized contrast function, 
\begin{equation}\label{eq:selecm1}
\widehat{m}_1 = \argmin_{m_1 \in \mathcal{M}_n} \left( \gamma_n (\widehat{S}_{m_1,\infty} ) + \pen(m_1) \right),
\end{equation}
where $\gamma_n$ is defined in Section \ref{LSE} and $ \pen$ is the penalty function defined as $\displaystyle \pen : m_1 \mapsto 8 (1 + \delta) \sigma_{\varepsilon}^2 m_1/n$, with $\delta > 0$ a numerical constant that will be tuned in practice, see Section \ref{sec:simus}. 

Remark that, when $m_1$ is fixed, $\gamma_n (\widehat{S}_{m_1,m_2} )$ decreases with $m_2$ by definition and $\pen(m_1)$ does not depend on $m_2$. Thus,  $(\widehat{m}_1,+\infty)$ is also a solution of the minimization problem
\[
\min_{(m_1,m_2)\in\mathcal{M}_n\times\mathbb N\backslash\{0\}\cup\{+\infty\}} \left( \gamma_n (\widehat{S}_{m_1,m_2} ) + \pen(m_1) \right). 
\]
With this writing, the selection procedure has strong similarities with the usual model selection procedures when two dimensions have to be selected (see e.g. \citealt{plancade2013adaptive,Lacour07}). Here, the specificity  is that the penalty criterion does not depend on $m_2$ (since the variance term only  depends on $m_1$). This makes it possible to consider, in an equivalent way, the criterion~\eqref{eq:selecm1} we have defined, which focuses on $m_1$ only. 

\subsection{Oracle-type inequality}
Theorem \ref{thm-UB-Emp-Risk} proves that the  penalty term introduced above has the good order of magnitude to automatically realize the best bias-variance trade-off. In the statement of the result, and in the sequel, $\Vert \cdot \Vert_n$ is the empirical norm defined for all operator $T$ as $\Vert T \Vert_n^2 = 1/n \sum_{i = 1}^n \Vert T(X_i) \Vert^2$ and $\widehat\Pi_{m_1, \infty}^{op}$ is the orthogonal projection onto the closure of $V_{m_1,\infty}=\Span\{\widehat\varphi_k\otimes\widehat\varphi_j,\;1\leq j\leq m_1,\; m_2 \geq 1 \}$. The proof can be found in Section \ref{sec:proof_sec_adapt}.

\begin{theorem}
\label{thm-UB-Emp-Risk}
Under Assumption $\mathcal{A}_{7}$, we have the following upper bounding, for all $\zeta>0$, 
$$\mathbb{E} \Vert S - \widehat{S}_{\widehat{m}_1, \infty} \Vert_n^2 \leq (1+\zeta) \inf_{m_1 \in \mathcal{M}_n} \left\{ \mathbb{E} \Vert S - \widehat{\Pi}^{op}_{m_1, \infty} S \Vert_n^2 + c(\zeta)\pen(m_1) \right\} + \frac{C'}{n}, $$
for a constant $C'>0$ which does not depend neither on $n$, nor on $m_1$ and $c(\zeta)=(2+\zeta)/(1+\zeta)$.
\end{theorem}

Theorem \ref{thm-UB-Emp-Risk} proves that the selected estimator achieves the best bias-variance compromise, up to a multiplicative constant, and the addition of the term $C'/n$, which is negligible.  Then it achieves the minimax rate and, since the dimension selection procedure does not require the knowledge of the unknown regularity $\alpha$, it is adaptive.  A similar result could be obtained for the risk MSPE, but at the price of additional technicalities. Indeed, to obtain such result it is necessary to prove that, with sufficiently large probability, the quantity $\|S\|^2_n/\MSPE(S)$ is lower bounded by a constant, for all $S\in V_{m_1,m_2}$ which is a random space (depending on the data $X_1,\hdots,X_n$). We could draw inspiration e.g. from the proof of \citet[Lemma 6]{brunel2016non}.

\section{Numerical study}\label{sec:simus}

The aim of this section is to assess the performance of the adaptive estimation method presented in Section \ref{sec:adaptation}. In Section \ref{Sdata}, we perform simulation studies for various functional models. Subsequently, we apply  the estimation method on two real data cases, in Section \ref{sec:Rdata}. All the study has been carried out with the free software \RR{}.

\subsection{Simulation study}
\label{Sdata}

\subsubsection{Simulated data}
To implement our estimation method, we consider three data generating mechanisms \ref{M1}, \ref{M2} and \ref{M3}. Each model is defined by the equation 
\begin{equation}\label{eq:simul_data}
Y^{(\ell)} = \int_0^1 \mathcal{S}_{\ell} (s, \cdot ) X^{(\ell)}(s) \, \mathrm{d} s + \varepsilon^{(\ell)},
\end{equation}
where $\ell = 1,2,3$. We also denote $S^{(\ell)}$ the integral operator with kernel $\mathcal{S}_{\ell}$. The analytical expressions of the kernels and noises are given below.  
\begin{enumerate}[label=(\roman*)]
\item \label{M1} The kernel is defined as $\mathcal{S}_1 : (s,t) \mapsto s^2 + t^2$ and the noise $\varepsilon^{(1)}$ is generated according to a standard Brownian motion divided by 20. In addition, the Karhunen–Lo\`eve decomposition of the covariate $X^{(1)}$ is written as $X^{(1)} = \sum_{j = 1}^{k_0} \sqrt{\lambda_j} \xi_j \varphi_j$, where $k_0 = 8$, $\lambda_j = 1/(\pi^2 (j- 0.5)^2)$,  $\varphi_j : t \mapsto \sqrt{2} \sin \left( (j - 0.5) \pi t \right)$, $j=1,\ldots,k_0$, and $(\xi_j)_j$ are independent standard Gaussian random variables.
\item \label{M2} The implementation of this simulation case is the same as \ref{M1} with only one difference, the error $\varepsilon^{(2)}$ is a Brownian motion divided by 2.
\item \label{M3} The model kernel is given by the equation $\mathcal{S}_3 : (s,t) \mapsto \sum_{j,l = 1}^{k_1} b_{j,l} \varphi_l (s)  \varphi_j (t)$, where $k_1 = 50$, for all $j,l$ in $\mathbb{N}\backslash\{0\}$, $\varphi_j : u \mapsto \sqrt{2} \cos(j \pi u)$ and $b_{j,l} = 4 (-1)^{j + l} j^{-\gamma} l^{-\beta}$, with $\beta = 3$ and $\gamma = 2.5$. The input is the random function $X^{(3)} = \sum_{j = 1}^{k_1} j^{-\alpha/2} U_j \varphi_j$, where $\alpha = 1.2$ and $U_j$ are independent uniform distributions over $[-\sqrt{3}, \sqrt{3}]$. Finally, the noise is defined as $\varepsilon^{(3)} = \sum_{j = 1}^{k_1} j^{-\delta/2} \xi_j \varphi_j$, with $\delta = 1.1$ and $(\xi_j)_{j\geq 1}$ are independent standard Gaussian random variables.
\end{enumerate}

The simulation cases \ref{M1} and \ref{M2} are drawn from \cite{crambes2013asymptotics}, while the model \ref{M3} is studied in \cite{imaizumi2018pca} with slight modifications.

\subsubsection{Implementation of the method}

To perform the model selection strategy for the examples described above, a first step is to compute the penalty term of the procedure in \eqref{eq:selecm1}. For simulation purposes, we keep the true value of $\sigma_{\varepsilon}$ (it will be replaced by an empirical counterpart for real data analysis in Section \ref{sec:Rdata}). But we have to wisely choose the calibration parameter $\kappa = 8 (1 + \delta)$, according to Section \ref{sec:adaptation} notations. Unlike the theoretical framework and for practical reasons, the chosen values of $\kappa$ are not necessarily greater than 8, as it is usual in model selection. More precisely, we compare the choices of $\kappa$ values in the range $[0.2, 2]$, with a step of 0.2 between each two successive values. For every $\kappa$ value and for each model, we generate $N = 500$ independent samples of inputs/outputs $(X_i,Y_i)_{i\in\{1,\ldots,n\}}$ of size $n = 600$ each. We estimate the three models $N$ times and for each iteration, we measure the prediction error by generating a new observation of the input/output pair. Thus, for any value of $\kappa$, and any $\ell=1,2,3$, we simulate $(X_{i,k}^{(\ell)},Y_{i,k}^{(\ell)})$ for $i=1,\ldots,n$, $k=1,\ldots,N$ from Model \eqref{eq:simul_data}, we compute the Empirical Mean Square Prediction error 
$$\EMSPE^{(\ell,\kappa)}=\frac{1}{N}\sum_{k=1}^N\left\|\widehat{S}_{\widehat m_1,\infty}^{(\ell,k)}(X_{n+1,k}^{(\ell)})-S^{(\ell)}(X_{n+1,k}^{(\ell)})\right\|^2$$
where $\widehat{S}_{\widehat m_1,\infty}^{(\ell,k)}$ is the penalized contrast estimator computed from $(X_{i,k}^{(\ell)},Y_{i,k}^{(\ell)})_{i=1,\ldots,n}$ and $X_{n+1,k}^{(\ell)}$ is distributed like $X^{(\ell)}$ and independent of the $(X_{i,k}^{(\ell)})_i$.

\begin{figure}[!h]
\centering
\includegraphics[scale=0.5]{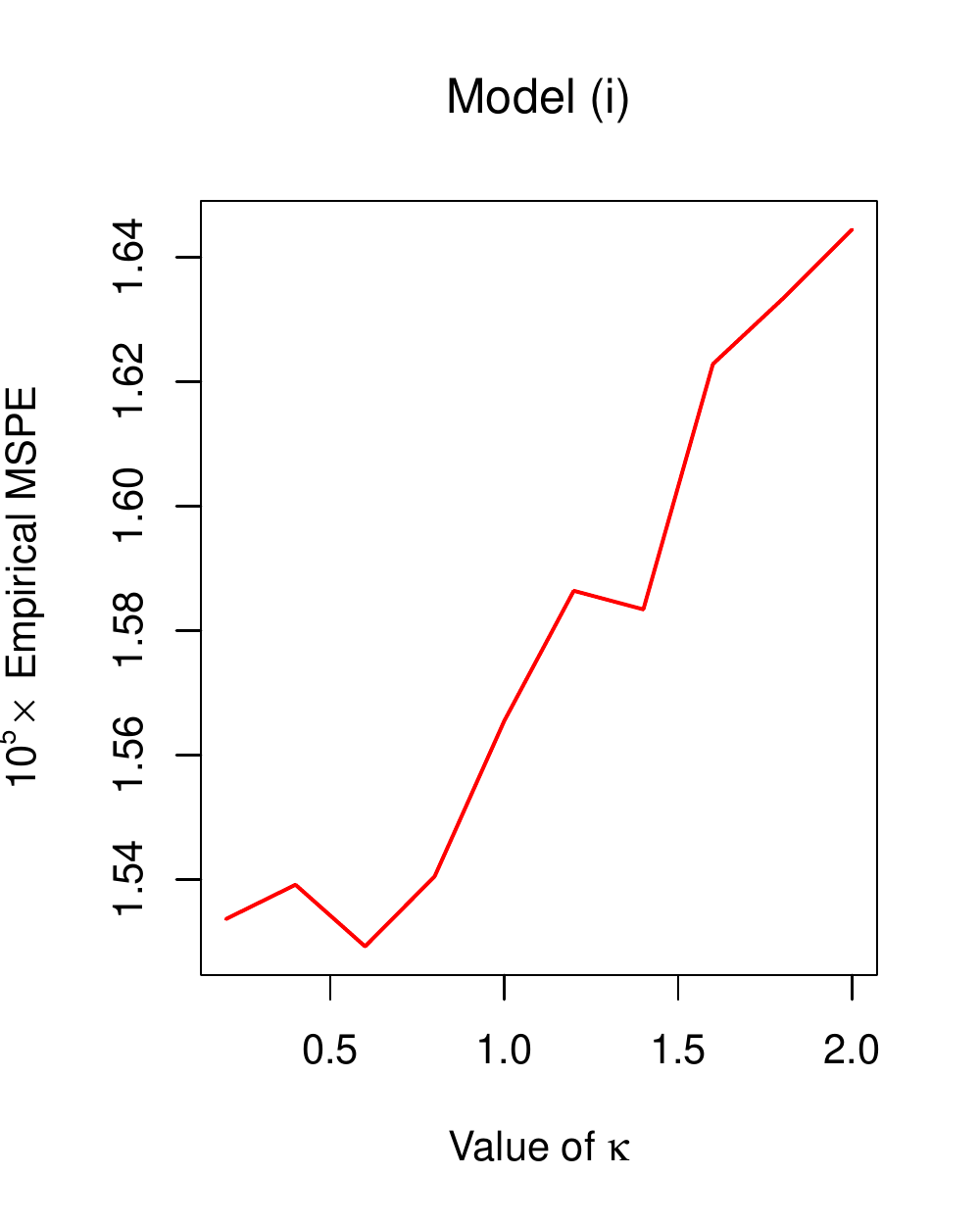} 
\includegraphics[scale=0.5]{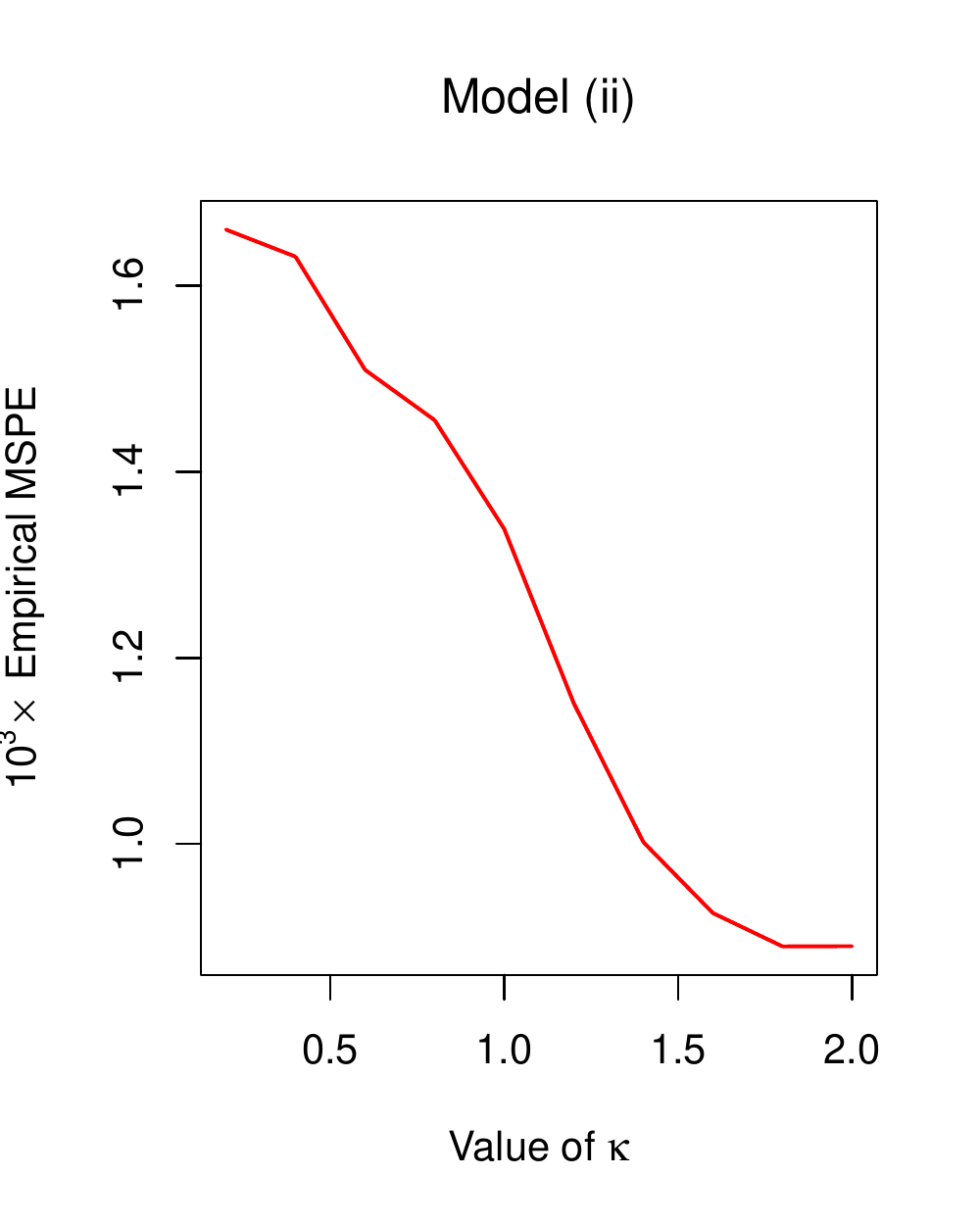} 
\includegraphics[scale=0.5]{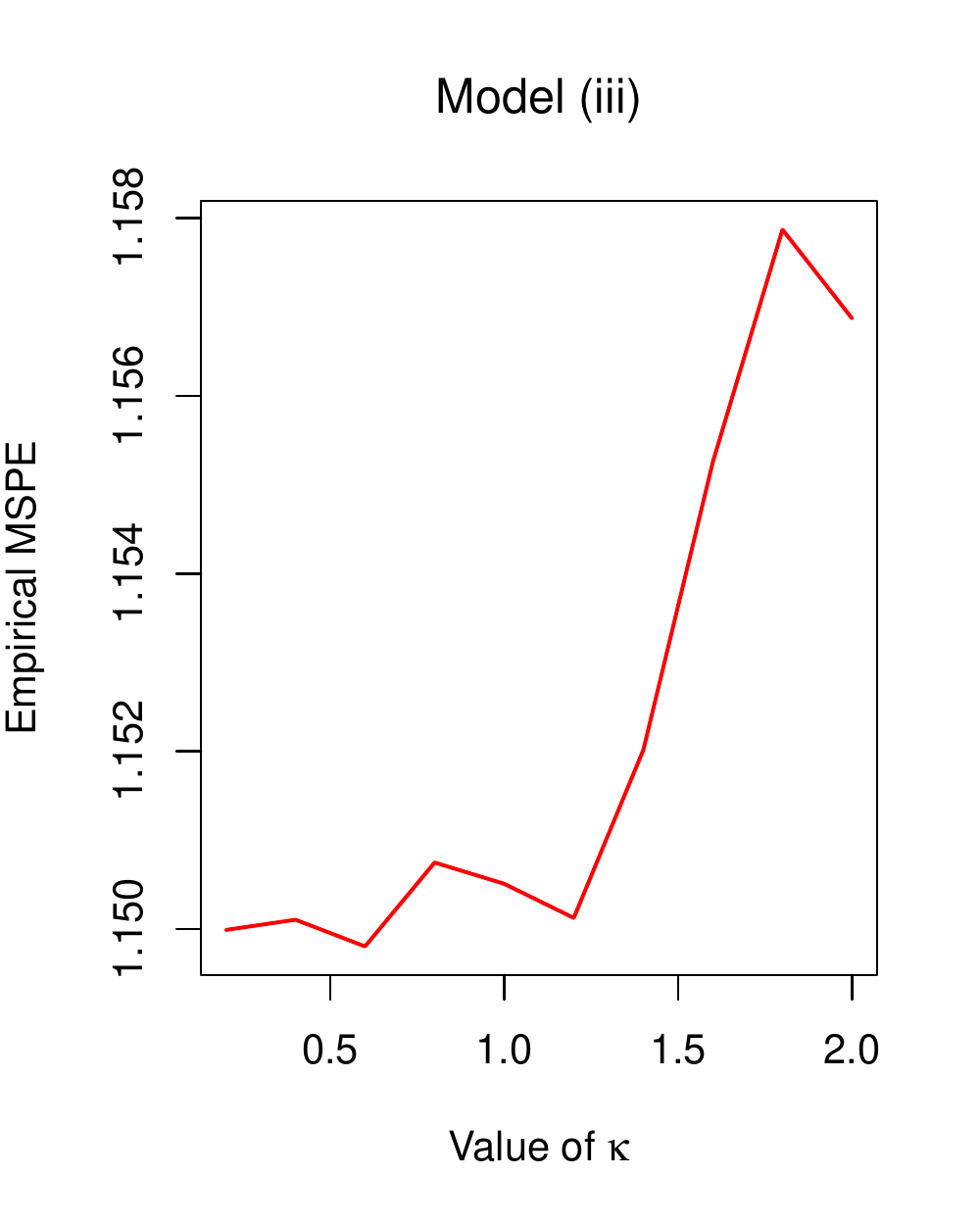} 
\caption{Empirical Mean Square Prediction Error ($\EMSPE$) with respect to $\kappa$ values for the models \ref{M1}, \ref{M2} and \ref{M3}.}
\label{MPE}
\end{figure}

\begin{figure}[!h]
\centering
\includegraphics[scale=0.5]{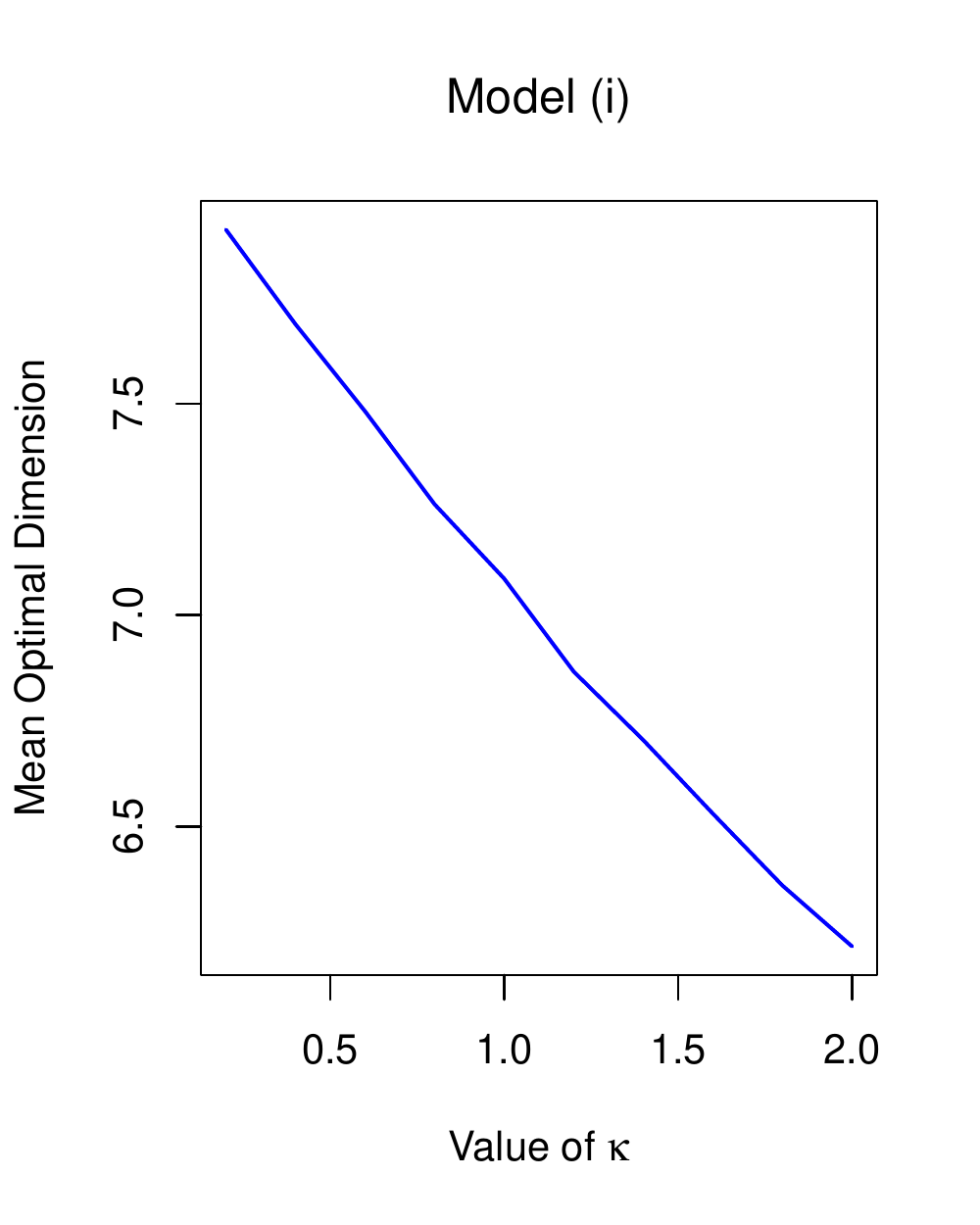} 
\includegraphics[scale=0.5]{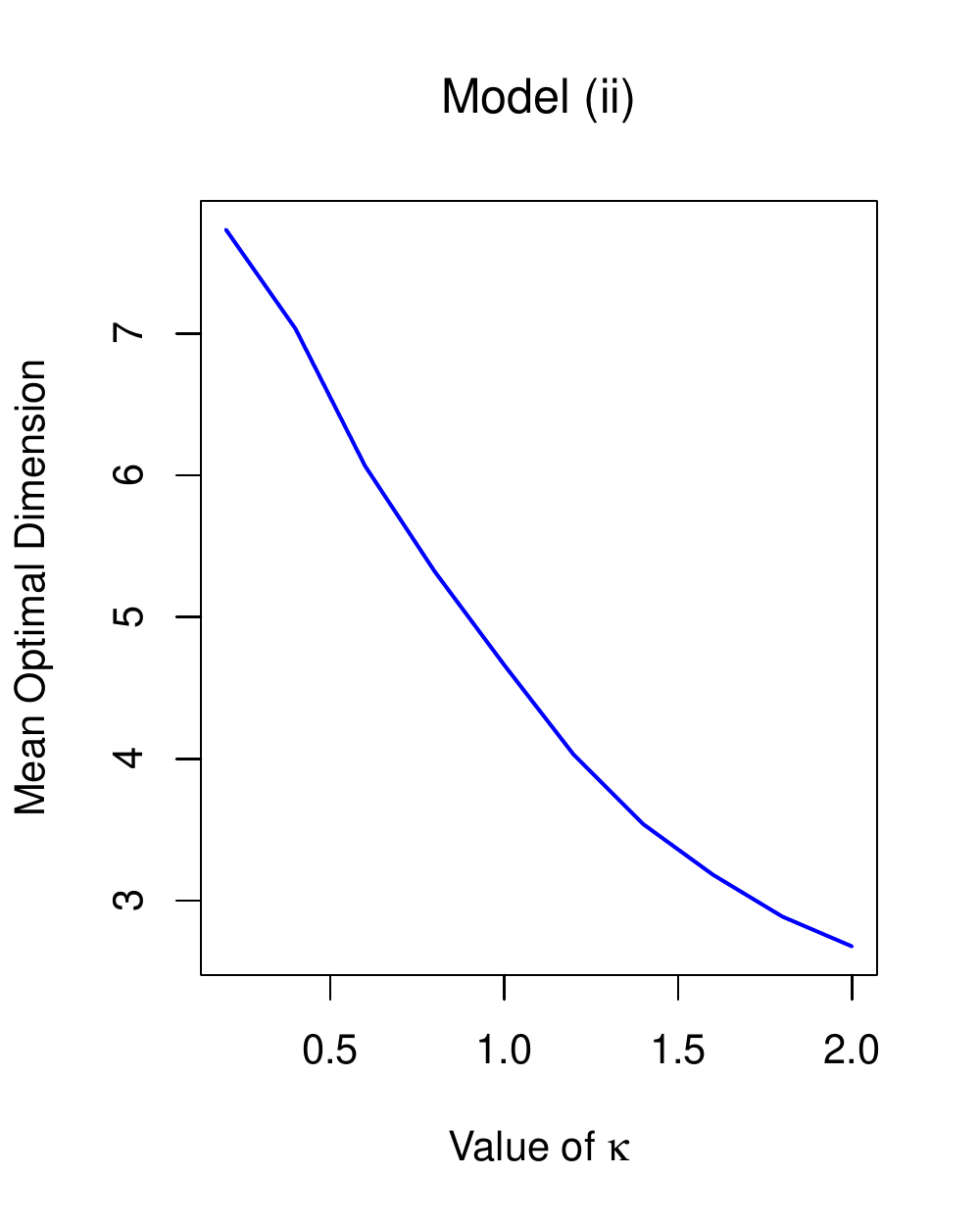} 
\includegraphics[scale=0.5]{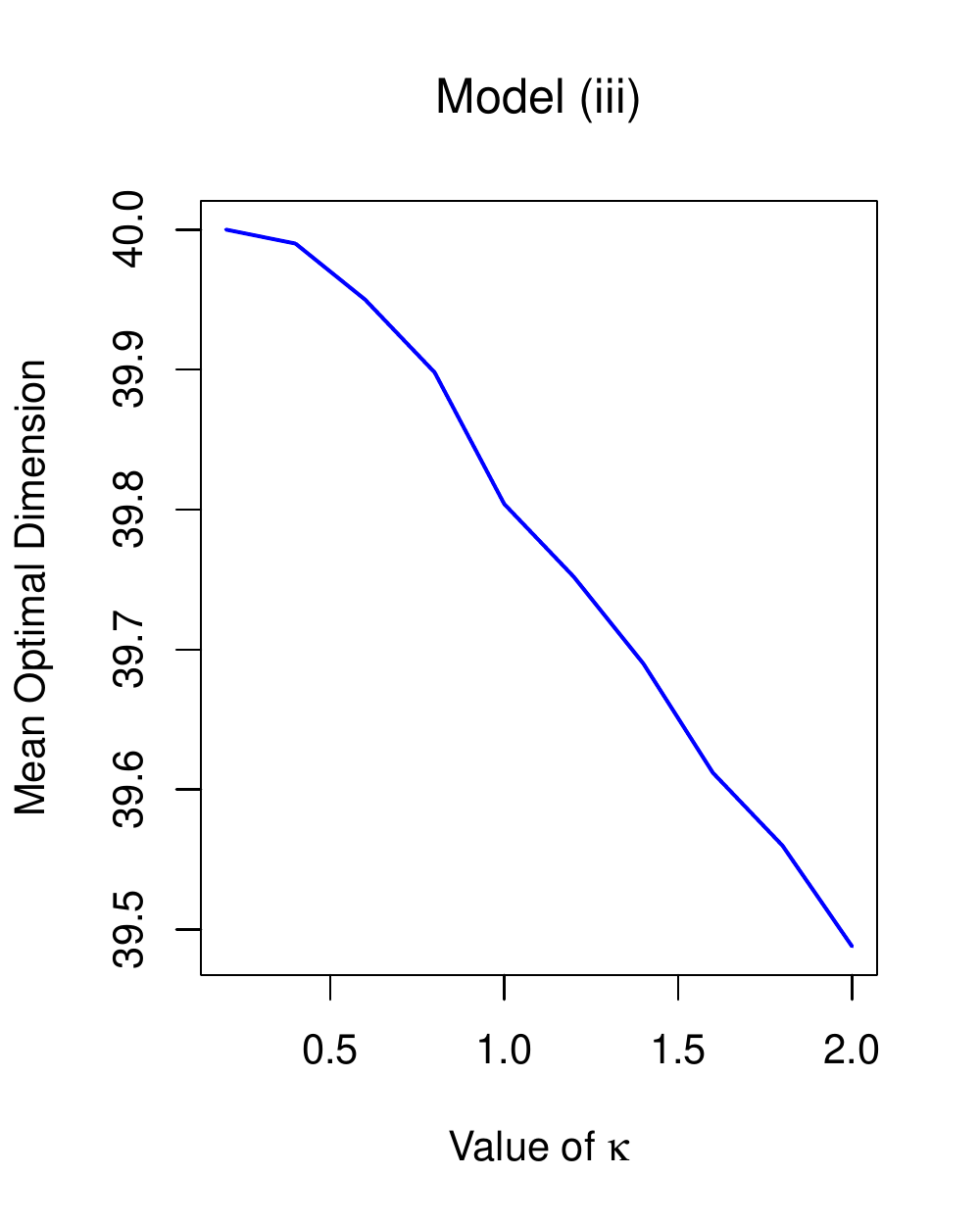} 
\caption{Mean Optimal Dimension with respect to $\kappa$ values for the models \ref{M1}, \ref{M2} and \ref{M3}.}
\label{OD}
\end{figure}

Note also that numerically, we discretize the input $X_i$ (respectively the output $Y_i$) realizations on a $[0,1]$ uniform grid of size $p$ (respectively $q$). The sizes of the grids are chosen to be $p =q =100$. Figure \ref{MPE} represents the Empirical Mean Square Prediction Error with respect to $\kappa$ value, while Figure \ref{OD} shows the mean optimal selected dimension for each $\kappa$ choice. A first general observation of curve shapes in Figure \ref{MPE} is a tendency of decrease then increase.  This reflects the fact that it is not recommended to choose neither too small nor too big calibration parameters. Indeed, small values favor the contrast term, while big values give the advantage to the penalty term, and in both cases the bias/variance compromise is missed. Another intuitive comment when comparing \ref{M1} and \ref{M2} curves in Figure \ref{MPE} is that the Empirical Mean Square  Prediction Error of the second model is much bigger than the first one, which is consistent with the fact that the only difference between these models is that the second one is too noisy compared to the first one. Similar arguments can be used in the comparison of the Mean Prediction Error of the model \ref{M3} with the two other ones. It is also worthwhile to point out that the optimal value is not necessarily unique, which can be suggested by the curve of \ref{M3} in Figure \ref{MPE}. Moreover, the exact numerical values of the optimal parameters for the three models are respectively 0.6, 1.8 and 0.6. In the sequel, $\kappa$ is set to the value 0.6. Furthermore, a simple overview of the graphics shows a systematic decrease of the Mean Optimal Dimension with respect to $\kappa$. This is due to the fact that high $\kappa$ values amplify the penalty, which induce small selected dimensions. Also, by comparison of \ref{M1} and \ref{M2} curves, the selected dimensions for the last model are much smaller than the first one. This is also a result of the noise variance magnitude. The numerical values of the Mean Optimal Dimension in the three cases are respectively 7.482, 2.886 and 39.95.     

\subsubsection{Simulation results}

\begin{figure}[!h]
\centering
\includegraphics[scale=0.5]{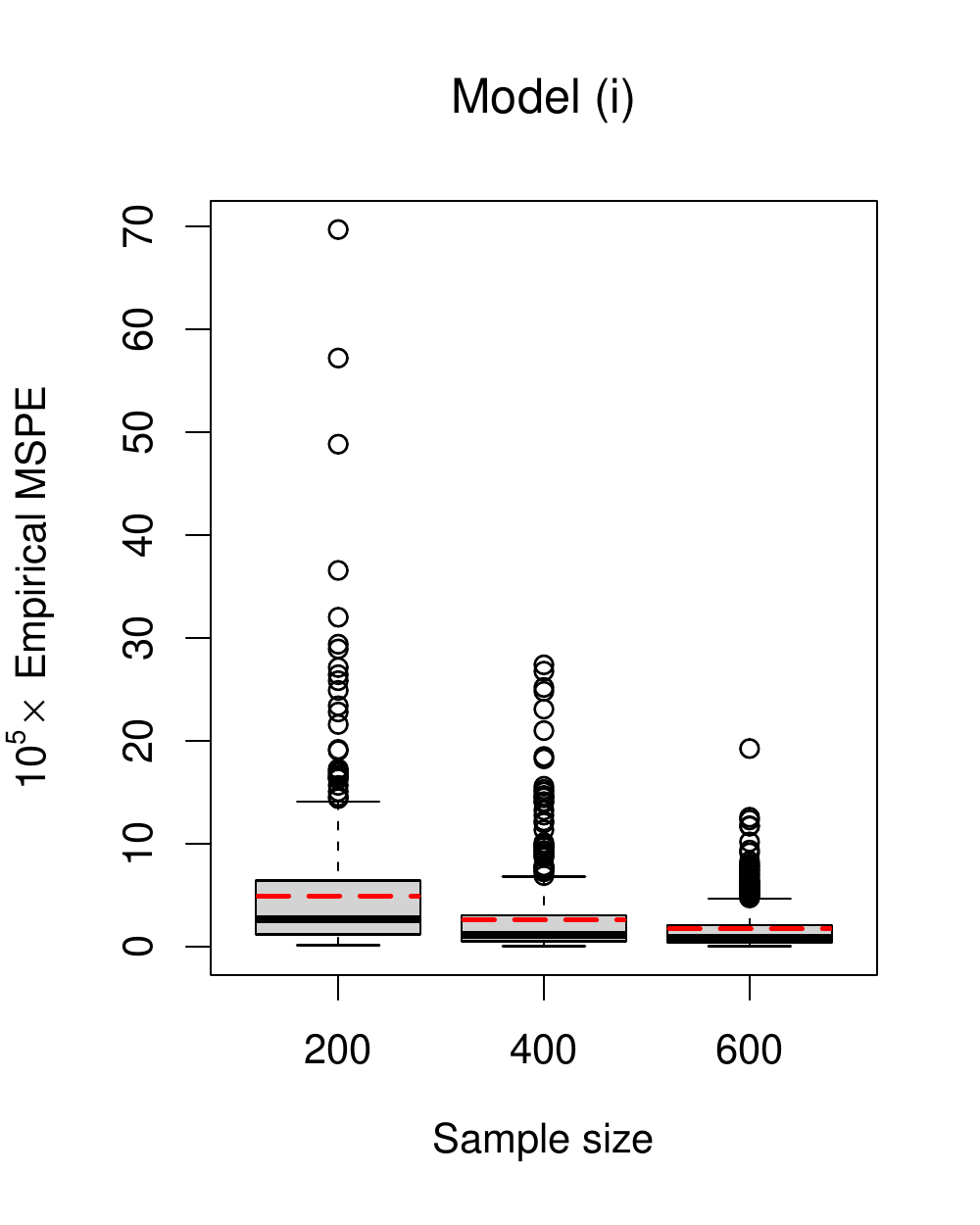} 
\includegraphics[scale=0.5]{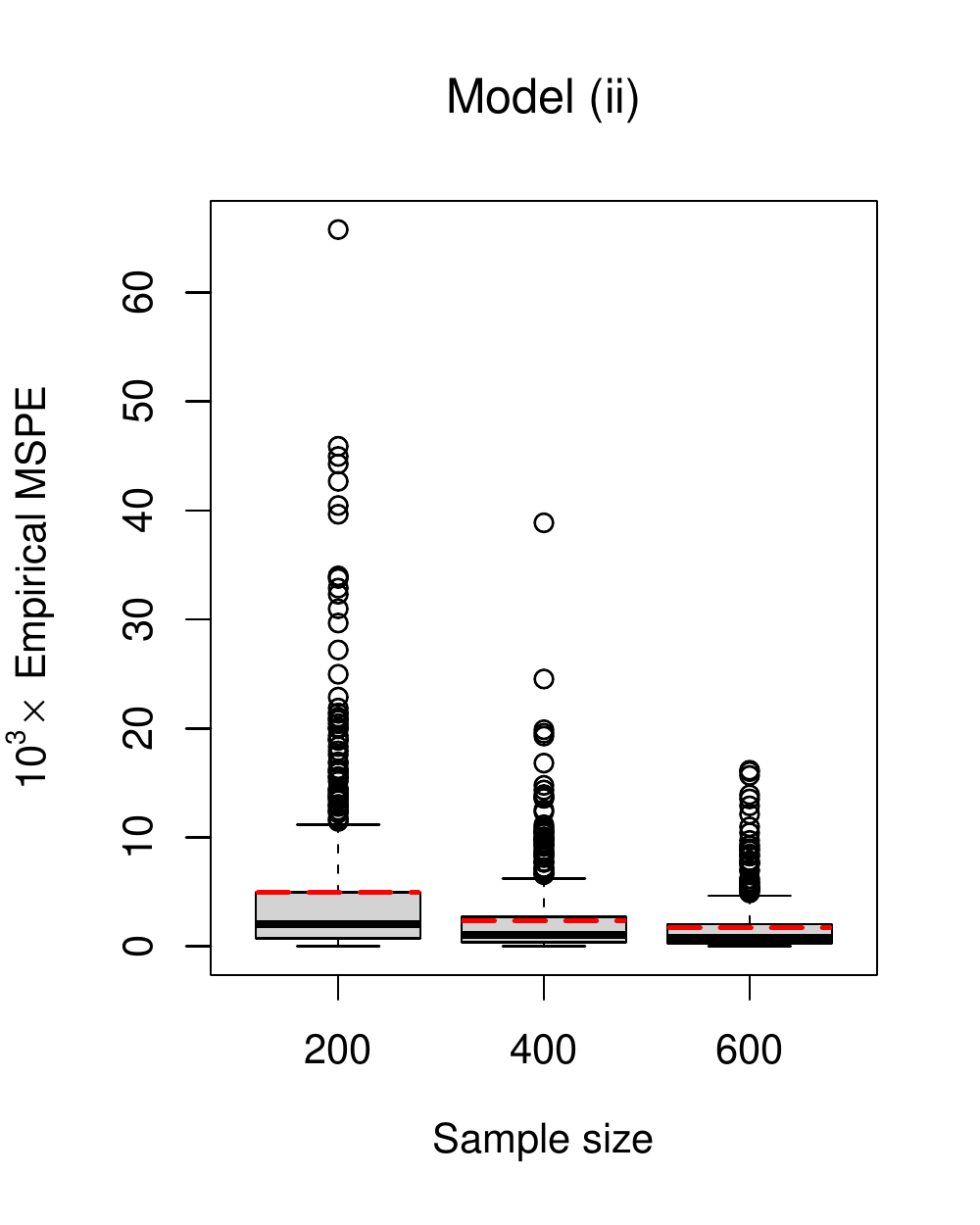} 
\includegraphics[scale=0.5]{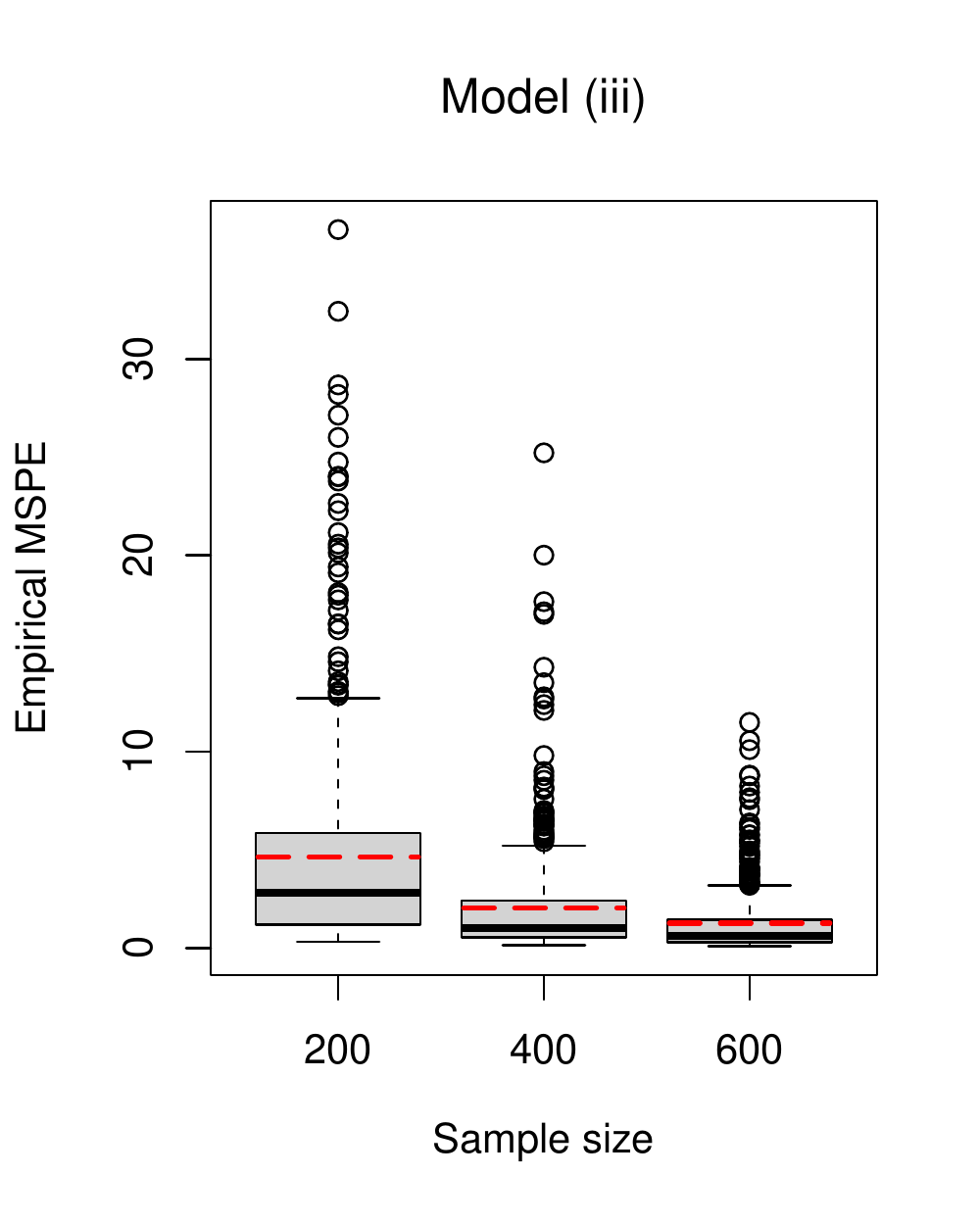} 
\caption{Boxplots of the Mean Square Prediction Errors for the models \ref{M1}, \ref{M2} and \ref{M3} for $\kappa = 0.6$. The mean values are represented in red dashed lines.}
\label{Boxplots}
\end{figure}

Now, we focus on the dispersion of the estimated prediction errors for different sample sizes. To do so, we consider three sample sizes $n = 200, 400, 600$ and we re-estimate $N = 500$ times the prediction errors for each model and $n$ value. As mentioned before, $\kappa$ is set to the value 0.6. The boxplots corresponding to each model are represented in Figure \ref{Boxplots}. As expected, as the sample sizes increase, the boxplots become tighter, the mean prediction errors get closer to zero and the outlier values decrease. This shows an improvement of the estimation accuracy with respect to the sample size, as expected. It is also noticeable that for equal sample sizes, the boxplots of the three models have the same form and a similar distribution of the outliers. This seems to suggest that the prediction quality is robust to the choice of the model and the noise magnitude.

\begin{figure}[!h]
\centering
\includegraphics[scale=0.5]{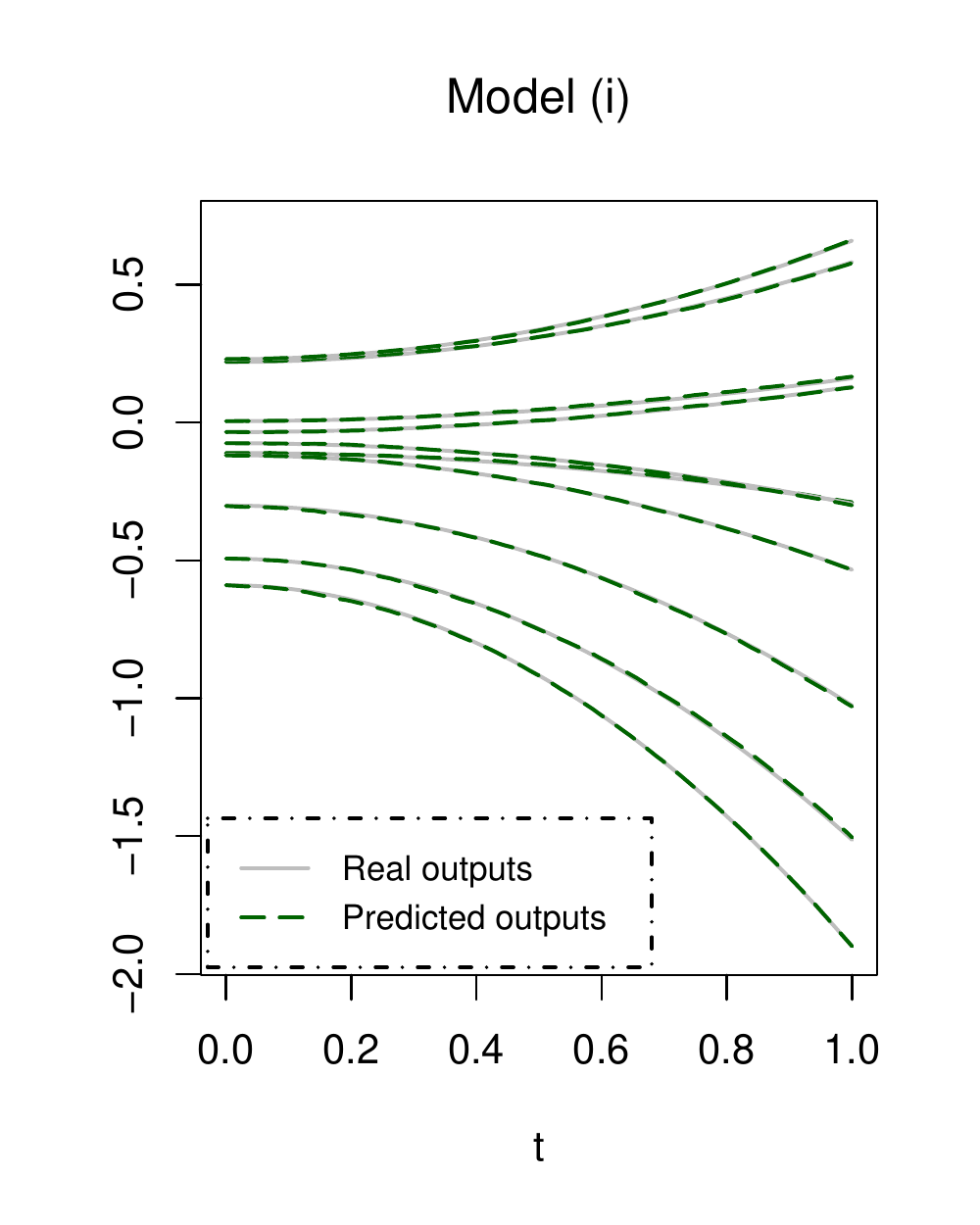} 
\includegraphics[scale=0.5]{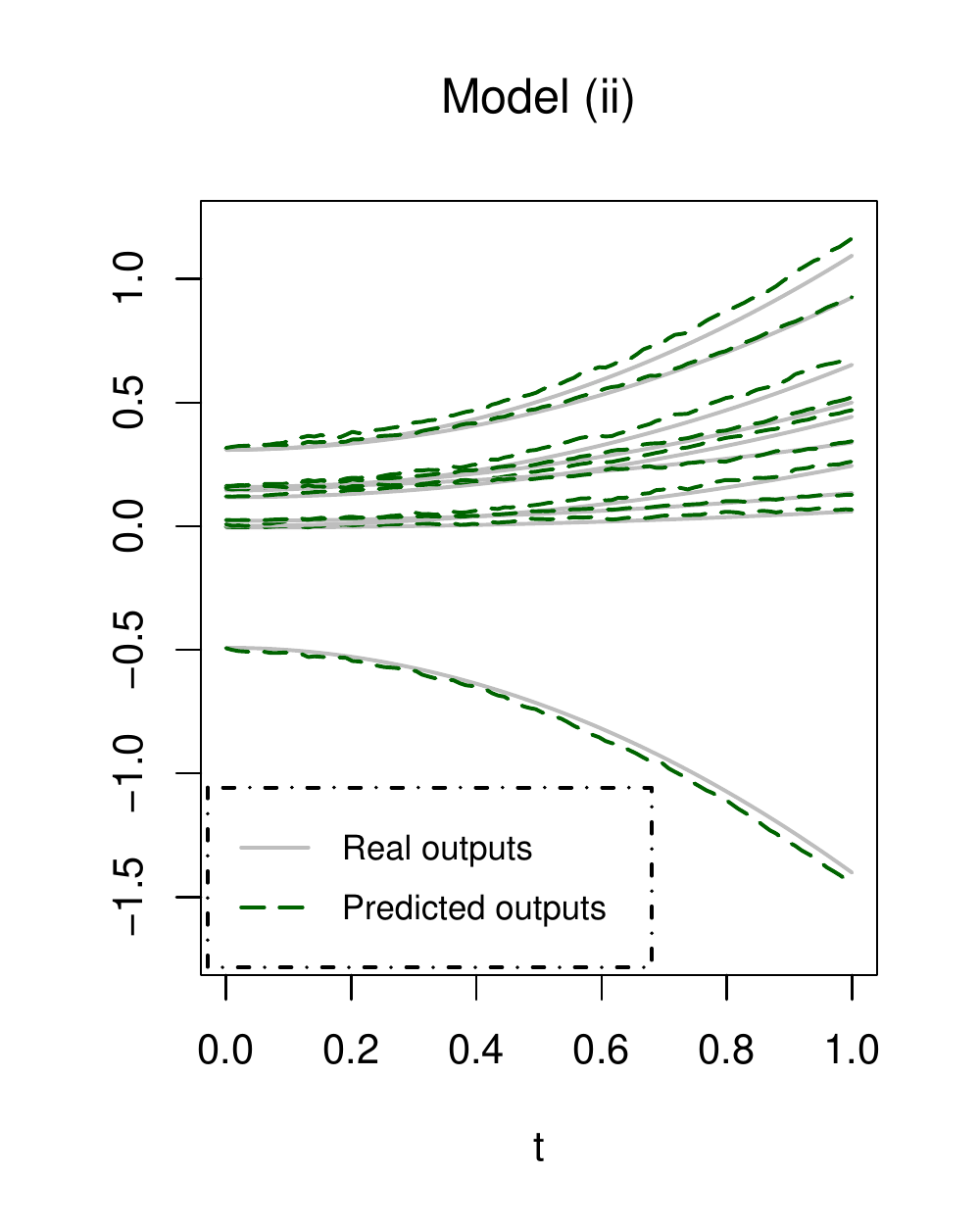} 
\includegraphics[scale=0.5]{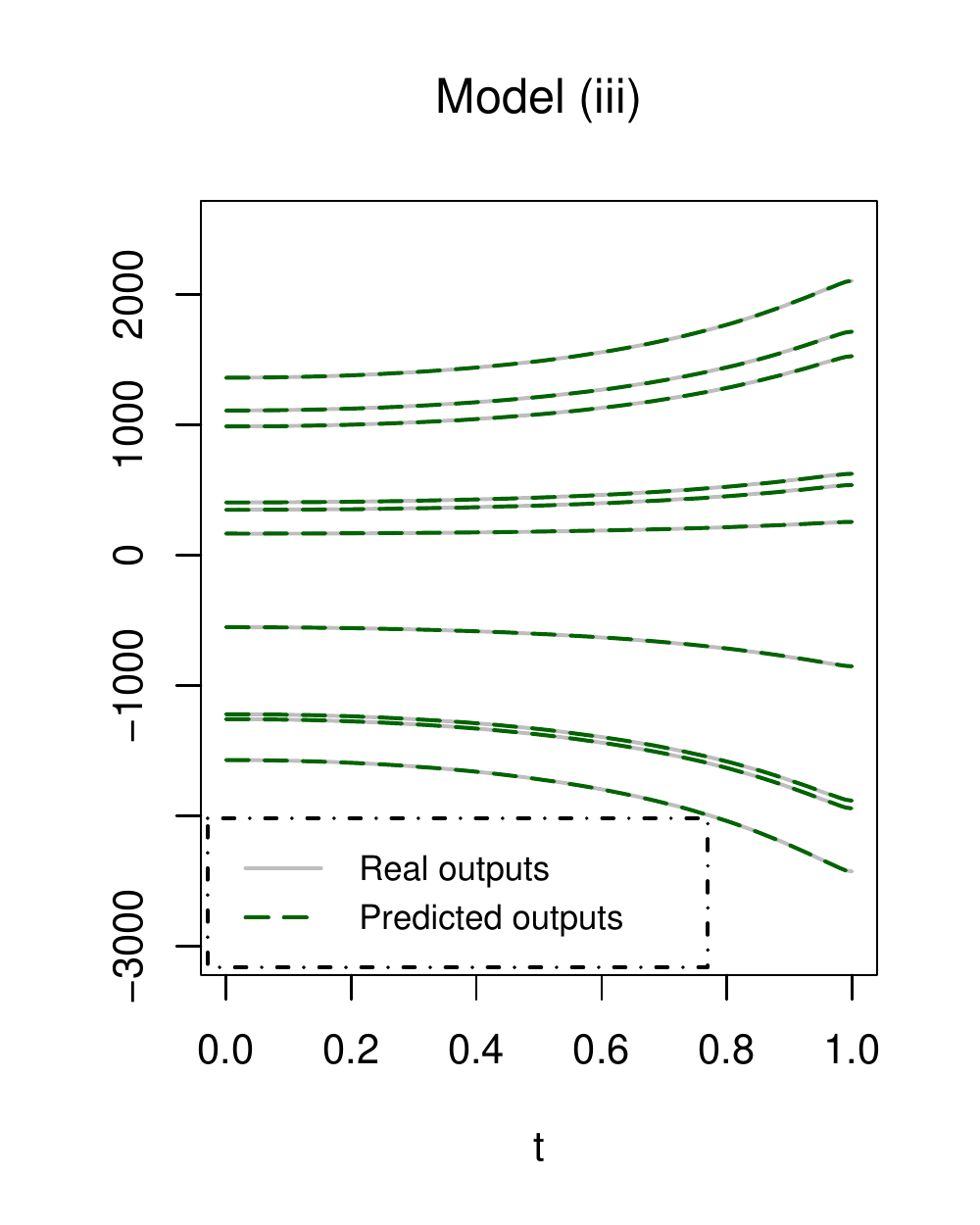} 
\caption{Prediction results : plots of $S^{(\ell)} (X_{n+j}^{(\ell)})$ (plain green lines) versus $\widehat{S}^{(\ell)} (X_{n+j}^{(\ell)})$ (dashed black lines) for $\ell=1,2,3$, $j=1,\ldots,10$, and where  $\widehat{S}_{\ell}$ is computed from a sample of size $n=600$.}
\label{Pred-Mod}
\end{figure}

\smallskip
To illustrate the prediction quality of the proposed adaptive estimators, we assume that for each model (i), (ii), and (iii), an input/output sample $(X_i^{(\ell)},Y_i^{(\ell)})_{i=1,\ldots,n}$ of size $n=600$ is available ($\ell=1,2,3$). These samples are used to estimate the operators $S_{\ell}$, $\ell=1,2,3$  and we predict the model output for 10 new independent inputs, denoted $X_{n+1}^{(\ell)}, \ldots, X_{n+10}^{(\ell)}$. Figure \ref{Pred-Mod} shows the obtained graphs for $S_{\ell} (X_{n+j}^{(\ell)})$ and $\widehat{S}_{\ell} (X_{n+j}^{(\ell)})$, with $\ell = 1,2,3$ and $j = 1, \ldots, 10$, while $S_{\ell}$, $\widehat{S}_{\ell}$ respectively denote the real and estimated slope operators. In general, the prediction is quite accurate. Once again, a large noise magnitude deteriorates the prediction quality, which can be observed by comparing the graphs of the first two models.

\subsection{Real data case} \label{sec:Rdata}
\subsubsection{Application to the prediction of electricity consumption}
\begin{figure}[!h]
\begin{tabular}{cc}
Original data & Transformed data\\
\includegraphics[width=0.45\textwidth]{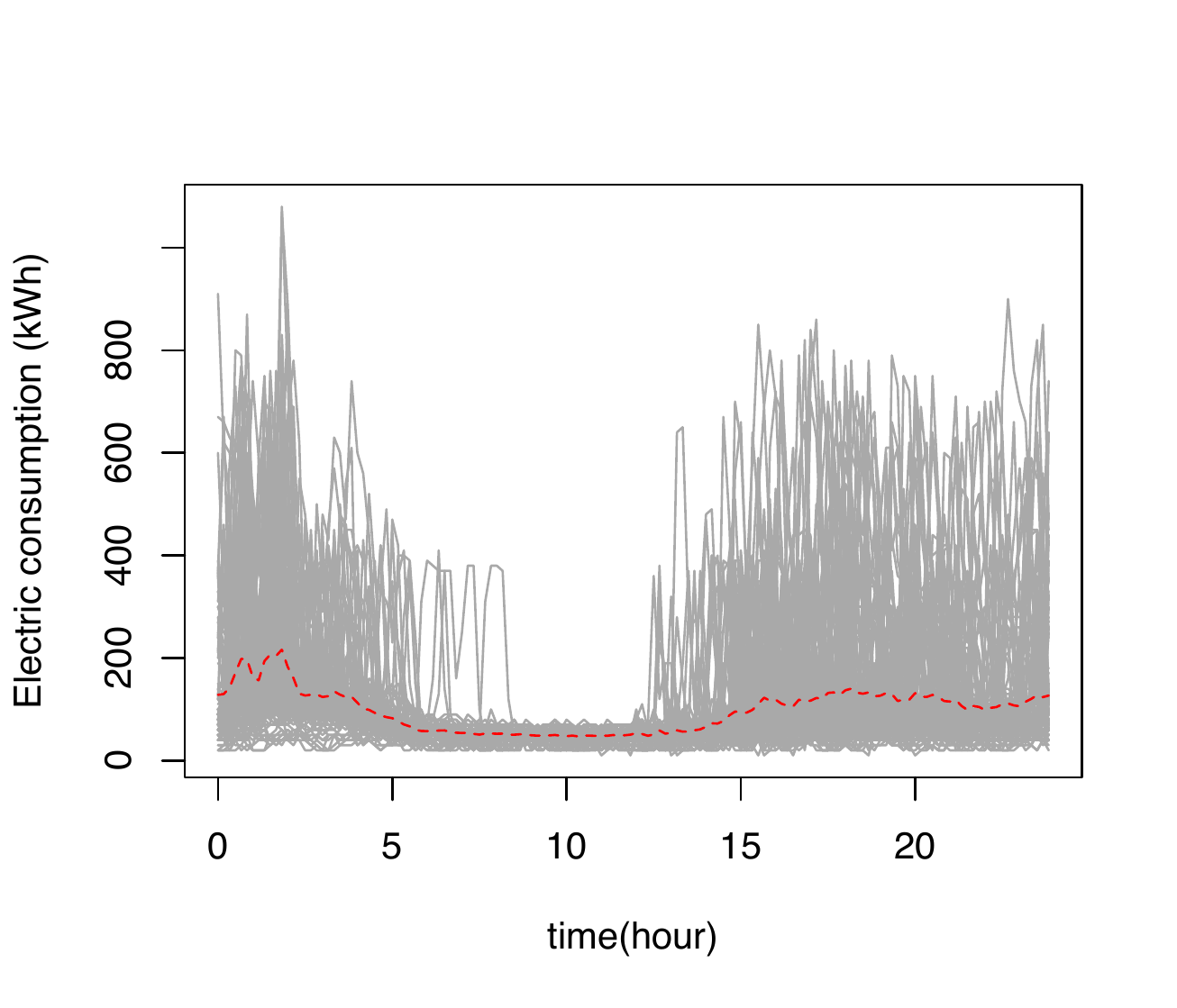}&\includegraphics[width=0.45\textwidth]{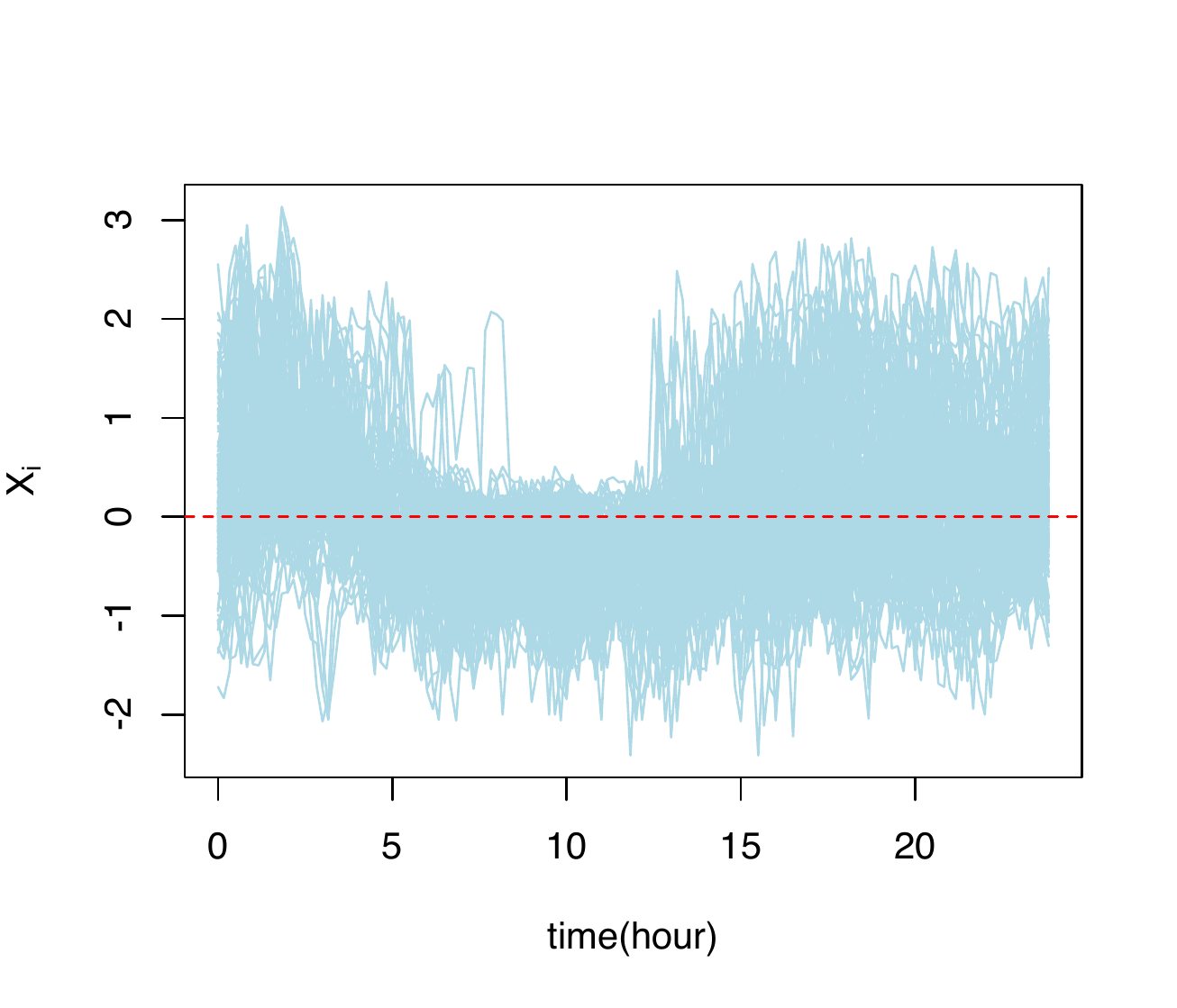}
\end{tabular}
\caption{\label{fig:energy_data} Evolution of electric consumption of appliances during $n=136$ days (original data, thin gray lines) and functions of the sample (transformed data : centered version of the logarithm of the original data, thin blue lines). The dashed red lines are the empirical mean of each sample.}
\end{figure}

 The data we study are the electricity consumption of appliances curve of a low energy house located in Stambrudge (Belgium). The dataset is freely available on UCI Machine Learning Repository \url{https://archive.ics.uci.edu/} and has been studied by \citet{CFD17}. It consists on measurements on 24 variables every 10 minutes from 11th january, 2016, 5pm to 27th may, 2016, 6pm. The variable of interest is the consumption of appliances, which is the main source of energy consumption. The data consists of a $d$-dimensional times series, with $d=24$. It is first transformed into a sample of functional data by splitting the data day by day. We can deduce from the variable selection study conducted in \citet{Roche_preprint} that the most important variable to predict  appliances electricity consumption of day $i$ is the appliances electricity consumption of day $i-1$, and that a $\ln$-transformation of the covariates seems to lead to better results. Then, in our study, the variable to predict $Y_i$ is the log of the appliances energy consumption of day $i$ and $X_i$ is the log of appliances energy consumption of day $i-1$. The data are also recentered. We present in Figure~\ref{fig:energy_data} the original and transformed data. 

The incorporation of other functional covariates could be of great interest for the application but is out of the scope of the paper. 

Another difficulty for the estimation procedure is that it requires the knowledge of the trace of the noise operator $\sigma_\varepsilon^2$, which is unknown in practice.  To get around this difficulty, we adapt the method proposed in \citet{brunel2016non}, consisting in replacing the unknown quantity $\sigma_\varepsilon^2$ in criterion~\eqref{eq:selecm1} by the contrast $\gamma_n(\widehat S_{m_1,\infty})$. In model selection in regression contexts, this method shows strong similarities with the one of \citet{BGH14}. In the context of the functional linear model with scalar output, it has been proven in \citet{brunel2016non}  that the estimator selected by this fully data-driven criterion verifies an oracle-type inequality, achieves the same minimax rates as the estimator selected by the criterion depending on the noise variance and that it does not change significantly the practical performances of the estimator. 

\begin{figure}[!h]
\begin{tabular}{cc}
Dimension selected & $ L^2$ prediction error of selected estimator\\
\includegraphics[width=0.45\textwidth]{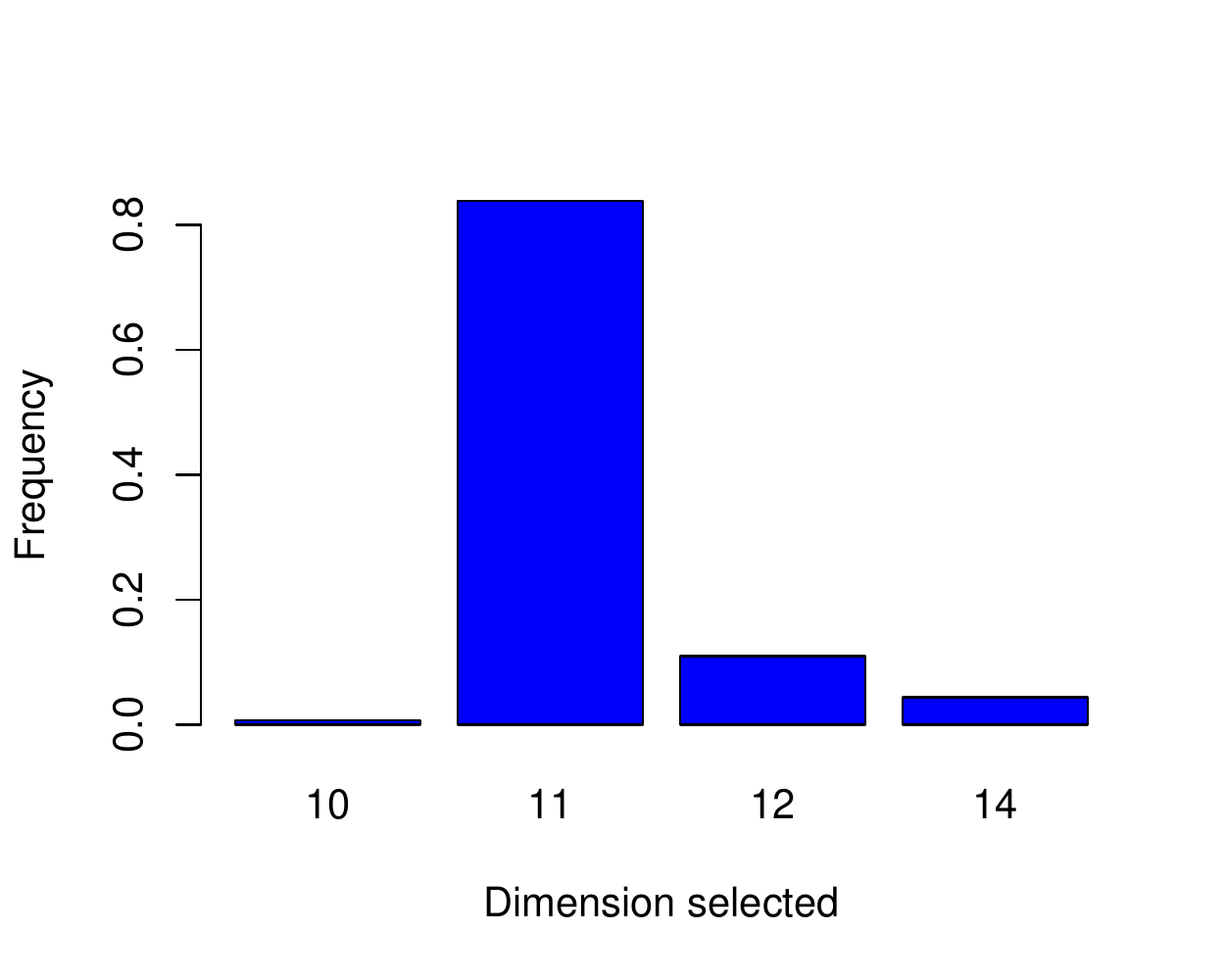}&\includegraphics[width=0.45\textwidth]{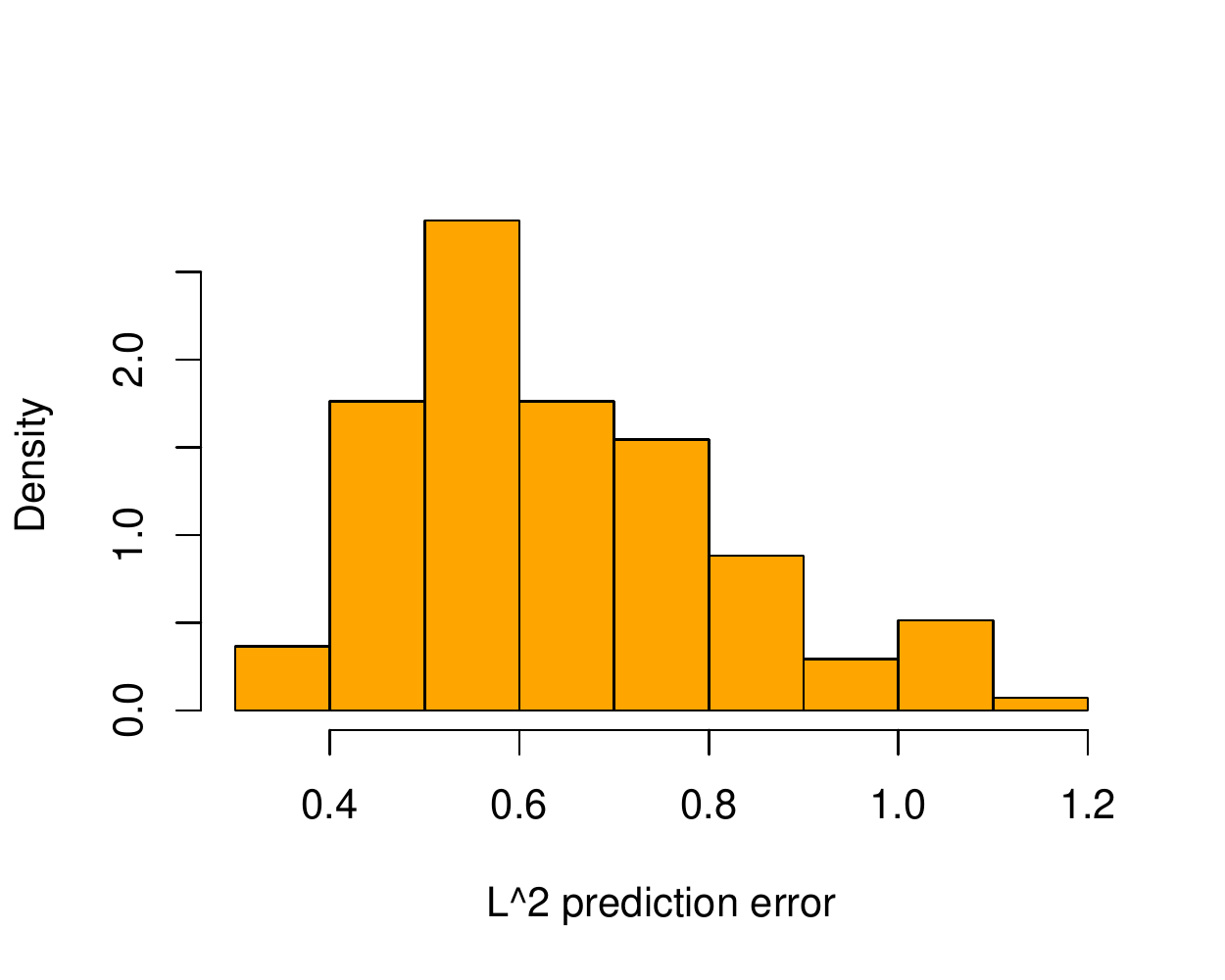}
\end{tabular}
\caption{\label{fig:energy_data_dimandrisk} Dimension selected and $L^2$ prediction error $\|Y_i-\widehat Y_i^{(-i)}\|$ of the estimator calculated from each cross-validated sample.}
\end{figure}

As suggested by the simulation study, the value of $\kappa$ is also fixed to $\kappa=0.6$. To study the selected dimension, the risk of the estimators and their stability, we perform cross-validation of the sample: for each day $i$, we calculate the selected dimension  $\widehat m_1^{(-i)}$ and the $L^2$-prediction error of the estimator $\widehat S_{\widehat m_1^{(-i)},\infty}^{-i}$ calculated from the sample $\{(X_j,Y_j), j\neq i\}$. The results are presented in Figure~\ref{fig:energy_data_dimandrisk}. The dimension selection procedure is quite stable, selecting more than 80\% of time the dimension $\widehat m_1=11$ and the $L^2$-prediction error does not explode for some observations. 

\begin{figure}[!h]
\includegraphics[width=\textwidth]{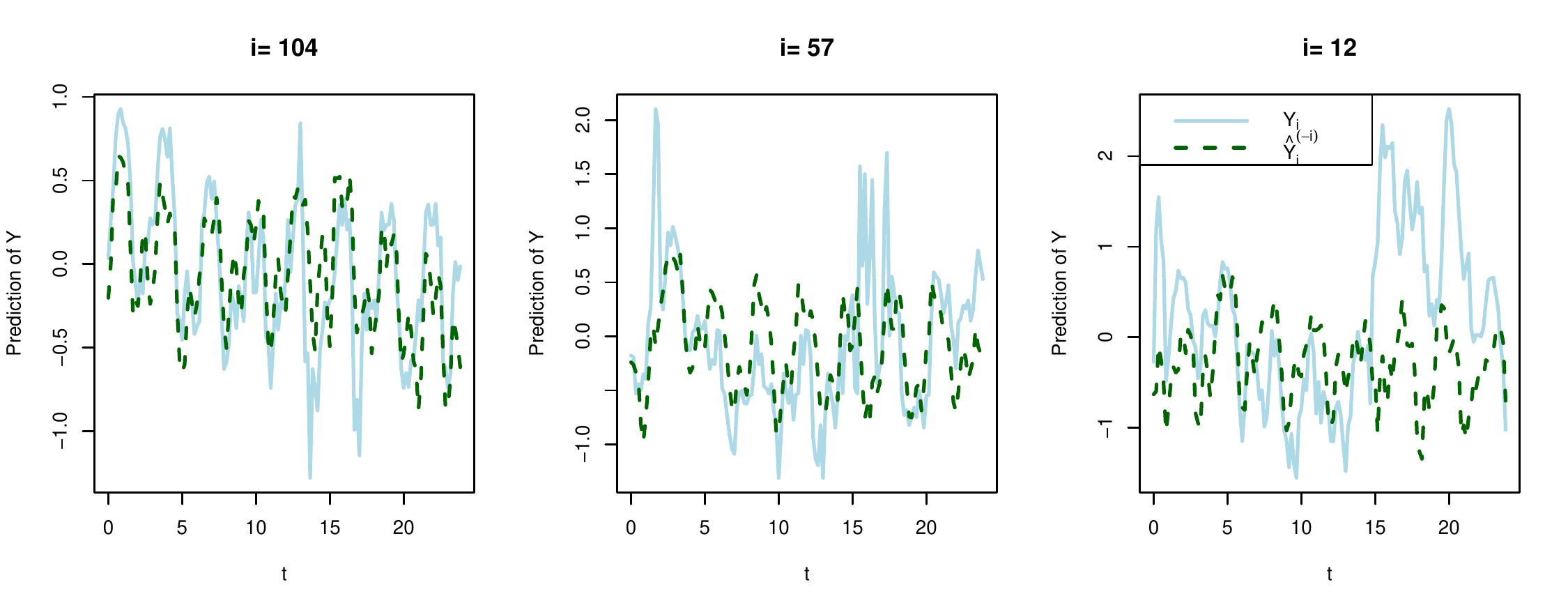}
\caption{\label{fig:energy_data_predY} Cross-validated prediction $\widehat Y_i^{-i}$ made for three days (the days where the prediction is best, median and worst).}
\end{figure}

We also plot in Figure~\ref{fig:energy_data_predY}, for three well-chosen days $i$ ($i=104$ is the day for which the distance $\|Y_i-\widehat Y_i^{-i}\|$ is minimal, $i=4$ corresponds to the median prediction error and $i=83$ to the maximal prediction error), the true value of $Y_i$ and its prediction $\widehat Y_i^{-i}=\widehat S_{\widehat m_1^{(-i)},\infty}^{-i}(X_i)$. 

\begin{figure}[!h]
\includegraphics[width=\textwidth]{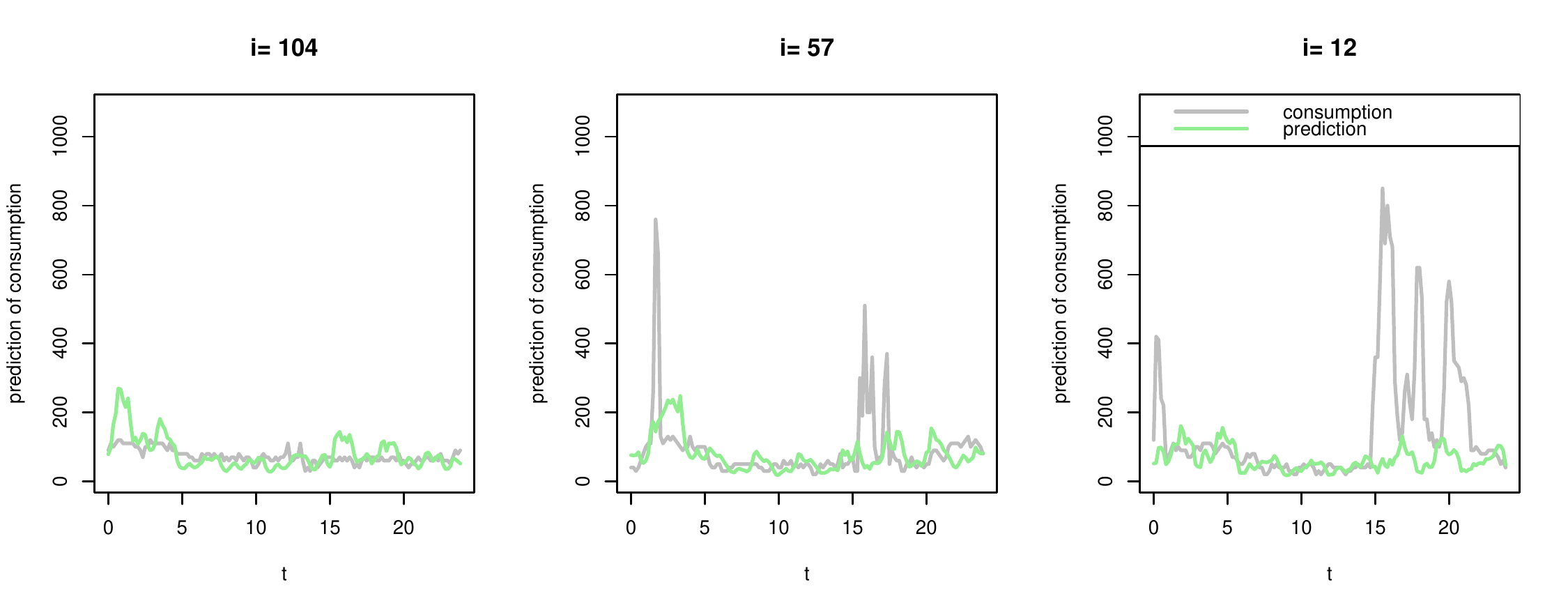}
\caption{\label{fig:energy_data_pred} Prediction of appliances energy consumption.}
\end{figure}

Figure~\ref{fig:energy_data_pred} represents, for the same days, the prediction of appliances energy consumption (after adding the mean and taking the exponential).

We see in Figure \ref{fig:energy_data_pred} that prediction captures trends well and that the worst prediction seems to be due to a brutal change of behavior of the appliances  electricity consumption which is quite hard to predict and may be due to external factors (hence unavoidable with our model).

\subsubsection{Application to the prediction of prices from wind power infeed}

\begin{figure}[!h]
\begin{tabular}{ccc}
&Original data & Transformed data\\
$X$&\includegraphics[width=0.45\textwidth]{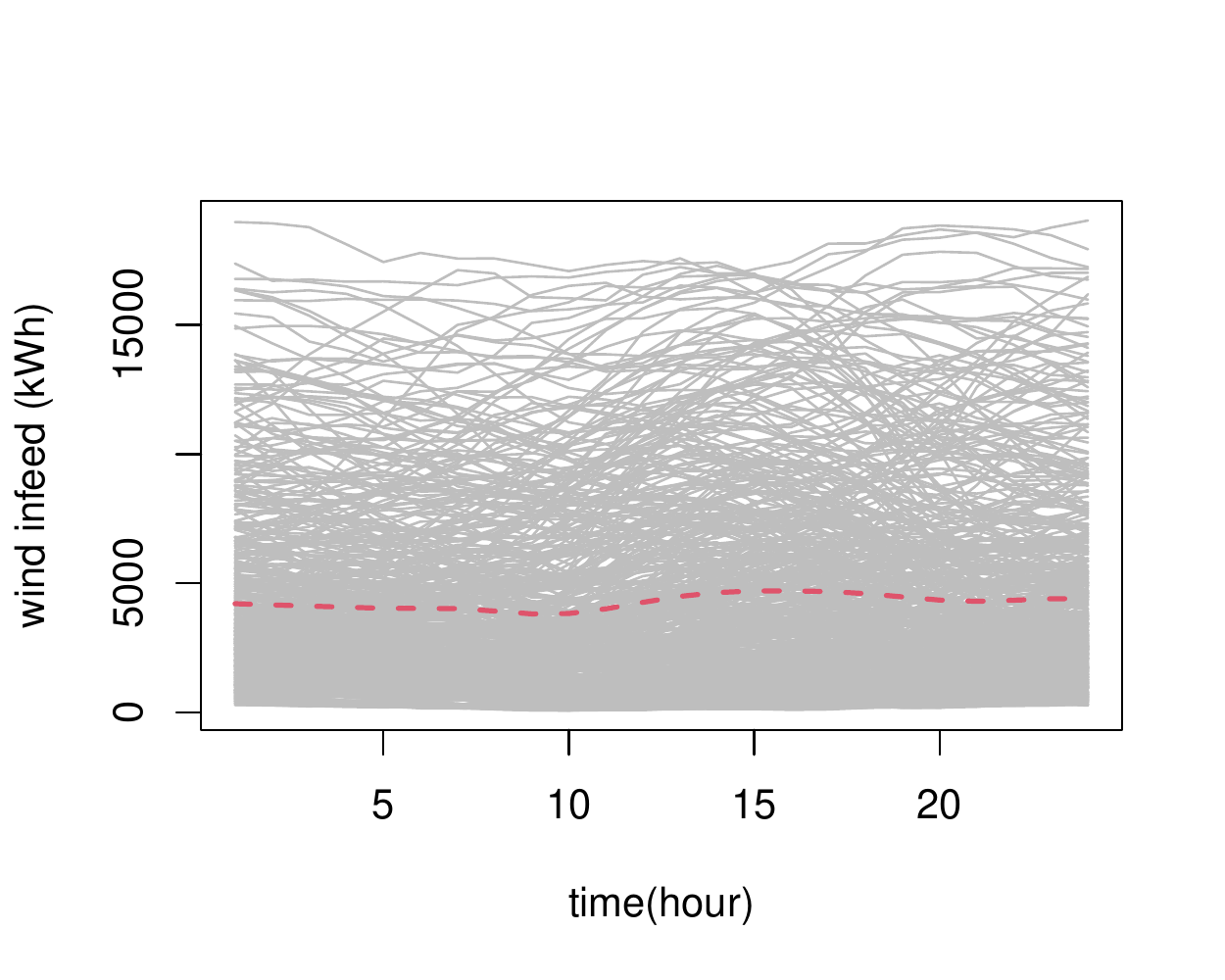}&\includegraphics[width=0.45\textwidth]{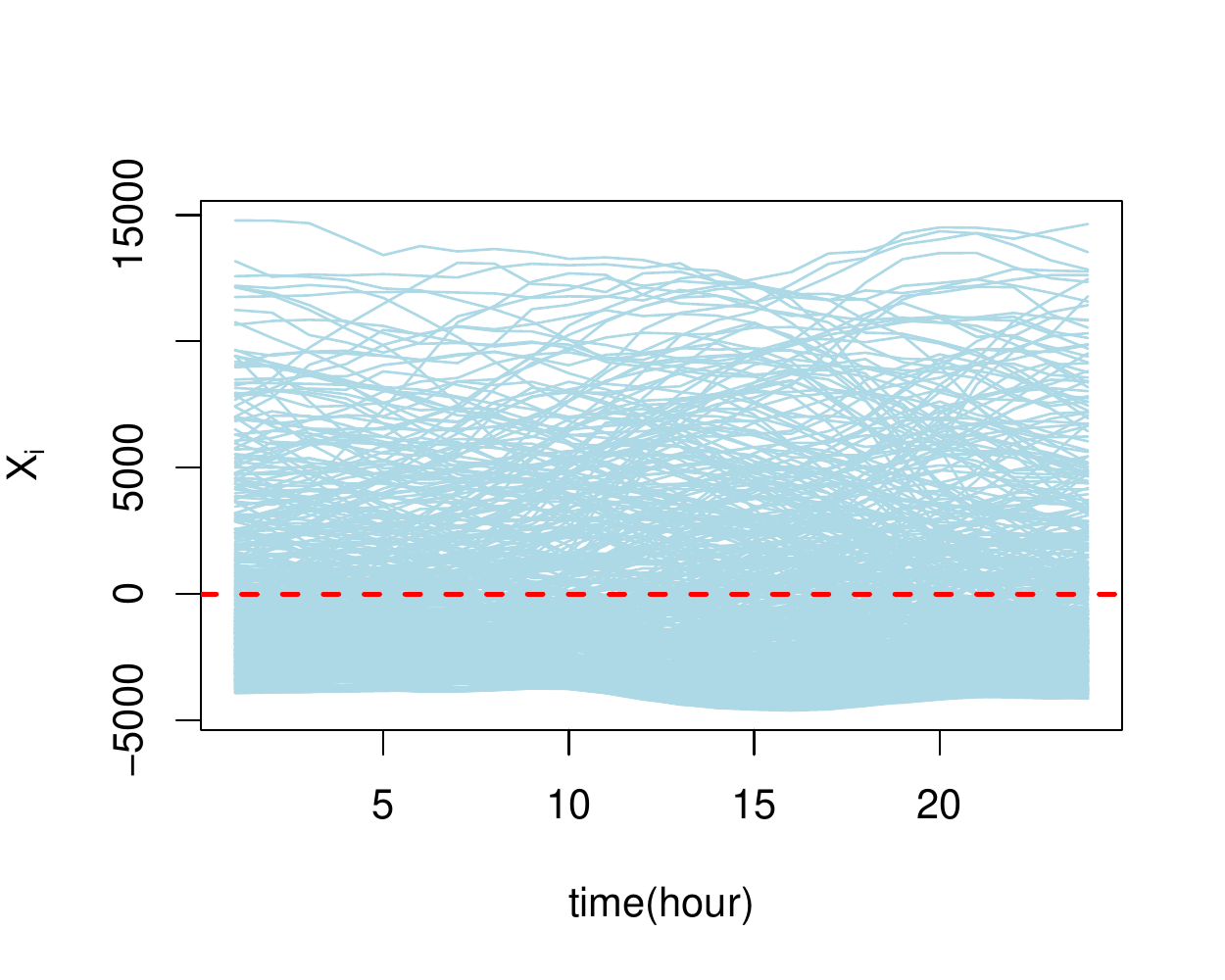}\\
$Y$&\includegraphics[width=0.45\textwidth]{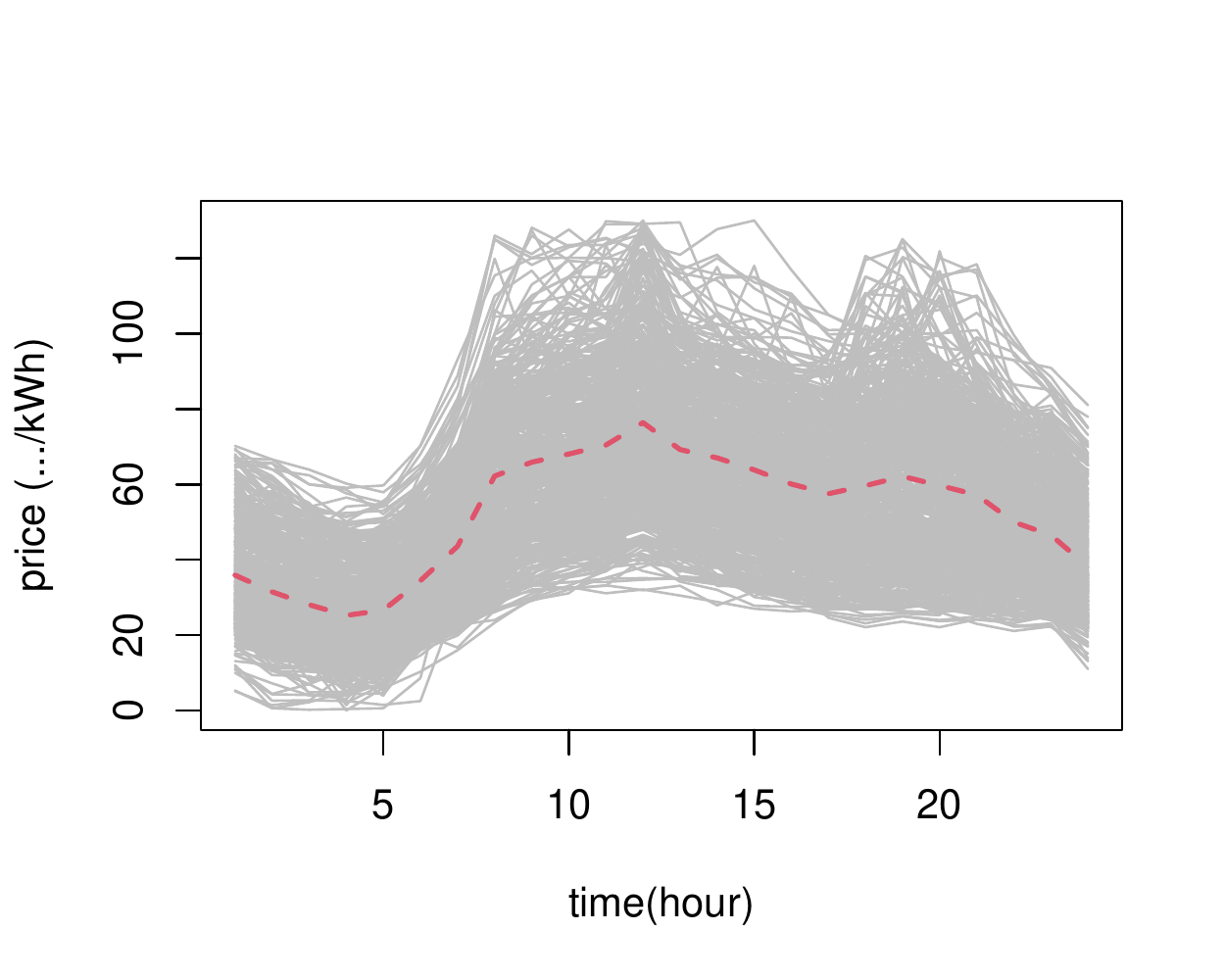}&\includegraphics[width=0.45\textwidth]{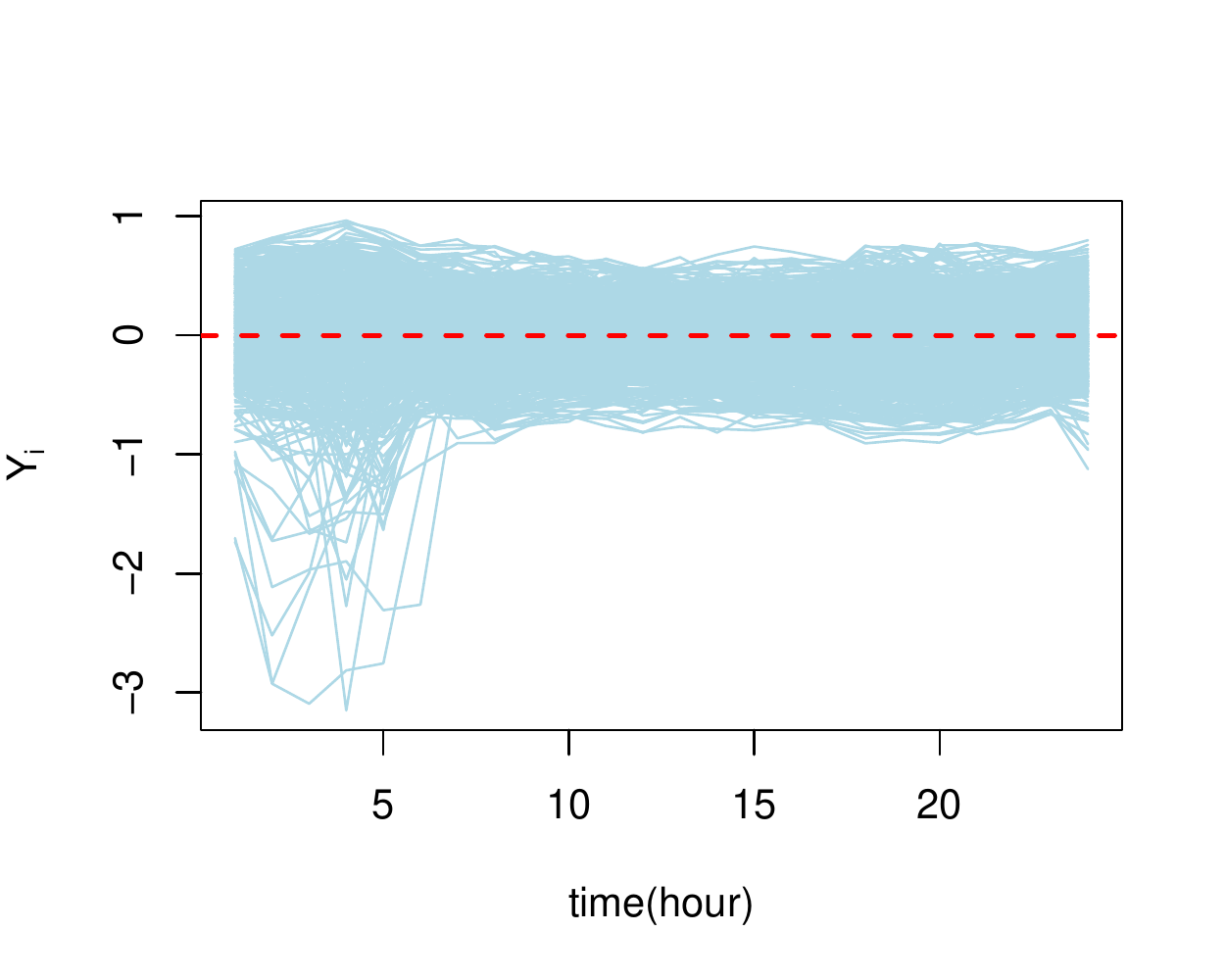}\\
\end{tabular}
\caption{\label{fig:WindPrices_data} First line: evolution of wind power in-feed during $n=516$ days (original data, thin gray lines) and functions of the sample (transformed data : centered version of the original data, thin blue lines). Second line: evolution of prices (original data, thin gray lines) and functions of the sample (transformed data : centered version of the log  the original data+1, thin blue lines). The dashed red lines are the empirical mean of each sample.  }
\end{figure}

We also apply now our estimation method to another, and more difficult prediction problem. The aim here is to predict the evolution of electricity prices in Germany from the wind power in-feed. This dataset has been extensively studied by \citet{liebl2013modeling,imaizumi2018pca} and  can be found at \url{https://www.dliebl.com/#publications}. We first remark that some observations exhibit non standard behaviors, in particular some prices are particularly elevated. Then, we start be removing the outliers that deserves a particular study which is out the scope of the paper. We consider a day to be an outlier if the maximal value of the price of the day is larger than $Q_3+1.5(Q_3-Q_1)$ where $Q_1$ (resp. $Q_3$) corresponds to the first (resp. third) quartile of the maximal prices of each day. Then the data are also centered. We present in Figure~\ref{fig:WindPrices_data} the original and transformed data. 

\begin{figure}[!h]
\begin{tabular}{cc}
Dimension selected & $L^2$ prediction error of selected estimator\\
\includegraphics[width=0.45\textwidth]{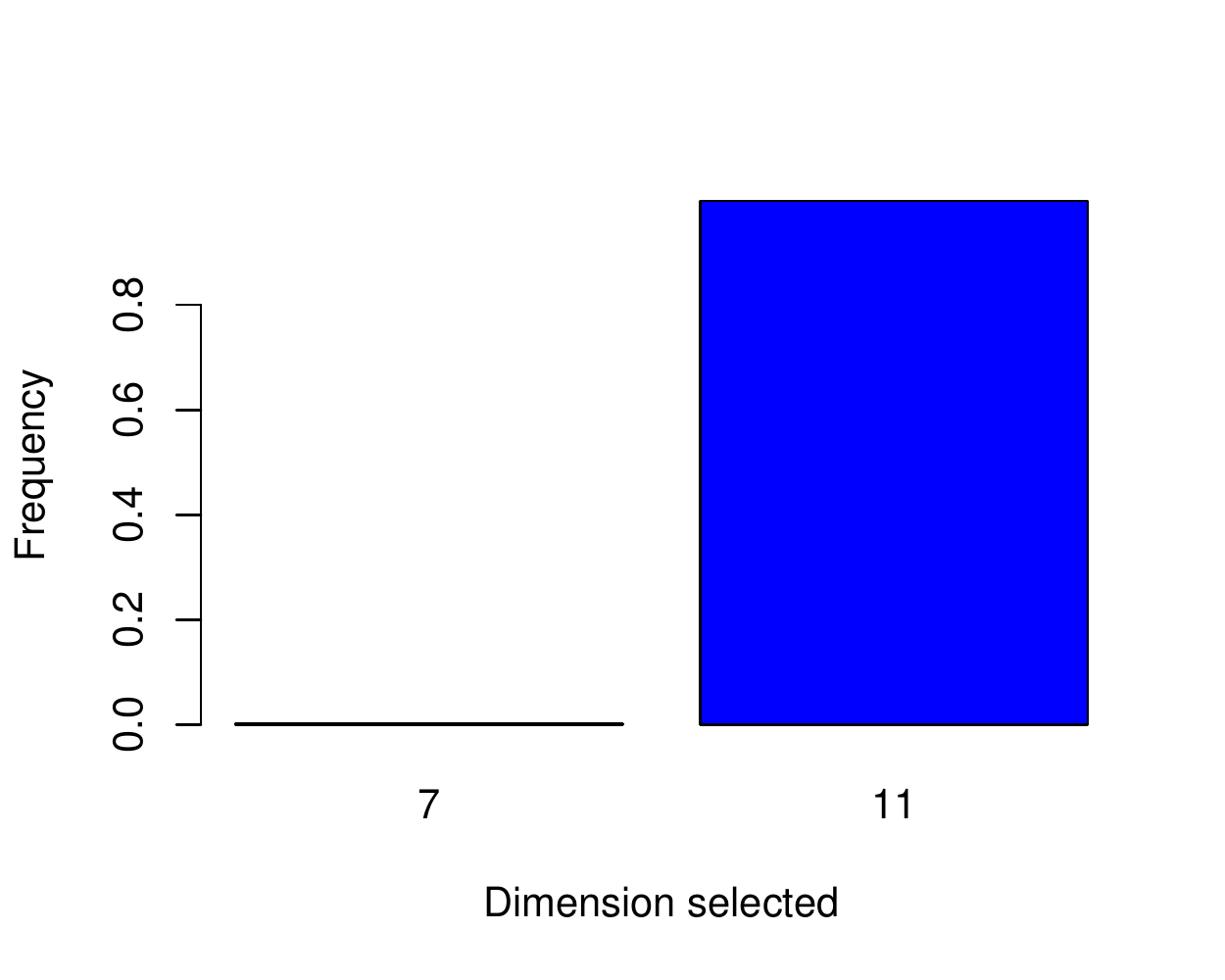}&\includegraphics[width=0.45\textwidth]{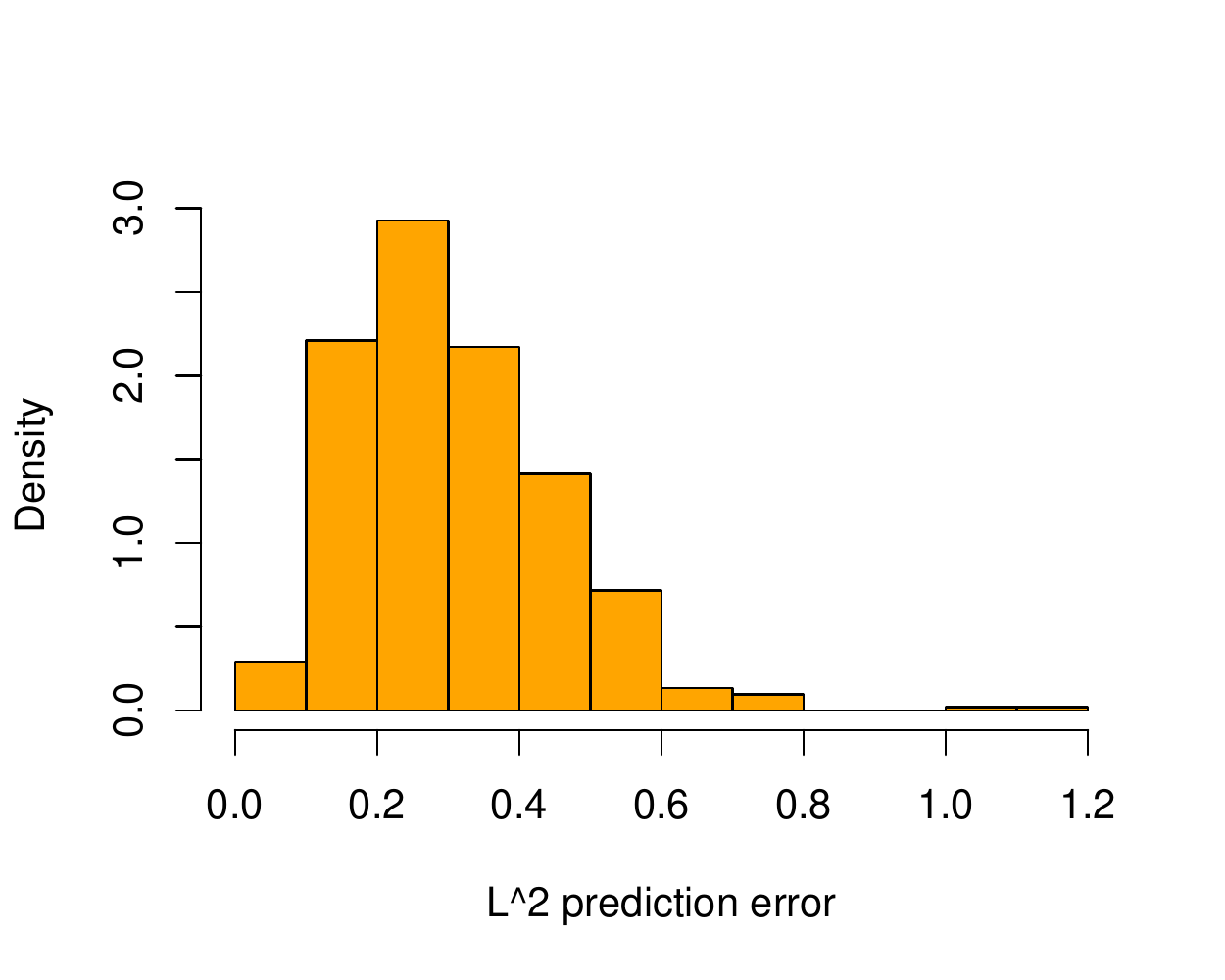}
\end{tabular}
\caption{\label{fig:WindPrices_data_dimandrisk} Dimension selected and $ L^2$ prediction error $\|Y_i-\widehat Y_i^{(-i)}\|$ of the estimator calculated from each cross-validated sample.}
\end{figure}

\smallskip
As in the previous section, we set $\kappa=0.6$ and replace the unknown quantity $\sigma_\varepsilon^2$ be $\gamma_n(\widehat S_{m_1,\infty})$. We also performed a cross-validation of selected dimensions and associated prediction risk. The results are presented in Figure~\ref{fig:WindPrices_data_dimandrisk}. 

\begin{figure}[!h]
\includegraphics[width=\textwidth]{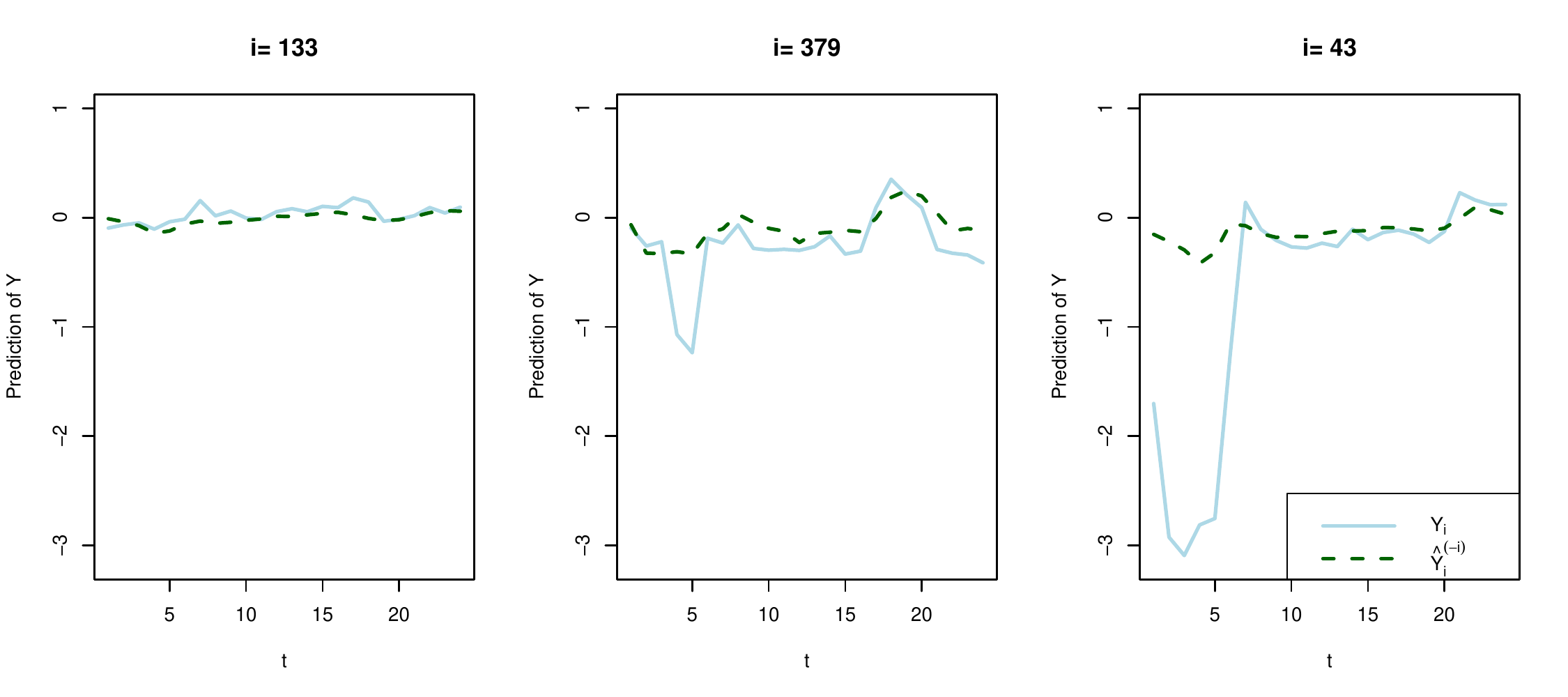}
\caption{\label{fig:WindPrices_data_predY} Cross-validated prediction $\widehat Y_i^{-i}$ made for three days.}
\end{figure}

\smallskip
We also plot in Figure~\ref{fig:WindPrices_data_predY}, for three well-chosen days $i$  the true value of $Y_i$ and its prediction $\widehat Y_i^{-i}=\widehat S_{\widehat m_1^{(-i)},\infty}^{-i}(X_i)$ : $i=133$ is the day for which the distance $\|Y_i-\widehat Y_i^{-i}\|$ is minimal, $i=379$ corresponds to the median prediction error and $i=43$ to the maximal prediction error.

\begin{figure}[!h]
\includegraphics[width=\textwidth]{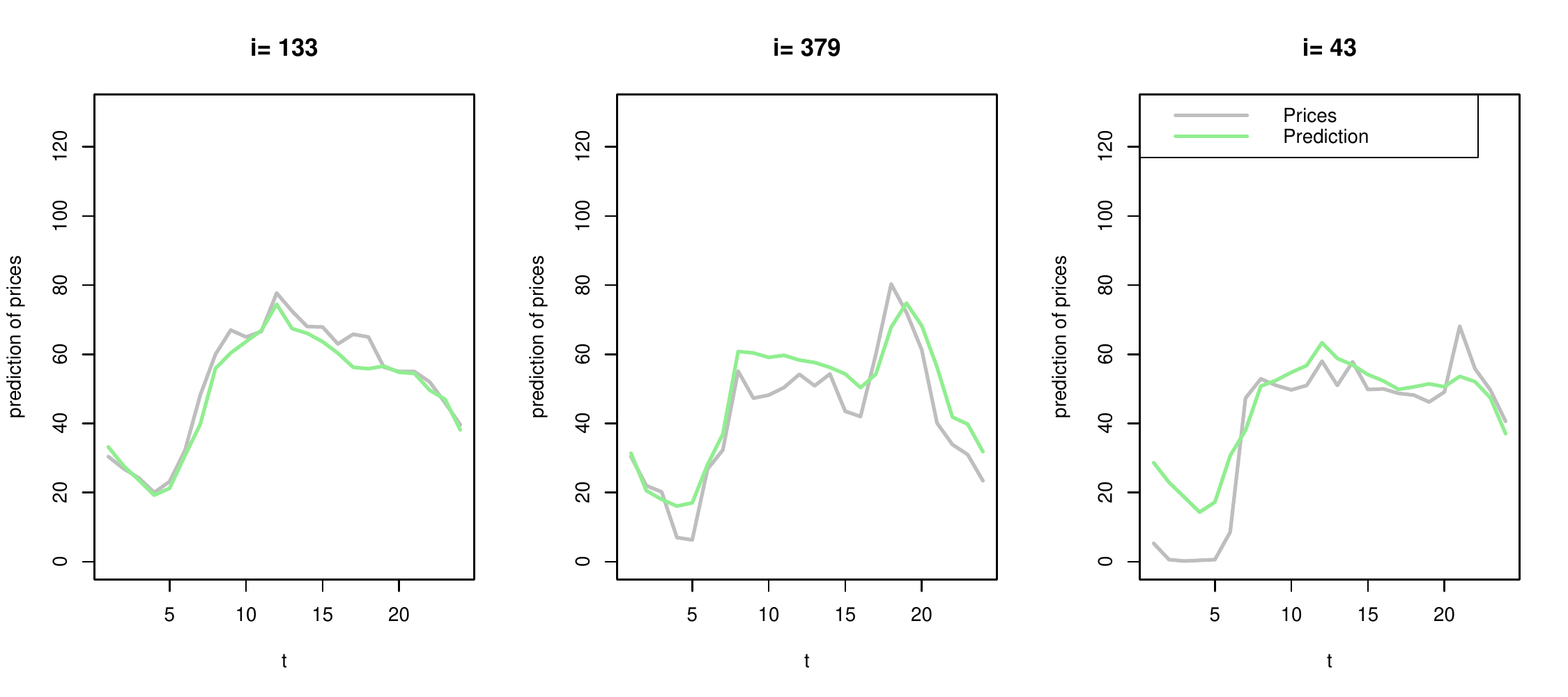}
\caption{\label{fig:WindPrices_data_pred} Cross-validated prediction $\widehat Y_i^{-i}$ made for three days taken randomly in the initial sample.}
\end{figure}

\smallskip
Similarly to what happen in the previous real data  problem (see  Figure \ref{fig:energy_data_pred} above), we see in Figure \ref{fig:WindPrices_data_pred} that the trends of each day is well captured and that the difficulty remains to predict the prices when there are  brutal changes of behavior in the curves.

\pagebreak
\section{Proofs}\label{sec:proofs}

All along the proofs, we denote by $C$ a positive constant which does not depend on $S$, $n$ or $m$ and whose value may change from line to line. Many proofs are based on technical results from the perturbation theory of bounded operators. A short account of the preliminary results we need is given in Section \ref{Perturbation-Theory} below. For more details, the reader can refers to \cite{dunford1965linear} and \cite{kato2013perturbation}.

\subsection{Proof of the results of  Section \ref{Estimation}}

\subsubsection{Proof of Proposition \ref{pro:estim_MC}}\label{sec:proof_prop_MC}

The remark following Proposition \ref{pro:estim_MC} implies that it is equivalent to reason either on $\gamma_n$, or on $\gamma_n'$. We choose $\gamma_n'$. We have 
$$\arg\min_{\mathcal T\in V_{m_1,m_2}'}\gamma_n'(\mathcal T)=\arg\min_{\mathcal T\in V_{m_1,m_2}'}\gamma_{n,1}'(\mathcal T)+\arg\min_{\mathcal T\in V_{m_1,m_2}'}\gamma_{n,2}'(\mathcal T),$$ where for any $\mathcal T\in V_{m_1,m_2}'$, 
$$\gamma_{n,1}'(\mathcal T)=\frac{1}{n}\sum_{i=1}^n\left\|\int_0^1\mathcal T(s,\cdot)X_i(s)ds\right\|^2,\;\;\gamma_{n,2}'(\mathcal T)=-\frac{2}{n}\sum_{i=1}^n\langle Y_i,\int_0^1\mathcal T(s,\cdot)X_i(s)ds\rangle.$$
For any $\mathcal T\in V_{m_1,m_2}'$, there also exists a unique sequence $b=(b_{j,k})_{j=1,\ldots,m_1,\;k=1,\ldots m_2}\in\mathbb R^{m_1m_2}$ such that $$\mathcal{T}(s,t)=\sum_{j=1}^{m_1}\sum_{k=1}^{m_2}b_{j,k}\phi_j(s)\phi_k(t),$$ with $(s,t)\in [0,1]^2$. Thus, 
$$\min_{\mathcal T\in V_{m_1,m_2}'} \gamma_{n}'(\mathcal T)=\min_{b\in\mathbb R^{m_1m_2}}\widetilde\gamma_{n}(b)\mbox{ with }\widetilde\gamma_{n}'=\widetilde\gamma_{n,1}'+\widetilde\gamma_{n,2}',$$ and
$$\widetilde\gamma_{n,1}'(b)=\frac{1}{n}\sum_{i=1}^n\sum_{k=1}^{m_2}\left(\sum_{j=1}^{m_1}b_{j,k}\langle\phi_j,X_i\rangle\right)^2,\;\;\widetilde\gamma_{n,2}'(b)=-\frac{2}{n}\sum_{i=1}^n\sum_{j=1}^{m_1}\sum_{k=1}^{m_2}b_{j,k}\langle\phi_j,X_i\rangle\langle Y_i,\phi_k\rangle.$$ Thus, we look for a minimum of the function $\widetilde\gamma_n$. The functions $\widetilde\gamma_{n,\ell}$, $\ell=1,2$, are differentiable and for any $(j_0,k_0)\in\{1,\ldots,m_1\}\times\{1,\ldots,m_2\}$,
$$\frac{\partial\widetilde\gamma_{n,1}(b)}{\partial b_{j_0,k_0}}=2\sum_{j=1}^{m_1}b_{j,k_0}\langle \Gamma_n \phi_{j_0},\phi_{j}\rangle,\;\;\frac{\partial\widetilde\gamma_{n,2}(b)}{\partial b_{j_0,k_0}}=-\langle \Delta_n \phi_{j_0},\phi_{k_0}\rangle.$$
This leads to $\nabla(\widetilde\gamma)(b)=2Ab-2Y_{\phi}$, with $b=(b_{j,k})_{j,k}\in\mathbb R^{m_1m_2}$. We have proved that $b=A^{-1}Y_{\phi}$ is a critical point. Moreover, the Hessian matrix can be computed as follows: 
$$\frac{\partial\widetilde\gamma_n(b)}{\partial b_{\ell,r}\partial b_{j_0,k_0}}=\delta_{r,k_0}\langle \Gamma_n \phi_{j_0},\phi_{k_0}\rangle,$$ 
where $\delta_{r,k_0}=1$ only if $r=k_0$, $\delta_{r,k_0}=0$ otherwise. By considering the indices $(j,k)\in\{1,\ldots,m_1\}\times\{1,\ldots,m_2\}$ of a vector $b\in\mathbb R^{m_1m_2}$, in the order $b=(b_{1,1},\ldots,b_{m_1,1},b_{1,2},\ldots,b_{m_1,2},\ldots,\ldots,b_{m_1,m_2})$, we obtain that the Hessian matrix of $\widetilde\gamma_n$ on $b$ is a block diagonal matrix, with $m_2$ blocks equal to $2A$. Thus, its determinant is $2^{m_2}(\mbox{det}(A))^{m_2}$. Since $A$ is the Gram matrix of a symmetric bilinear form, if it is invertible, the critical point is a global minimum, which proves Proposition \ref{pro:estim_MC}.

\subsubsection{Proof of Equality \eqref{eq:egalite_proj}}\label{sec:proof_egalite_proj}

For any $T\in \mathcal{L}_2(\mathbb H)$, 
$$\Pi_{m_1,m_2}^{op}T=\sum_{j=1}^{m_1}\sum_{k=1}^{m_2}\langle T,\varphi_k\otimes\varphi_j\rangle_{HS}\varphi_k\otimes\varphi_j,$$
with $\langle\cdot,\cdot\rangle_{HS}$ the scalar product associated to the Hilbert-Schmidt norm. Then, for any $r\in\mathbb N\backslash\{0\}$, 
\begin{align*}\Pi_{m_1,m_2}^{op}T(\varphi_r)&=\sum_{j=1}^{m_1}\sum_{k=1}^{m_2}\langle T,\varphi_k\otimes\varphi_j\rangle_{HS}\varphi_k\otimes\varphi_j(\varphi_r),\\
&=\sum_{j=1}^{m_1}\sum_{k=1}^{m_2}\langle T,\varphi_k\otimes\varphi_j\rangle_{HS}\delta_{j,r}\varphi_k,\mbox{ since the basis is orthonormal,}\\
&=\delta_{r\leq m_1}\sum_{k=1}^{m_2}\langle T,\varphi_k\otimes\varphi_r\rangle_{HS}\varphi_k,\\
&=\delta_{r\leq m_1}\sum_{k=1}^{m_2}\sum_{\ell=1}^\infty\langle T\varphi_{\ell},\varphi_k\otimes\varphi_r(\varphi_{\ell})\rangle\varphi_k,\mbox{ by definition of the scalar product},\\
&=\delta_{r\leq m_1}\sum_{k=1}^{m_2}\langle T\varphi_r,\varphi_k\rangle\varphi_k=\delta_{r\leq m_1}\sum_{k=1}^{m_2}\varphi_k\otimes\varphi_k(T\varphi_r)=\delta_{r\leq m_1}\Pi_{m_2}T(\varphi_r),\\
&=\Pi_{m_2}T\Pi_{m_1}(\varphi_r).
\end{align*}
This ends the proof.

\subsubsection{Proof of Proposition \ref{pro:estim_MC_ACP}}\label{sec:proof_prop_MC_ACP}
Let us start with the proof of \eqref{eq:MC_ACP_kernel}. Considering the result of Proposition \ref{pro:estim_MC}, we begin with the computation of the elements of the matrix $A$, when we consider the PCA basis. We write
$$ \langle\Gamma_n\widehat\varphi_j,\widehat \varphi_k\rangle=\langle\widehat\lambda_j\widehat\varphi_j,\widehat\varphi_k\rangle=\delta_{j,k}\widehat\lambda_j,$$
with $\delta_{j,k}=0$ if $j\neq k$, $\delta_{j,k}=1$ otherwise. Thus, $A$ is a diagonal matrix, and in this case, if $\widehat\lambda_j>0$, for any $j=1,\ldots,m_1$, we obtain the existence and uniqueness of the least-squares estimator, since $A^{-1}$ exists and is equal to the diagonal matrix with diagonal elements $\widehat\lambda_j^{-1}$. The coefficients of the estimators are
$$\widehat b_{j,k}=\frac{1}{\widehat\lambda_j}\langle\Delta_n\widehat\varphi_j,\widehat\varphi_k\rangle,\;\;j\in\{1,\ldots,m_1\},\;k\in\{1,\ldots,m_2\}.$$ This ends the proof of \eqref{eq:MC_ACP_kernel}.

To prove \eqref{eq:MC_ACP_op}, we start from \eqref{eq:estim_MC_op}, and the previous expression for $\widehat{\mathcal S}_m$. We immediatly get 
\begin{align*}
\widehat S_{m_1,m_2}&=\sum_{j=1}^{m_1}\sum_{k=1}^{m_2}\frac{1}{\widehat\lambda_j}\langle\Delta_n\widehat\varphi_j,\widehat\varphi_k\rangle\widehat\varphi_k\otimes\widehat\varphi_j,\\
&=\sum_{j=1}^{m_1}\frac{1}{\widehat\lambda_j}\langle\Delta_n\widehat\varphi_j,\widehat\varphi_j\rangle\widehat\varphi_j\otimes\widehat\varphi_j+\sum_{\substack{1\leq j\leq m_1 \\1\leq k\leq m_2 \\ j\neq k}}\frac{1}{\widehat\lambda_j}\langle\Delta_n\widehat\varphi_j,\widehat\varphi_k\rangle\widehat\varphi_k\otimes\widehat\varphi_j.
\end{align*}
It remains to apply the following lemma to the operator $T=\Delta_n\Gamma_{n,m_1}^{\dag} $ of the first part of the right-hand-side of the last equality, and to remark that $\Gamma_{n,m_1}^{\dag} \widehat\varphi_j=0$ as soon as $j\geq m$. 

\begin{lemma}\label{lm:op_base} Let $T$ be a linear operator of a separable Hilbert space $(\mathbb H,\langle\cdot,\cdot\rangle)$, self-adjoint and compact.  Let $(e_j)_{j\geq 1}$ be an orthonormal basis of eigenvectors of $T$. Then, 
$$T=\sum_{j=1}^\infty\langle Te_j,e_j\rangle e_j\otimes e_j.$$
\end{lemma}

\noindent\textit{Proof of Lemma \ref{lm:op_base}}

For  any Hilbertian basis $(e_j)_{j\geq 1}$, it is well known that any operator $T$ can be written 
$$T=\sum_{j,k=1}^\infty\langle Te_j,e_k\rangle e_k\otimes e_j.$$
Now, if the $e_j$'s are orthonormal eigenfunctions of $T$, there exist eigenvalues $\theta_j$ such that
$$\langle Te_j,e_k\rangle=\theta_j\langle e_j,e_k\rangle=\theta_j\delta_{j,k},$$
with $\delta_{j,k}=1$ if $j=k$, $\delta_{j,k}=0$ otherwise. This proves the result.

\subsection{Perturbation theory background}
\label{Perturbation-Theory}

We provide here a digest of some key results of the perturbation theory, which will be helpful within the proofs. These results are largely sourced from \cite{roche2014modelisation,brunel2016non,MR15}, but also from \cite{crambes2013asymptotics} and \cite{cardot2007clt}. In a nutshell, the aim of the perturbation theory is to control the proximity between the eigenfunctions of $\Gamma$ and those of the random operator $\Gamma_n$. We remind that the operator $\Pi_m$ (resp. $\widehat{\Pi}_m$) stands for the orthonormal projector onto $\Span (\varphi_1, \ldots, \varphi_m)$ (resp. $\Span (\widehat{\varphi}_1, \ldots, \widehat{\varphi}_m)$).

Let us denote by $\mathcal{B}_j$ the oriented circle of the complex plane of center $\lambda_j$ and radius $\delta_j/2$, 
where $\delta_j  = \lambda_j - \lambda_{j+1}.$ We also define $\mathcal{C}_m = \bigcup^m_{j=1} \mathcal{B}_j$ which is a union of disjoint circles since, by Assumption~ $\boldsymbol{\mathcal{A}_{2}}$, we also have $\delta_j=\min \{ \lambda_j - \lambda_{j+1}, \lambda_{j-1} - \lambda_j \}$. Let also $\displaystyle \bm{a_j} = \frac{\lambda_j}{\delta_j} + \sum_{r \neq j} \frac{\lambda_r}{\vert \lambda_r - \lambda_j \vert},$ for all $j \geq 1$, we define the set 
\begin{equation*}
\mathcal{A}_n = \bigcap^m_{j=1} \left\{ \vert \widehat{\lambda}_j - \lambda_j \vert < \frac{\delta_j}{2} \right\} \bigcap \left\{ \sup_{z \in \supp(\mathcal{C}_m)} \Vert T_n (z) \Vert_{\infty} < \frac{\mathbf{a}_j}{\sqrt{n}} \ln(n)\right\}. 
\end{equation*}

The following lemma is the keystone of the results related to perturbation theory. It provides a link between the difference of the empirical and theoretical projectors $\widehat{\Pi}_m - \Pi_m$, which we want to control, and the difference between empirical and theoretical covariance operators $\Gamma_n - \Gamma$, which can be controlled with the Bernstein inequality. 
\begin{lemma}
\label{lemma-hatPim-Pim}
Under Assumption $\boldsymbol{\mathcal{A}_{5}}$, there exists a set $\mathcal{A}_n$ such that 
$$\mathbb{P} \left( \mathcal{A}_n^{\complement} \right) \leq \exp(-c^* \ln(n)^2),$$ 
where $c^*$ is some positive constant depending on $(\lambda_j)_{j \geq 1}$ and  
\begin{equation*}
(\widehat{\Pi}_m - \Pi_m ) \bm{1}_{\mathcal{A}_n} = \frac{1}{2 i \pi} \sum_{k = 1}^m \int_{\mathcal{B}_k} R(z) (\Gamma_n - \Gamma) R(z) \mathrm{d} z \bm{1}_{\mathcal{A}_n} + \frac{1}{2 i \pi} \sum_{k = 1}^m \int_{\mathcal{B}_k} R^{1/2}(z) [I - T_n(z)]^{-1} T_n(z)^2 R^{1/2}(z) \mathrm{d} z \bm{1}_{\mathcal{A}_n},
\end{equation*}
with $T_n (z) = R^{1/2} (z) (\Gamma_n - \Gamma) R^{1/2} (z)$ and $R(z) = (z I - \Gamma)^{-1}$.
\end{lemma}
Lemma~\ref{lemma-hatPim-Pim} is proved in \cite{brunel2016non} (see Lemma 12 and Remark 4, p.224).


Throughout the proofs, we will also need some results on the behavior of the eigenvalues $(\lambda_j)_{j \geq 1}$. Lemma \ref{lemma-behavior-vap} and Lemma \ref{lemma-behavior-ak} the main results we will need. 

\begin{lemma} 
\label{lemma-behavior-vap}\citep[Lemma 1]{cardot2007clt}
Assume that Assumption $\boldsymbol{\mathcal{A}_{3}}$ is satisfied. Then, for all positive integers $j$ and $k$, such that $k > j$, we have 
$$j \lambda_j \geq k \lambda_k \quad \mbox{ and } \quad \lambda_j - \lambda_k \geq \left( 1 - \frac{j}{k} \right) \lambda_j.$$
In addition, 
$$\sum_{r \geq k} \lambda_r \leq (k+1) \lambda_k.$$
\end{lemma} 
 
\begin{lemma}\citep[Lemma 10.1]{hilgert2013minimax}
\label{lemma-behavior-ak}
Under Assumptions $\boldsymbol{\mathcal{A}_{4}}$ and $\boldsymbol{\mathcal{A}_{6}}$, we have
$$\mathbf{a}_k \leq C k \ln(k).$$
\end{lemma}

\subsection{Proofs of the results of Section \ref{Risque&Optimality}}\label{sec:proof_bv}

In order to achieve the bias-variance decomposition for the prediction risk presented in Section \ref{Risque&Optimality}, we will need to formulate $\widehat S_{m_1,m_2}$ in terms of $S$. We remind that 
\begin{equation}
\label{eq-Sn}
\widehat S_{m_1,m_2} = \widehat{\Pi}_{m_2} \Delta_n \Gamma^{\dagger}_{n,m_1}.
\end{equation}
Due to the linearity of $S$, it is straightforward that
$$\Delta_n = S \Gamma_n + \frac{1}{n} \sum_{i = 1}^n \varepsilon_i \otimes X_i.$$
Now, given that $\Gamma_{n,m_1}^{\dagger} $ is self-adjoint, one can easily see that
\begin{equation}
\label{eq_Delta_Gamma_dagger}
\Delta_n \Gamma_{n,m_1}^{\dagger}  = S \widehat{\Pi}_{m_1} + U_n,
\end{equation}
where $U_n =n^{-1} \sum_{i = 1}^n \varepsilon_i \otimes \Gamma_{n,m_1}^{\dagger}  (X_i)$. Combining Equations \eqref{eq-Sn} and \eqref{eq_Delta_Gamma_dagger} allows us to write
\begin{equation}
\label{eq-Sn-wrt-S}
\widehat S_{m_1,m_2} = \widehat{\Pi}_{m_2} S \widehat{\Pi}_{m_1} + \widehat{\Pi}_{m_2} U_n.
\end{equation}
Besides, Lemma \ref{lemma-ps-Xn+1-HS} stated just below will be very helpful for onward proofs.   

\begin{lemma}
\label{lemma-ps-Xn+1-HS}
Let $V$ and $W$ be random bounded linear operators, independent of $X_{n+1}$. We assume that $V$ and $W$ are Hilbert-Schmidt operators. Then, 
$$\mathbb{E} \langle V(X_{n+1}) , W(X_{n+1}) \rangle = \mathbb{E} \langle V \Gamma^{1/2} , W \Gamma^{1/2} \rangle_{\HS},$$
where $\langle \cdot , \cdot \rangle_{\HS}$ refers to the Hilbert-Schmidt scalar product and
\[
\mathbb{E}\|V(X_{n+1})\|^2=\mathbb E\|V\Gamma^{1/2}\|_{\HS}^2. 
\]
\end{lemma}

\medskip

\noindent\textit{Proof of Lemma \ref{lemma-ps-Xn+1-HS}.} We start by
\begin{align*}
\mathbb{E} \langle V(X_{n+1}) , W(X_{n+1}) \rangle &= \sum_{j = 1}^{+\infty} \mathbb{E} \langle X_{n+1} , \varphi_j \rangle \langle V(\varphi_j) , W(X_{n+1}) \rangle \\
&= \sum_{j = 1}^{+ \infty} \mathbb{E} \left\langle V(\varphi_j) , W \left(\langle X_{n+1} , \varphi_j \rangle X_{n+1} \right) \right\rangle.
\end{align*}
From here, we first compute the expectation with respect to $X_{n+1}$. Given both the linearity of $W$ and that $X_{n+1}$ is independent of $V$ and $W$, we obtain
\begin{align*}
\mathbb{E} \langle V(X_{n+1}) , W(X_{n+1}) \rangle &= \sum_{j = 1}^{+ \infty} \mathbb{E} \left\langle V(\varphi_j) , W \Gamma (\varphi_j) \right\rangle \\
&= \sum_{j = 1}^{+ \infty} \lambda_j \mathbb{E} \left\langle V(\varphi_j) , W(\varphi_j) \right\rangle.
\end{align*}
Therefore, 
\begin{align*}
\mathbb{E} \langle V(X_{n+1}) , W(X_{n+1}) \rangle &= \sum_{j = 1}^{+ \infty} \mathbb{E} \left\langle V \Gamma^{1/2} (\varphi_j) , W \Gamma^{1/2} (\varphi_j) \right\rangle \\
&= \mathbb{E} \langle V \Gamma^{1/2} , W \Gamma^{1/2} \rangle_{\HS}.
\end{align*}
This proves the first equality of Lemma \ref{lemma-ps-Xn+1-HS}, and the second one is a direct consequence, taking $V=W$.

\subsubsection{Proof of Theorem \ref{thm-upper-bound-MSPE}}

The starting point of the proof is to achieve a first bias-variance decomposition of the prediction risk. As a result of Lemma \ref{lemma-ps-Xn+1-HS}, we note
\begin{align*}
\mathbb{E} \Vert \widehat S_{m_1,m_2} (X_{n+1}) - S (X_{n+1}) \Vert^2 &= \mathbb{E} \Vert (\widehat S_{m_1,m_2} - S) \Gamma^{1/2} \Vert_{\HS}^2. 
\end{align*}
In addition, from Equation \eqref{eq-Sn-wrt-S}, we write 
\begin{align}
\label{eq-first-decomp-MSPE}
\mathbb{E} \Vert \widehat S_{m_1,m_2} (X_{n+1}) - S (X_{n+1}) \Vert^2 &= \mathbb{E} \Vert ( \widehat{\Pi}_{m_2} S \widehat{\Pi}_{m_1} - S) \Gamma^{1/2} + \widehat{\Pi}_{m_2} U_n \Gamma^{1/2} \Vert_{\HS}^2 \nonumber \\
&= \mathbb{E} \Vert (S - \widehat{\Pi}_{m_2} S \widehat{\Pi}_{m_1}) \Gamma^{1/2} \Vert_{\HS}^2 + \mathbb{E} \Vert \widehat{\Pi}_{m_2} U_n \Gamma^{1/2} \Vert_{\HS}^2 \nonumber \\
&- 2 \mathbb{E} \langle (S - \widehat{\Pi}_{m_2} S \widehat{\Pi}_{m_1}) \Gamma^{1/2} , \widehat{\Pi}_{m_2} U_n \Gamma^{1/2} \rangle_{\HS}.
\end{align}
The last right-side expectation of Equation \eqref{eq-first-decomp-MSPE} is null. Indeed, recall that
\begin{equation*}
\mathbb{E} \langle (S - \widehat{\Pi}_{m_2} S \widehat{\Pi}_{m_1}) \Gamma^{1/2} , \widehat{\Pi}_{m_2} U_n \Gamma^{1/2} \rangle_{\HS} = \sum_{j = 1}^{+ \infty} \lambda_j \mathbb{E} \langle (S - \widehat{\Pi}_{m_2} S \widehat{\Pi}_{m_1}) (\varphi_j) , \widehat{\Pi}_{m_2} U_n (\varphi_j) \rangle.
\end{equation*}
Notice that, $\widehat{\Pi}_{m_2} S \widehat{\Pi}_{m_1}$ only depends on $(X_1, \ldots, X_n)$. It remains then to show that the expectation of $\widehat{\Pi}_{m_2} U_n (\varphi_j)$ conditionally to $(X_1, \ldots, X_n)$ is zero. We write, 
\begin{align*}
\mathbb{E} [ \widehat{\Pi}_{m_2} U_n (\varphi_j) \vert X_1, \ldots, X_n ] &= \frac{1}{n} \sum_{i=1}^n  \mathbb{E} [ \langle \Gamma^{\dagger}_{n,m_1} (X_i) , \varphi_j\rangle \widehat{\Pi}_{m_2} (\varepsilon_i) \vert X_1, \ldots, X_n ] \\
&= \frac{1}{n} \sum_{i=1}^n \sum_{r = 1}^{m_2} \langle \Gamma^{\dagger}_{n,m_1} (X_i) ,\varphi_j \rangle  \mathbb{E} [ \langle \varepsilon_i , \varphi_r \rangle \varphi_r \vert X_1, \ldots, X_n ] \\
 &= \frac{1}{n} \sum_{i=1}^n \sum_{r = 1}^{m_2} \langle \Gamma^{\dagger}_{n,m_1} (X_i) , \varphi_j  \rangle \langle \mathbb{E} [ \varepsilon_i ] , \varphi_r \rangle \varphi_r \\
&= 0_{\mathbb{H}}.
\end{align*}   
In other words, 
$$\mathbb{E} \langle (S - \widehat{\Pi}_{m_2} S \widehat{\Pi}_{m_1}) \Gamma^{1/2} , \widehat{\Pi}_{m_2} U_n \Gamma^{1/2} \rangle_{\HS} = 0.$$
Therefore, 
\begin{equation}
\label{eq-gross-decomp-MSPE}
\mathbb{E} \Vert \widehat S_{m_1,m_2} (X_{n+1}) - S (X_{n+1}) \Vert^2 = \mathbb{E} \Vert (S - \widehat{\Pi}_{m_2} S \widehat{\Pi}_{m_1}) \Gamma^{1/2} \Vert_{\HS}^2 + \mathbb{E} \Vert \widehat{\Pi}_{m_2} U_n \Gamma^{1/2} \Vert_{\HS}^2.
\end{equation}
In Equation \eqref{eq-gross-decomp-MSPE}, we recognize the common decomposition of the estimation risk as a compromise of a bias and variance terms. The bias, namely the first term in the right side, is decreasing with respect to $m_1$ and $m_2$ and is related to the regularity of the regression model $S$. While the variance term is increasing with respect to $m_2$ and depending in particular on the observation errors $\varepsilon_1, \ldots, \varepsilon_n$. In the following propositions \ref{prop-variance-MSPE} and \ref{prop-bias-MSPE}, we give sharp upper bounds of these two terms, which ends the proof.

\begin{proposition}
\label{prop-variance-MSPE}
The variance term of Equation \eqref{eq-gross-decomp-MSPE} can be upper bounded as $$\mathbb{E} \Vert \widehat{\Pi}_{m_2} U_n \Gamma^{1/2} \Vert_{\HS}^2 \leq \sigma_{\varepsilon}^2 \frac{m_1}{n} + A_{n,m_1},$$
where $\sigma_{\varepsilon} = \mathbb{E} \Vert \varepsilon \Vert^2$ and $\displaystyle A_{n,m_1} = \sigma_{\varepsilon}^2 \frac{C m_1^2 \ln^2 (m_1)}{n^2}$.
\end{proposition}

\begin{proposition}
\label{prop-bias-MSPE}
The bias term of Equation \eqref{eq-gross-decomp-MSPE} can be upper bounded as 
\begin{equation*}
\mathbb{E} \Vert (S - \widehat{\Pi}_{m_2} S \widehat{\Pi}_{m_1}) \Gamma^{1/2} \Vert_{\HS}^2 \leq 3 \sum_{j = m_1 + 1}^{+\infty} \lambda_j \Vert S (\varphi_j) \Vert^2 + 3 \sum_{j = 1}^{m_1} \lambda_j \Vert (\Id - \Pi_{m_2}) S (\varphi_j) \Vert^2 + B_{n,m_1} + D_{n,m_2} + E_n. 
\end{equation*}
where 
\begin{align*}
B_{n,m_1} &= \frac{C m_1^2 \lambda_{m_1} \Vert S \Vert_{\HS}}{n}, \quad E_n = \frac{C}{n^2} \Vert S \Gamma^{1/2} \Vert_{\HS}^2 \\
D_{n,m_2} &=  \frac{C m_2 \lambda_{m_2}}{n} + \frac{C m_2^2 \ln(m_2)}{n \psi_{\beta} (m_2)} + \frac{C m_2^3}{n \psi_{\beta} (\lfloor m_2/2 \rfloor)} + \frac{C \ln^4 (n)}{n^2} \left( \sum_{k = 1}^{m_2} \frac{k^2 \ln^2 (k)}{\sqrt{\psi_{\beta} (k)}} \right)^2.
\end{align*}
\end{proposition}

\medskip

\noindent\textit{Proof of Proposition \ref{prop-variance-MSPE}.} It is obvious that
\begin{align*}
\mathbb{E} \Vert \widehat{\Pi}_{m_2} U_n \Gamma^{1/2} \Vert_{\HS}^2 &= \sum_{j = 1}^{+\infty} \sum_{r = 1}^{m_2} \lambda_j \mathbb{E} \langle U_n (\varphi_j) , \varphi_r \rangle^2 \\
&= \frac{1}{n^2} \sum_{j = 1}^{+\infty} \sum_{r = 1}^{m_2} \lambda_j \mathbb{E} \left[ \sum_{i = 1}^n \langle \Gamma^{\dagger}_{n,m_1} (X_i) , \varphi_j \rangle \langle \varepsilon_i , \varphi_r \rangle \right]^2.
\end{align*}
Then, 
\begin{align*}
\mathbb{E} \Vert \widehat{\Pi}_{m_2} U_n \Gamma^{1/2} \Vert_{\HS}^2 &= \frac{1}{n^2} \sum_{j = 1}^{+\infty} \sum_{r = 1}^{m_2} \sum_{i = 1}^n \lambda_j \mathbb{E} \langle \Gamma^{\dagger}_{n,m_1} (X_i) , \varphi_j \rangle^2 \langle \varepsilon_i , \varphi_r \rangle^2 \\
&+ \frac{1}{n^2} \sum_{j = 1}^{+\infty} \sum_{r = 1}^{m_2} \sum_{\substack{i \neq i' \\ (i,i') \in [\![1, n]\!]^2}} \lambda_j \mathbb{E} \langle \Gamma^{\dagger}_{n,m_1} (X_i) , \varphi_j \rangle \langle \Gamma^{\dagger}_{n,m_1} (X_{i'}) , \varphi_j \rangle \langle \varepsilon_i , \varphi_r \rangle \langle \varepsilon_{i'} , \varphi_r \rangle.
\end{align*}
All the terms of the second right-side sum are null. Indeed, by independence of $(\varepsilon_1, \ldots, \varepsilon_n)$ and $(X_1, \ldots, X_n)$, we have for $i \neq i'$, 
\begin{multline*}
\mathbb{E} \left[ \langle \Gamma^{\dagger}_{n,m_1} (X_i) , \varphi_j \rangle \langle \Gamma^{\dagger}_{n,m_1} (X_{i'}) , \varphi_j \rangle \langle \varepsilon_i , \varphi_r \rangle \langle \varepsilon_{i'} , \varphi_r \rangle \vert X_1, \ldots, X_n \right] \\
= \langle \Gamma^{\dagger}_{n,m_1} (X_i) , \varphi_j \rangle \langle \Gamma^{\dagger}_{n,m_1} (X_{i'}) , \varphi_j \rangle \mathbb{E} \left[ \langle \varepsilon_i , \varphi_r \rangle \langle \varepsilon_{i'} , \varphi_r \rangle \vert X_1, \ldots, X_n \right].
\end{multline*}
Thus, knowing that $\varepsilon_1, \ldots, \varepsilon_n$ are independent and centered, we have
\begin{multline*}
\mathbb{E} \left[ \langle \Gamma^{\dagger}_{n,m_1} (X_i) , \varphi_j \rangle \langle \Gamma^{\dagger}_{n,m_1} (X_{i'}) , \varphi_j \rangle \langle \varepsilon_i , \varphi_r \rangle \langle \varepsilon_{i'} , \varphi_r \rangle \vert X_1, \ldots, X_n \right] \\
= \langle \Gamma^{\dagger}_{n,m_1} (X_i) , \varphi_j \rangle \langle \Gamma^{\dagger}_{n,m_1} (X_{i'}) , \varphi_j \rangle \langle \mathbb{E} [\varepsilon_i] , \varphi_r \rangle \langle \mathbb{E} [\varepsilon_{i'}] , \varphi_r \rangle  = 0.
\end{multline*}
Therefore, the variance term is simply given by 
\begin{align*}
\mathbb{E} \Vert \widehat{\Pi}_{m_2} U_n \Gamma^{1/2} \Vert_{\HS}^2 &= \frac{1}{n^2} \sum_{i = 1}^n \sum_{j = 1}^{+\infty} \lambda_j \mathbb{E} \langle \Gamma^{\dagger}_{n,m_1} (X_i) , \varphi_j \rangle^2 \mathbb{E} \Vert \Pi_{m_2} (\varepsilon_i) \Vert^2 \\
&\leq \frac{\sigma_{\varepsilon}^2}{n^2} \sum_{i = 1}^n \sum_{j = 1}^{+\infty} \lambda_j \mathbb{E} \langle \Gamma^{\dagger}_{n,m_1} (X_i) , \varphi_j \rangle^2,
\end{align*}
where $\sigma_{\varepsilon}^2 = \mathbb{E} \Vert \varepsilon \Vert^2$. Furthermore, giving that $\Gamma_{n,m_1}^{\dagger} $ is self-adjoint, we get
\begin{align*}
\frac{1}{n} \sum_{i = 1}^n \sum_{j = 1}^{+\infty} \lambda_j \mathbb{E} \langle \Gamma^{\dagger}_{n,m_1} (X_i) , \varphi_j \rangle^2 &= \frac{1}{n} \sum_{i = 1}^n \sum_{j = 1}^{+\infty} \lambda_j \mathbb{E} \langle X_i , \Gamma^{\dagger}_{n,m_1} (\varphi_j) \rangle^2 \\
&= \frac{1}{n} \sum_{i = 1}^n \sum_{j = 1}^{+\infty} \lambda_j \mathbb{E} \left\langle \langle X_i , \Gamma^{\dagger}_{n,m_1} (\varphi_j) \rangle X_i , \Gamma^{\dagger}_{n,m_1} (\varphi_j) \right\rangle \\
&= \sum_{j = 1}^{+\infty} \lambda_j \mathbb{E} \langle \Gamma_n \Gamma^{\dagger}_{n,m_1} (\varphi_j) , \Gamma^{\dagger}_{n,m_1} (\varphi_j) \rangle.
\end{align*}
Yet, it is easy to see that $\Gamma_{n,m_1}^{\dagger}  \Gamma_n \Gamma^{\dagger}_{n,m_1} = \Gamma_{n,m_1}^{\dagger}  \widehat{\Pi}_{m_1} = \Gamma_{n,m_1}^{\dagger} $. Thus, 
\begin{equation*}
\mathbb{E} \Vert \widehat{\Pi}_{m_2} U_n \Gamma^{1/2} \Vert_{\HS}^2 \leq \frac{\sigma_{\varepsilon}^2}{n} \sum_{j = 1}^{+\infty} \lambda_j \mathbb{E} \langle \Gamma^{\dagger}_{n,m_1} (\varphi_j) , \varphi_j \rangle.
\end{equation*}
Now, it is fairly easy to show the following equation. This can be done by diagonalization of the self-adjoint operator $\Gamma_{n,m_1}^{\dagger} $ in the orthonormal basis of its eigenfunctions, we write 
\begin{equation*}
\langle \Gamma^{\dagger}_{n,m_1} (\varphi_j) , \varphi_j \rangle = \Tr \left( \Gamma^{\dagger}_{n,m_1} \cdot \varphi_j \otimes \varphi_j \right), 
\end{equation*}
where the notation '$\cdot$' refers to the Hadamard product for operators (i.e. for two operators $S,T\in \mathcal L(H)$, $S\cdot T = \sum_{j,k\geq 1}\langle S,\varphi_j\rangle\langle T,\varphi_k\rangle\varphi_j\otimes\varphi_k$). Then, 
\begin{align*}
\mathbb{E} \Vert \widehat{\Pi}_{m_2} U_n \Gamma^{1/2} \Vert_{\HS}^2 &\leq  \frac{\sigma_{\varepsilon}^2}{n} \mathbb{E} \left[ \Tr \left( \Gamma^{\dagger}_{n,m_1} \cdot \Gamma \right) \right] \\
&= \frac{\sigma_{\varepsilon}^2}{n} \left[ \Tr \left( \mathbb{E} \left[ ( \Gamma^{\dagger}_{n,m_1} - \Gamma^{\dagger} ) \cdot \Gamma \right] \right) + \Tr \left( \Gamma^{\dagger} \cdot \Gamma \right) \right].
\end{align*}
Moreover, based on \citet[Lemma 19]{crambes2013asymptotics}, we have
\begin{equation*}
\Tr \left( \mathbb{E} \left[ ( \Gamma^{\dagger}_{n,m_1} - \Gamma^{\dagger} ) \cdot \Gamma \right] \right) \leq \frac{C m_1^2 \ln^2 (m_1)}{n}.
\end{equation*}
At last, remark that $\Tr \left( \Gamma^{\dagger} \cdot \Gamma \right) = \Tr \left( \Pi_{m_1} \right)=m_1$, which finishes the proof. 

\medskip 

\noindent\textit{Proof of Proposition \ref{prop-bias-MSPE}.} To achieve a sharp upper bound for the bias term, we intensively rely on the perturbation theory for bounded operators presented in Section \ref{Perturbation-Theory}. We also use some results already established in \cite{crambes2013asymptotics}. We begin with the following plain decomposition, 
\begin{equation*}
S - \widehat{\Pi}_{m_2} S \widehat{\Pi}_{m_1} = S - \Pi_{m_2} S \Pi_{m_1} - \widehat{\Pi}_{m_2} S (\widehat{\Pi}_{m_1} - \Pi_{m_1}) - (\widehat{\Pi}_{m_2} - \Pi_{m_2}) S \Pi_{m_1}. 
\end{equation*}
Then, 
\begin{align}
\label{eq-first-decomp-bias}
\mathbb{E} \Vert (S - \widehat{\Pi}_{m_2} S \widehat{\Pi}_{m_1}) \Gamma^{1/2} \Vert_{\HS}^2 &\leq 3 \mathbb{E} \Vert (S - \Pi_{m_2} S \Pi_{m_1}) \Gamma^{1/2} \Vert_{\HS}^2 + 3 \mathbb{E} \Vert \widehat{\Pi}_{m_2} S (\widehat{\Pi}_{m_1} - \Pi_{m_1}) \Gamma^{1/2} \Vert_{\HS}^2 \nonumber \\
&+ 3 \mathbb{E} \Vert (\widehat{\Pi}_{m_2} - \Pi_{m_2}) S \Pi_{m_1} \Gamma^{1/2} \Vert_{\HS}^2. 
\end{align}
In the remainder of the proof, we upper bound each term of Equation \eqref{eq-first-decomp-bias}. We have
\begin{align*}
\mathbb{E} \Vert (S - \Pi_{m_2} S \Pi_{m_1}) \Gamma^{1/2} \Vert_{\HS}^2 &= \mathbb{E} \Vert S (\Id - \Pi_{m_1}) \Gamma^{1/2} + (\Id - \Pi_{m_2}) S \Pi_{m_1} \Gamma^{1/2} \Vert_{\HS}^2 \\
&= \mathbb{E} \Vert S (\Id - \Pi_{m_1}) \Gamma^{1/2} \Vert_{\HS}^2 + \mathbb{E} \Vert (\Id - \Pi_{m_2}) S \Pi_{m_1} \Gamma^{1/2} \Vert_{\HS}^2 \\
&+ 2 \mathbb{E} \langle S (\Id - \Pi_{m_1}) \Gamma^{1/2} , (\Id - \Pi_{m_2}) S \Pi_{m_1} \Gamma^{1/2} \rangle_{\HS}. 
\end{align*}
The last right-side term of the equation just above is zero. Indeed,  
\begin{align*}
\langle S (\Id - \Pi_{m_1}) \Gamma^{1/2} , (\Id - \Pi_{m_2}) S \Pi_{m_1} \Gamma^{1/2} \rangle_{\HS} &= \sum_{j = 1}^{+\infty} \lambda_j \langle S (\Id - \Pi_{m_1}) (\varphi_j) , (\Id - \Pi_{m_2}) S \Pi_{m_1} (\varphi_j) \rangle \\
&= 0.
\end{align*}
Besides, 
\begin{align*}
\mathbb{E} \Vert S (\Id - \Pi_{m_1}) \Gamma^{1/2} \Vert_{\HS}^2 &= \sum_{j = 1}^{+\infty} \lambda_j \Vert S (\Id - \Pi_{m_1}) (\varphi_j) \Vert^2 \\
&= \sum_{j = m_1 + 1}^{+\infty} \lambda_j \Vert S (\varphi_j) \Vert^2. 
\end{align*}
In addition, 
\begin{align*}
\mathbb{E} \Vert (\Id - \Pi_{m_2}) S \Pi_{m_1} \Gamma^{1/2} \Vert_{\HS}^2 &= \sum_{j = 1}^{+ \infty} \lambda_j \Vert (\Id - \Pi_{m_2}) S \Pi_{m_1} (\varphi_j) \Vert^2 \\
&= \sum_{j = 1}^{m_1} \lambda_j \Vert (\Id - \Pi_{m_2}) S (\varphi_j) \Vert^2.
\end{align*}
This means that
\begin{equation}
\label{eq-term1-bias}
\mathbb{E} \Vert (S - \Pi_{m_2} S \Pi_{m_1}) \Gamma^{1/2} \Vert_{\HS}^2 = \sum_{j = m_1 + 1}^{+\infty} \lambda_j \Vert S (\varphi_j) \Vert^2 + \sum_{j = 1}^{m_1} \lambda_j \Vert (\Id - \Pi_{m_2}) S (\varphi_j) \Vert^2.
\end{equation}
We now upper bound the second right-side term of Equation \eqref{eq-first-decomp-bias}. It is easy to see that
\begin{equation*}
\mathbb{E} \Vert \widehat{\Pi}_{m_2} S (\widehat{\Pi}_{m_1} - \Pi_{m_1}) \Gamma^{1/2} \Vert_{\HS}^2 \leq \mathbb{E} \Vert S (\widehat{\Pi}_{m_1} - \Pi_{m_1}) \Gamma^{1/2} \Vert_{\HS}^2.
\end{equation*}
An upper bound of the right expectation in the previous equation is given in \citet[Proposition 15]{crambes2013asymptotics}. However, we believe that following the article notations, a '$k$' is missing in the upper bound. This happens on the page 2644, in the equation block below Formula (17), while switching from the second to the third equation. We would like to point out that this error does not alter the optimality of the results demonstrated in the paper, but it constrains the regularity of $X$. We also believe that thanks to Assumption $\boldsymbol{\mathcal{A}_{6}}$, this loss can be avoided. This make it possible to replace the square of the sum by the sum of squares. We then obtain, 
\begin{equation}
\label{eq-term2-bias}
\mathbb{E} \Vert \widehat{\Pi}_{m_2} S (\widehat{\Pi}_{m_1} - \Pi_{m_1}) \Gamma^{1/2} \Vert_{\HS}^2 \leq \frac{C m_1^2 \lambda_{m_1} \Vert S \Vert_{\HS}}{n}.
\end{equation}
Upper-bounding the last term of \eqref{eq-first-decomp-bias} requires a wide use of perturbation theory. Lemma \ref{lemma-term3-bias} gives such an upper bound. The proof of Proposition~\ref{prop-bias-MSPE} is then a direct result of Equations \eqref{eq-first-decomp-bias}, \eqref{eq-term1-bias}, \eqref{eq-term2-bias} and Lemma \ref{lemma-term3-bias}.

\begin{lemma}
\label{lemma-term3-bias}
The following inequality holds
\begin{align*}
\mathbb{E} \Vert (\widehat{\Pi}_{m_2} - \Pi_{m_2}) S \Pi_{m_1} \Gamma^{1/2} \Vert_{\HS}^2 &\leq \frac{C m_2 \lambda_{m_2}}{n} + \frac{C m_2^2 \ln(m_2)}{n \psi_{\beta} (m_2)} + \frac{C m_2^3}{n \psi_{\beta} (\lfloor m_2/2 \rfloor)} + \frac{C \ln^2 (n)}{n^2} \left( \sum_{k = 1}^{m_2} \frac{k^2 \ln^2 (k)}{\sqrt{\psi_{\beta} (k)}} \right)^2 \\
&+ \frac{C}{n^2} \Vert S \Gamma^{1/2} \Vert_{\HS}^2.
\end{align*}
\end{lemma}

\medskip
\noindent\textit{Proof of Lemma \ref{lemma-term3-bias}.} According to Lemma \ref{lemma-hatPim-Pim}, we have
$$(\widehat{\Pi}_{m_2} - \Pi_{m_2}) \bm{1}_{\mathcal{A}_n} = \mathcal{H}_n (z) + \mathcal{G}_n (z),$$
where the operators $\mathcal{H}_n (z)$ and $\mathcal{G}_n (z)$ are defined as
\begin{align*}
\mathcal{H}_{n}(z) &= \frac{1}{2 i \pi} \sum_{k = 1}^{m_2} \int_{\mathcal{B}_k} R(z) (\Gamma_n - \Gamma) R(z) \mathrm{d} z \bm{1}_{\mathcal{A}_n}, \\
\mathcal{G}_n(z) &= \frac{1}{2 i \pi} \sum_{k = 1}^{m_2} \int_{\mathcal{B}_k} R^{1/2}(z) [I - T_n(z)]^{-1} T_n(z)^2 R^{1/2}(z) \mathrm{d} z \bm{1}_{\mathcal{A}_n},
\end{align*}
with $T_n (z) = R^{1/2} (z) (\Gamma_n - \Gamma) R^{1/2} (z)$ and $R(z) = (z I - \Gamma)^{-1}$. Then, 
\begin{align}
\label{eq-decomp-term3-Hn-Gn}
\mathbb{E} \Vert (\widehat{\Pi}_{m_2} - \Pi_{m_2}) S \Pi_{m_1} \Gamma^{1/2} \Vert_{\HS}^2 &\leq 2 \mathbb{E} \Vert \mathcal{H}_n (z) S \Pi_{m_1} \Gamma^{1/2} \Vert^2_{\HS} + 2 \mathbb{E} \Vert \mathcal{G}_n (z) S \Pi_{m_1} \Gamma^{1/2} \Vert^2_{\HS} \nonumber \\
&+ \mathbb{E} \Vert (\widehat{\Pi}_{m_2} - \Pi_{m_2}) S \Pi_{m_1} \Gamma^{1/2} \bm{1}_{\mathcal{A}_n^{\complement}} \Vert_{\HS}^2.
\end{align}
Upper bounding the expectations of Equation \eqref{eq-decomp-term3-Hn-Gn} is technical and the remainder of the proof is quite long. Since $\mathcal{H}_n (z)$ is self-adjoint, we write
\begin{align}
\label{eq-first-EHn(z)}
\mathbb{E} \Vert \mathcal{H}_n (z) S \Pi_{m_1} \Gamma^{1/2} \Vert^2_{\HS} &= \sum_{j = 1}^{m_1} \sum_{l = 1}^{+\infty} \lambda_j \mathbb{E} \langle \mathcal{H}_n (z) S (\varphi_j) , \varphi_l  \rangle^2 \nonumber \\
&= \sum_{j = 1}^{m_1} \sum_{l = 1}^{+\infty} \lambda_j \mathbb{E} \langle S (\varphi_j) , \mathcal{H}_n (z) (\varphi_l)  \rangle^2.
\end{align} 
Furthermore, for $j,l$ in $\mathbb{N}\backslash\{0\}$, we have
\begin{align*}
\langle S (\varphi_j) , \mathcal{H}_n (z) (\varphi_l) \rangle &= \sum_{r = 1}^{+\infty} \langle S (\varphi_j) , \varphi_r \rangle \langle \mathcal{H}_n (z) (\varphi_l) , \varphi_r \rangle \\
&= \sum_{r = 1}^{+\infty} \langle S (\varphi_j) , \varphi_r \rangle \times \frac{1}{2i\pi}  \sum_{k = 1}^{m_2} \int_{\mathcal{B}_k} \langle R(z) (\Gamma_n - \Gamma) R(z) (\varphi_l) , \varphi_r \rangle \mathrm{d} z \bm{1}_{\mathcal{A}_n}.
\end{align*}
Now, it is straightforward to show that $R(z) (\varphi_\ell) = (z - \lambda_\ell)^{-1} \varphi_\ell$. This together with the fact that $R(z)$ is self-adjoint, gives
\begin{align*}
\langle S (\varphi_j) , \mathcal{H}_n (z) (\varphi_l) \rangle &= \sum_{r = 1}^{+\infty} \langle S (\varphi_j) , \varphi_r \rangle \langle (\Gamma_n - \Gamma) (\varphi_l) , \varphi_r \rangle \times \frac{1}{2i\pi} \sum_{k = 1}^{m_2} \int_{\mathcal{B}_k} (z-\lambda_l)^{-1}(z-\lambda_r)^{-1} \mathrm{d} z \bm{1}_{\mathcal{A}_n} \\
&= \sum_{r = 1}^{+\infty} \langle S (\varphi_j) , \varphi_r \rangle \langle (\Gamma_n - \Gamma) (\varphi_l) , \varphi_r \rangle \mathcal{I}_{m_2,l,r} \bm{1}_{\mathcal{A}_n},
\end{align*}
where $\displaystyle \mathcal{I}_{m_2,l,r} = \frac{1}{2 i \pi} \sum_{k = 1}^{m_2} \int_{\mathcal{B}_k} (z-\lambda_l)^{-1}(z-\lambda_r)^{-1} \mathrm{d} z$. The integral sum $\displaystyle \mathcal{I}_{m_2,l,r}$ can be computed using Cauchy's integral formula \citep{rudin1987real}. Indeed, 
\begin{align*}
\mathcal{I}_{m_2,l,r} &= \frac{1}{\lambda_l - \lambda_r} \times \frac{1}{2i\pi} \sum_{k = 1}^{m_2} \int_{\mathcal{B}_k} \left[ \frac{1}{z - \lambda_l} - \frac{1}{z - \lambda_r} \right] \mathrm{d} z \\
&= \frac{1}{\lambda_l - \lambda_r} \sum_{k = 1}^{m_2} \left[ \Ind_{\mathcal{B}_k} (\lambda_l) - \Ind_{\mathcal{B}_k} (\lambda_r) \right],
\end{align*}
where $\Ind_{\mathcal{B}_k}$ is an integer-valued function, taking 1 if its argument belongs to the surface drawn by $\mathcal{B}_k$ and 0 otherwise. Then, 
\begin{equation*}
\mathcal{I}_{m_2,l,r} = \frac{1}{\lambda_l - \lambda_r} \left[ \bm{1}_{l \in \{1, \ldots, m_2 \}} - \bm{1}_{r \in \{1, \ldots, m_2 \}}  \right].
\end{equation*}
We turn back to Equation \eqref{eq-first-EHn(z)}. Since $\langle\Gamma\varphi_j,\varphi_r\rangle=0$ if $j\neq r$,  we can write
\begin{equation*}
\mathbb{E} \Vert \mathcal{H}_n (z) S \Pi_{m_1} \Gamma^{1/2} \Vert^2_{\HS} = \mathbb{E} \left[ A \bm{1}_{\mathcal{A}_n} \right] + \mathbb{E} \left[ B \bm{1}_{\mathcal{A}_n} \right],
\end{equation*}
where $A$ and $B$ are defined as
\begin{align*}
A &= \sum_{j = 1}^{m_1} \lambda_j \sum_{l = 1}^{m_2} \left[ \sum_{r = m_2 + 1}^{+ \infty} \frac{\left\langle \Gamma_n (\varphi_l), \varphi_r \right\rangle}{\lambda_l - \lambda_r} \left\langle  S (\varphi_j) , \varphi_r \right\rangle \right]^2, \\
B &= \sum_{j = 1}^{m_1} \lambda_j \sum_{l = m_2 + 1}^{+\infty} \left[ \sum_{r = 1}^{m_2} \frac{\left\langle \Gamma_n (\varphi_l), \varphi_r \right\rangle}{\lambda_r - \lambda_l} \left\langle S (\varphi_j) , \varphi_r \right\rangle \right]^2.
\end{align*}
This implies that
\begin{equation}
\label{eq-EHn(z)-EA+EB}
\mathbb{E} \Vert \mathcal{H}_n (z) S \Pi_{m_1} \Gamma^{1/2} \Vert^2_{\HS} \leq \mathbb{E} [ A ] + \mathbb{E} [ B ]. 
\end{equation}
Thanks to Assumption $\boldsymbol{\mathcal{A}_{6}}$, while developing the sum squares in $A$ and $B$, the expectation of the cross terms equals zero. In fact, for $r, r'$ in $\mathbb{N}\backslash\{0\}$ and $\ell$ in $\mathbb{N}\backslash\{0,r , r'\}$, we write 
\begin{align*}
\mathbb{E} \langle \Gamma_n (\varphi_l), \varphi_r \rangle \langle \Gamma_n (\varphi_l), \varphi_{r'} \rangle &= \frac{1}{n^2} \sum_{i = 1}^n \sum_{i' = 1}^n \mathbb{E} \langle X_i \otimes X_i (\varphi_l) , \varphi_r \rangle \langle X_{i'} \otimes X_{i'} (\varphi_l) , \varphi_{r'} \rangle \\
&= \frac{1}{n^2} \sum_{i = 1}^n \mathbb{E} \langle X_i , \varphi_r \rangle \langle X_i , \varphi_{r'} \rangle \mathbb{E} \langle X_i , \varphi_l \rangle^2 \\
&= \mathbf 1_{\{r=r'\}}\frac{\lambda_r \lambda_{l}}{n}.
\end{align*} 
Therefore, 
\begin{equation}
\label{eq-EA-sum-squares}
\mathbb{E} \left[ A \right] = \frac{1}{n} \sum_{j = 1}^{m_1} \sum_{l = 1}^{m_2} \sum_{r = m_2 + 1}^{+ \infty} \frac{\lambda_l \lambda_r}{(\lambda_l - \lambda_r)^2} \langle  S \Gamma^{1/2} (\varphi_j) , \varphi_r \rangle^2.
\end{equation}
From here, we split the sum above in two parts as follows
\begin{align}
\label{eq-split-EA}
\mathbb{E} \left[ A \right] &= \frac{1}{n} \sum_{j = 1}^{m_1} \sum_{l = 1}^{m_2} \sum_{r = m_2 + 1}^{2 m_2} \frac{\lambda_l \lambda_r}{(\lambda_l - \lambda_r)^2} \langle  S \Gamma^{1/2} (\varphi_j) , \varphi_r \rangle^2 \nonumber \\
&+ \frac{1}{n} \sum_{j = 1}^{m_1} \sum_{l = 1}^{m_2} \sum_{r = 2m_2 + 1}^{+ \infty} \frac{\lambda_l \lambda_r}{(\lambda_l - \lambda_r)^2} \langle  S \Gamma^{1/2} (\varphi_j) , \varphi_r \rangle^2.
\end{align}
To upper bound the first right-side term of Equation \eqref{eq-split-EA}, we begin with 
\begin{align}
\label{eq-interm-EA1}
&\frac{1}{n} \sum_{j = 1}^{m_1} \sum_{l = 1}^{m_2} \sum_{r = m_2 + 1}^{2 m_2} \frac{\lambda_l \lambda_r}{(\lambda_l - \lambda_r)^2} \langle  S \Gamma^{1/2} (\varphi_j) , \varphi_r \rangle^2 \nonumber \\
\leq &\frac{1}{n} \sum_{l = 1}^{m_2} \sum_{u = m_2 + 1}^{2 m_2} \frac{\lambda_l \lambda_u}{(\lambda_l - \lambda_u)^2} \times \sum_{j = 1}^{m_1} \sum_{r = m_2 + 1}^{2 m_2} \langle  S \Gamma^{1/2} (\varphi_j) , \varphi_r \rangle^2.
\end{align}
Now, according to Lemma \ref{lemma-behavior-vap}, we obtain
\begin{eqnarray*}
\sum_{l = 1}^{m_2} \sum_{u = m_2 + 1}^{2 m_2} \frac{\lambda_l \lambda_u}{(\lambda_l - \lambda_u)^2} & = & \sum_{l = 1}^{m_2} \sum_{u = m_2 + 1}^{2 m_2} \frac{\lambda_l^2}{(\lambda_l - \lambda_u)^2} \times \frac{\lambda_u}{\lambda_l} \\
& \leq & \sum_{l = 1}^{m_2} \sum_{u = m_2 + 1}^{2 m_2} \frac{u l}{(u - l)^2} \\
& \leq & \sum_{l = 1}^{m_2} \sum_{u = m_2 + 1}^{2 m_2} \frac{2m_2^2}{(u - l)^2}.
\end{eqnarray*}
By making in the last sum the substitution $v := u - l$, we get
\begin{eqnarray*}
\sum_{l = 1}^{m_2} \sum_{u = m_2 + 1}^{2 m_2} \frac{1}{(u - l)^2} & = & \sum_{l = 1}^{m_2} \sum_{u = 1}^{m_2} \frac{1}{(m_2 + u - l)^2} = \sum_{v = 1 - m_2}^{m_2 - 1} \frac{m_2 - \vert v \vert}{(m_2 + v)^2}\leq  \sum_{v = 1}^{m_2} \frac{1}{v} + \sum_{w = 1}^{m_2 - 1} \frac{m_2 - w}{(m_2 + w)^2} \\
& \leq &  \sum_{v = 1}^{m_2} \frac{1}{v} +1 \leq  1 + \ln(m_2).
\end{eqnarray*}
We turn back to Equation \eqref{eq-interm-EA1} and we write
\begin{align*}
\frac{1}{n} \sum_{j = 1}^{m_1} \sum_{l = 1}^{m_2} \sum_{r = m_2 + 1}^{2 m_2} \frac{\lambda_l \lambda_r}{(\lambda_l - \lambda_r)^2} \langle  S \Gamma^{1/2} (\varphi_j) , \varphi_r \rangle^2 &\leq \frac{C m_2^2 \ln(m_2)}{n} \sum_{j = 1}^{m_1} \sum_{r = m_2 + 1}^{2 m_2} \langle  S \Gamma^{1/2} (\varphi_j) , \varphi_r \rangle^2 \\
&\leq \frac{C m_2^2 \ln(m_2)}{n} \sum_{j = 1}^{m_1} \sum_{r = m_2 + 1}^{2 m_2} \frac{\psi_{\beta} (r)}{\psi_{\beta} (m_2)} \langle  S \Gamma^{1/2} (\varphi_j) , \varphi_r \rangle^2.
\end{align*}
In other words, 
\begin{equation}
\label{eq-upper-bound-EA1}
\frac{1}{n} \sum_{j = 1}^{m_1} \sum_{l = 1}^{m_2} \sum_{r = m_2 + 1}^{2 m_2} \frac{\lambda_l \lambda_r}{(\lambda_l - \lambda_r)^2} \langle  S \Gamma^{1/2} (\varphi_j) , \varphi_r \rangle^2 \leq \frac{C m_2^2 \ln(m_2)}{n \psi_{\beta}(m_2)}.
\end{equation}
We now deal with the second term of Equation \eqref{eq-split-EA}. From Lemma \ref{lemma-behavior-vap}, for $l \leq m_2$ and $r \geq 2 m_2 +1$, we have 
\begin{eqnarray}
\label{eq-upper-bound-f(lambda-lr)}
\frac{\lambda_l^2}{(\lambda_l - \lambda_r)^2} & \leq & \frac{1}{(1 - l/r)^2}  \leq  \frac{1}{\left[ 1 - m_2/(2m_2 + 1) \right]^2} \leq  4.
\end{eqnarray}
Thus, 
\begin{align*}
\frac{1}{n} \sum_{j = 1}^{m_1} \sum_{l = 1}^{m_2} \sum_{r = 2m_2 + 1}^{+ \infty} \frac{\lambda_l \lambda_r}{(\lambda_l - \lambda_r)^2} \langle  S \Gamma^{1/2} (\varphi_j) , \varphi_r \rangle^2 &\leq \frac{4}{n} \sum_{l = 1}^{m_2} \sum_{j = 1}^{m_1} \sum_{r = 2m_2 + 1}^{+ \infty} \frac{\lambda_r}{\lambda_l} \langle  S \Gamma^{1/2} (\varphi_j) , \varphi_r \rangle^2 \\
&\leq \frac{C m_2}{n} \sum_{l = 1}^{m_2} \frac{1}{\lambda_l} \sum_{u = 2m_2 + 1}^{+\infty} \lambda_u \times \sum_{j = 1}^{m_1} \sum_{r = 2m_2 + 1}^{+ \infty} \langle  S \Gamma^{1/2} (\varphi_j) , \varphi_r \rangle^2 \\
&\leq \frac{C m_2}{n} \sum_{l = 1}^{m_2} \frac{\lambda_{2m_2 + 1}}{\lambda_l} \times \sum_{j = 1}^{m_1} \sum_{r = 2m_2 + 1}^{+ \infty} \langle  S \Gamma^{1/2} (\varphi_j) , \varphi_r \rangle^2,
\end{align*}
where the latter inequality is achieved thanks to Lemma \ref{lemma-behavior-vap}. Subsequently,
\begin{align*}
\frac{1}{n} \sum_{j = 1}^{m_1} \sum_{l = 1}^{m_2} \sum_{r = 2m_2 + 1}^{+ \infty} \frac{\lambda_l \lambda_r}{(\lambda_l - \lambda_r)^2} \langle  S \Gamma^{1/2} (\varphi_j) , \varphi_r \rangle^2 &\leq \frac{C m_2^2}{n} \sum_{j = 1}^{m_1} \sum_{r = 2m_2 + 1}^{+ \infty} \langle  S \Gamma^{1/2} (\varphi_j) , \varphi_r \rangle^2 \\
&\leq \frac{C m_2^2}{n \psi_{\beta} (m_2)} \sum_{j = 1}^{m_1} \sum_{r = 2m_2 + 1}^{+ \infty} \psi_{\beta} (r) \langle  S \Gamma^{1/2} (\varphi_j) , \varphi_r \rangle^2.
\end{align*}
Simply put that
\begin{equation}
\label{eq-upper-bound-EA2}
\frac{1}{n} \sum_{j = 1}^{m_1} \sum_{l = 1}^{m_2} \sum_{r = 2m_2 + 1}^{+ \infty} \frac{\lambda_l \lambda_r}{(\lambda_l - \lambda_r)^2} \langle  S \Gamma^{1/2} (\varphi_j) , \varphi_r \rangle^2 \leq \frac{C m_2^2}{n \psi_{\beta} (m_2)}.
\end{equation}
As a result of Equations \eqref{eq-split-EA}, \eqref{eq-upper-bound-EA1} and \eqref{eq-upper-bound-EA2}, we get
\begin{equation}
\label{eq-upper-bound-EA}
\mathbb{E} [A] \leq \frac{C m_2^2 \ln(m_2)}{n \psi_{\beta} (m_2)}.
\end{equation}
We now upper bound the second right-side term of Equation \eqref{eq-EHn(z)-EA+EB}. Similarly to Equation \eqref{eq-EA-sum-squares}, we write
\begin{equation*}
\mathbb{E} [B] = \frac{1}{n} \sum_{l = m_2 + 1}^{+\infty} \sum_{r = 1}^{m_2} \frac{\lambda_r \lambda_l}{(\lambda_r - \lambda_l)^2} \sum_{j=1}^{m_1} \langle S \Gamma^{1/2} (\varphi_j) , \varphi_r \rangle^2.
\end{equation*}
Next, we break the equation down into two parts as shown below
\begin{align}
\label{eq-split-EB}
\mathbb{E} \left[ B \right] &= \frac{1}{n} \sum_{l = m_2 + 1}^{+\infty} \sum_{r = 1}^{\lfloor m_2/2 \rfloor} \frac{\lambda_r \lambda_l}{(\lambda_r - \lambda_l)^2} \sum_{j=1}^{m_1} \langle S \Gamma^{1/2} (\varphi_j) , \varphi_r \rangle^2 \nonumber \\
&+ \frac{1}{n} \sum_{l = m_2 + 1}^{+\infty} \sum_{r = \lfloor m_2/2 \rfloor}^{m_2} \frac{\lambda_r \lambda_l}{(\lambda_r - \lambda_l)^2} \sum_{j=1}^{m_1} \langle S \Gamma^{1/2} (\varphi_j) , \varphi_r \rangle^2.
\end{align}
We now focus on the first right-side term of Equation \eqref{eq-split-EB}. From the assumption $\lambda_r \psi_{\beta} (r) \geq 1$, we obtain
\begin{multline*}
\frac{1}{n} \sum_{l = m_2 + 1}^{+\infty} \sum_{r = 1}^{\lfloor m_2/2 \rfloor} \frac{\lambda_r \lambda_l}{(\lambda_r - \lambda_l)^2} \sum_{j=1}^{m_1} \langle S \Gamma^{1/2} (\varphi_j) , \varphi_r \rangle^2 \\
\leq \frac{1}{n} \sum_{l = m_2 + 1}^{+\infty} \lambda_l \sum_{r = 1}^{\lfloor m_2/2 \rfloor} \frac{\lambda_r^2}{(\lambda_r - \lambda_l)^2} \sum_{j=1}^{m_1} \psi_{\beta} (r) \langle S \Gamma^{1/2} (\varphi_j) , \varphi_r \rangle^2.
\end{multline*}
Similarly to Equation \eqref{eq-upper-bound-f(lambda-lr)}, one can show that for $l > m_2$ and $r \leq m_2/2$, we have
\begin{equation*}
\frac{\lambda_r^2}{(\lambda_r - \lambda_l)^2} \leq C.
\end{equation*}
According to Lemma \ref{lemma-behavior-vap}, we obtain
\begin{equation}
\label{eq-upper-bound-EB1}
\frac{1}{n} \sum_{l = m_2 + 1}^{+\infty} \sum_{r = 1}^{\lfloor m_2/2 \rfloor} \frac{\lambda_r \lambda_l}{(\lambda_r - \lambda_l)^2} \sum_{j=1}^{m_1} \langle S \Gamma^{1/2} (\varphi_j) , \varphi_r \rangle^2 \leq \frac{C m_2 \lambda_{m_2}}{n}.
\end{equation}
We next deal with the second term of Equation \eqref{eq-split-EB}. We start with 
\begin{eqnarray*}
\frac{1}{n} \sum_{l = m_2 + 1}^{+\infty} \sum_{r = \lfloor m_2/2 \rfloor}^{m_2} \frac{\lambda_r \lambda_l}{(\lambda_r - \lambda_l)^2} \sum_{j=1}^{m_1} \langle S \Gamma^{1/2} (\varphi_j) , \varphi_r \rangle^2 \\
&\hspace{-10cm}\leq&\hspace{-5cm} \frac{1}{n \psi_{\beta} (\lfloor m_2/2 \rfloor)} \sum_{l = m_2 + 1}^{+\infty} \lambda_l \times \sum_{r = \lfloor m_2/2 \rfloor}^{m_2}  \frac{\lambda_r}{(\lambda_r - \lambda_{m_2 + 1})^2} \sum_{j=1}^{m_1} \psi_{\beta} (r) \langle S \Gamma^{1/2} (\varphi_j) , \varphi_r \rangle^2 \\
&\hspace{-10cm}\leq&\hspace{-5cm} \frac{C m_2}{n \psi_{\beta} (\lfloor m_2/2 \rfloor)} \sum_{r = \lfloor m_2/2 \rfloor}^{m_2}  \frac{\lambda_r \lambda_{m_2 + 1}}{(\lambda_r - \lambda_{m_2 + 1})^2} \sum_{j=1}^{m_1} \psi_{\beta} (r) \langle S \Gamma^{1/2} (\varphi_j) , \varphi_r \rangle^2.
\end{eqnarray*}
Still using Lemma \ref{lemma-behavior-vap}, for $l > m_2$ and $\lfloor m_2/2 \rfloor \leq r \leq m_2$ we have that
\begin{equation*}
\frac{\lambda_r \lambda_{m_2 + 1}}{(\lambda_r - \lambda_{m_2 + 1})^2} \leq C m_2^2.
\end{equation*}
Consequently,  
\begin{equation}
\label{eq-upper-bound-EB2}
\frac{1}{n} \sum_{l = m_2 + 1}^{+\infty} \sum_{r = \lfloor m_2/2 \rfloor}^{m_2} \frac{\lambda_r \lambda_l}{(\lambda_r - \lambda_l)^2} \sum_{j=1}^{m_1} \langle S \Gamma^{1/2} (\varphi_j) , \varphi_r \rangle^2 \leq \frac{C m_2^3}{n \psi_{\beta} (\lfloor m_2/2 \rfloor)}.
\end{equation}
From Equations \eqref{eq-upper-bound-EB1} and \eqref{eq-upper-bound-EB2}, we conclude that 
\begin{equation}
\label{eq-upper-bound-EB}
\mathbb{E} [B] \leq \frac{C m_2 \lambda_{m_2}}{n} + \frac{C m_2^3}{n \psi_{\beta} (\lfloor m_2/2 \rfloor)}.
\end{equation}
Subsequently, from Equations \eqref{eq-EHn(z)-EA+EB}, \eqref{eq-upper-bound-EA} and \eqref{eq-upper-bound-EB}, we write
\begin{equation}
\label{eq-upper-bound-EHn(z)}
\mathbb{E} \Vert \mathcal{H}_n (z) S \Pi_{m_1} \Gamma^{1/2} \Vert^2_{\HS} \leq \frac{C m_2 \lambda_{m_2}}{n} + \frac{C m_2^2 \ln(m_2)}{n \psi_{\beta} (m_2)} + \frac{C m_2^3}{n \psi_{\beta} (\lfloor m_2/2 \rfloor)}.
\end{equation}
To upper bound the second right-side term of Equation \eqref{eq-decomp-term3-Hn-Gn}, we draw inspiration from the proof of Lemmma 15 of \cite{brunel2016non}. We write
\begin{align*}
\Vert \mathcal{G}_n (z) S \Pi_{m_1} \Gamma^{1/2} \Vert_{\HS}^2 &= \frac{1}{4 \pi^2} \sum_{j = 1}^{m_1} \Vert \mathcal{G}_n (z) S \Gamma^{1/2} (\varphi_j) \Vert^2 \\
&\leq \frac{1}{4 \pi^2} \sum_{j = 1}^{m_1} \left( \sum_{k = 1}^{m_2} \int_{\mathcal{B}_k} \Vert R^{1/2}(z) [I - T_n(z)]^{-1} T_n(z)^2  R^{1/2}(z) S \Gamma^{1/2} (\varphi_j) \Vert \right)^2
\end{align*}
We introduce the diagonal operator $P_{\beta}$ defined for all $j$ in $\mathbb{N}\backslash\{0\}$ by 
$$P_{\beta} (\varphi_j) = \psi_{\beta} (j)^{1/2} \varphi_j.$$ 
Then, 
\begin{multline*}
\Vert \mathcal{G}_n (z) S \Pi_{m_1} \Gamma^{1/2} \Vert_{\HS}^2 
\\ \leq \frac{1}{4 \pi^2} \sum_{j = 1}^{m_1} \left( \sum_{k = 1}^{m_2} \int_{\mathcal{B}_k} \Vert R^{1/2}(z) \Vert_{\infty} \Vert [I - T_n(z)]^{-1} \Vert_{\infty} \Vert T_n(z) \Vert_{\infty}^2 \Vert R^{1/2}(z) P_{\beta}^{-1} \Vert_{\infty} \Vert P_{\beta} S \Gamma^{1/2} (\varphi_j) \Vert \right)^2
\end{multline*}
From here we use some results of \cite{brunel2016non}, page 226. In particular, on the set $\mathcal{A}_n$ and for $z$ in $\mathcal{B}_k$, we write
$$\Vert [I - T_n(z)]^{-1} \Vert_{\infty} < 2 \quad \mbox{ and } \quad \Vert T_n(z) \Vert_{\infty} \leq \frac{\mathbf{a}_k}{\sqrt{n}} \ln(n).$$ 
Also, remark that for $z$ in $\mathcal{B}_k$, we have 
\begin{align*}
\Vert R^{1/2}(z) P_{\beta}^{-1} \Vert_{\infty} &= \sup_{l \in \mathbb{N}\backslash\{0\}} \psi_{\beta} (l)^{-1/2} \vert z - \lambda_l \vert^{- 1/2} \\
&= \psi_{\beta} (k)^{-1/2} \sqrt{2/\delta_k}.
\end{align*}
Hence, 
\begin{align*}
\Vert \mathcal{G}_n (z) S \Pi_{m_1} \Gamma^{1/2} \Vert_{\HS}^2 &\leq \frac{C \ln^4 (n)}{n^2} \sum_{j = 1}^{m_1} \Vert P_{\beta} S \Gamma^{1/2} (\varphi_j) \Vert^2  \left( \sum_{k = 1}^{m_2} \mathbf{a}_k^2 \psi_{\beta} (k)^{-1/2} \right)^2 \\
&= \frac{C \ln^4 (n)}{n^2} \Vert P_{\beta} S \Gamma^{1/2} \Vert_{\HS}^2  \left( \sum_{k = 1}^{m_2} \mathbf{a}_k^2 \psi_{\beta} (k)^{-1/2} \right)^2. 
\end{align*} 
Besides, according to Lemma \ref{lemma-behavior-ak}, under Assumptions $\boldsymbol{\mathcal{A}_{2}}$ and $\boldsymbol{\mathcal{A}_{6}}$, we have $\mathbf{a}_k \leq C k \ln(k)$. Furthermore, 
\begin{equation*}
\Vert P_{\beta} S \Gamma^{1/2} \Vert_{\HS}^2 = \sum_{j = 1}^{+\infty} \sum_{l = 1}^{+\infty} \psi_{\beta} (l) \langle S \Gamma^{1/2} (\varphi_j) , \varphi_l \rangle^2 < + \infty.
\end{equation*}
Therefore, 
\begin{equation}
\label{eq-upper-bound-EGn(z)}
\mathbb{E} \Vert \mathcal{G}_n (z) S \Pi_{m_1} \Gamma^{1/2} \Vert_{\HS}^2 \leq \frac{C \ln^4 (n)}{n^2} \left( \sum_{k = 1}^{m_2} \frac{k^2 \ln^2 (k)}{\sqrt{\psi_{\beta} (k)}} \right)^2. 
\end{equation} 
To end this proof, we provide an upper bound for the remaining term of Equation \eqref{eq-decomp-term3-Hn-Gn}. It is plain that
\begin{equation*}
\mathbb{E} \Vert (\widehat{\Pi}_{m_2} - \Pi_{m_2}) S \Pi_{m_1} \Gamma^{1/2} \bm{1}_{\mathcal{A}_n^{\complement}} \Vert_{\HS}^2 \leq 4 \Vert S \Gamma^{1/2} \Vert_{\HS}^2 \mathbb{P} ( \bm{1}_{\mathcal{A}_n^{\complement}}).
\end{equation*}
It follows from Lemma \ref{lemma-hatPim-Pim} that
\begin{equation}
\label{eq-upper-bound-last-term-bias}
\mathbb{E} \Vert (\widehat{\Pi}_{m_2} - \Pi_{m_2}) S \Pi_{m_1} \Gamma^{1/2} \bm{1}_{\mathcal{A}_n^{\complement}} \Vert_{\HS}^2 \leq \frac{C}{n^2} \Vert S \Gamma^{1/2} \Vert_{\HS}^2.
\end{equation}
At last, combining Equations \eqref{eq-decomp-term3-Hn-Gn}, \eqref{eq-upper-bound-EHn(z)},\eqref{eq-upper-bound-EGn(z)} and \eqref{eq-upper-bound-last-term-bias} leads to the required result.

\subsubsection{Proof of Corollary \ref{corr-upper-bound-MSPE}}

From Equation \eqref{eq-thm-upper-bound-MSPE}, the minimax prediction rate is lower than the infimum of the right-side term with respect to $m_1$ and $m_2$. Let us consider the terms of the upper-bound successively.

\begin{itemize}
\item We keep the first one, $\sigma^2_{\varepsilon} m_1/n$.
\item For the second term, under Assumption $\mathbf{\mathcal{A}_1}$, we obtain
\begin{align*}
\sum_{j = m_1 + 1}^{+\infty} \Vert S \Gamma^{1/2} (\varphi_j) \Vert^2 &= \sum_{j = m_1 + 1}^{+\infty} \sum_{r = 1}^{+\infty} \langle S \Gamma^{1/2}(\varphi_j) , \varphi_r \rangle^2 \\
&\leq \sum_{j = m_1 + 1}^{+\infty} \sum_{r = 1}^{+\infty} \frac{\eta_{\alpha} (j)}{\eta_{\alpha} (m_1)} \langle S \Gamma^{1/2}(\varphi_j) , \varphi_r \rangle^2 \\
&\leq c/\eta_{\alpha} (m_1),
\end{align*}
where $c$ is some positive constant.
\item The third term, $\sum_{j=1}^{m_1}\|(I-\Pi_{m_2})S\Gamma^{1/2}(\varphi_j)\|^2$ goes to zero when $m_2$ goes to $+ \infty$.
\item We have $A_{n,m_1}\leq \sigma_{\varepsilon}^2m_1/n$ under the constraint $m_1\leq n/\ln^2(n)$. 
\item The additional assumption $\lambda_{m_1}\leq m_1^{-1-\nu}$ permits to obtain that $B_{n,m_1}$ is negligible with respect to $\sigma_{\varepsilon}^2m_1/n$. 
\item The term $E_n$ is immediately negligible with respect to $1/n$.
\item From the convergence of $\sum_{m_2 \geq 1} \lambda_{m_2}$, the first term of $D_{n,m_2}$ goes to 0 when $m_2$ goes to $+\infty$. The same result applies to the second and third terms of $D_{n,m_2}$ in view of the assumptions on $\psi_{\beta}$.
\item The last term of $D_{n,m_2}$ is lower to $c \ln^4 (n)/n^2$ (and thus to $c/n$). Indeed, from the assumptions on $\psi_{\beta}$, the serie $\sum_{k \geq 1} k^2 \ln^2 (k) / \sqrt{\psi_{\beta} (k)}$ is convergent. 
\end{itemize}

\subsubsection{Proof of Theorem \ref{thm-lower-bound-MSPE}}

To lower bound the minimax risk, we follow the general scheme of reduction to a finite hypotheses number, as described in \cite{tsybakov2008introduction}. Our approach is similar to the one applied in \cite{crambes2013asymptotics}, but the regularity assumption differs, making it necessary to adapt the proof significantly. We start by considering a family of hypothesis $S^{\theta}$ of $S$, indexed by $\theta$ in $\Omega_{m_1^*} = \{ 0 , 1 \}^{m_1^*}$. More precisely, for all $\theta = (\omega_1, \ldots, \omega_{m_1^*})$ in $\Omega_{m_1^*}$, we define 
$$S^{\theta} = \sum_{j = 1}^{m_1^*} \mu_j \omega_j \varphi_1 \otimes \varphi_j,$$ 
where $(\mu_j)_{1 \leq j \leq m_1^*}$ are chosen in such a way that $S^{\theta} \Gamma^{1/2}$ belongs to $\mathcal{W}_{\alpha, \beta}^{R}$. Elementary computations show that we can set for all $j$ in $\Omega_{m_1^*}$, 
$$\mu_j^2 = R^2/e \times \left[ \lambda_j m_1^* \eta_{\alpha} (m_1^*) \right]^{-1}.$$
Besides, it is straightforward to see that
\begin{equation}
\label{eq-minimax-MSPE-finit-hyp}
\inf_{\widehat S_{n}} \sup_{S \Gamma^{1/2} \in \mathcal{W}_{\alpha,\beta}^R} \mathbb{E} \Vert \widehat S_{n} (X_{n+1}) - S (X_{n+1}) \Vert^2\geq \inf_{\widehat S_{n}} \max_{\theta \in \Omega_{m_1^*}} \mathbb{E} \Vert \widehat S_{n} (X_{n+1}) - S^{\theta} (X_{n+1}) \Vert^2.
\end{equation}
For a given estimator $\widehat S_{n}$, let $\widehat{\theta} (\widehat S_{n})$ be a random vector verifying that $$\widehat{\theta} (\widehat S_{n}) \in \argmin\limits_{\theta \in \Omega_{m_1^*}} \mathbb{E} \Vert \widehat S_{n} (X_{n+1}) - S^{\theta} (X_{n+1}) \Vert^2.$$ The model $S^{\widehat{\theta} (\widehat S_{n})}$ is one of the nearest to $\widehat S_{n}$, among the collection $\left\{ S^{\theta}; \theta \in \Omega_{m_1^*} \right\}$. Then, for all $\theta$ in $\Omega_{m_1^*}$, 
\begin{align*}
\mathbb{E} \Vert S^{\widehat{\theta} (\widehat S_{n})} (X_{n+1}) - S^{\theta} (X_{n+1}) \Vert^2 &\leq 2 \mathbb{E} \Vert \widehat S_{n} (X_{n+1}) - S^{\widehat{\theta} (\widehat S_{n})} (X_{n+1}) \Vert^2 + 2 \mathbb{E} \Vert \widehat S_{n} (X_{n+1}) - S^{\theta} (X_{n+1}) \Vert^2 \\
&\leq 4 \mathbb{E} \Vert \widehat S_{n} (X_{n+1}) - S^{\theta} (X_{n+1}) \Vert^2.
\end{align*}
Therefore, using Equation \eqref{eq-minimax-MSPE-finit-hyp}, we get
\begin{align}
\label{eq-minimax-MSPE-theta-hat}
\inf_{\widehat S_{n}} \sup_{S \Gamma^{1/2} \in \mathcal{W}_{\alpha,\beta}^R} \mathbb{E} \Vert \widehat S_{n} (X_{n+1}) - S (X_{n+1}) \Vert^2 &\geq \frac{1}{4} \inf_{\widehat S_{n}} \max_{\theta \in \Omega_{m_1^*}} \mathbb{E} \Vert S^{\widehat{\theta} (\widehat S_{n})} (X_{n+1}) - S^{\theta} (X_{n+1}) \Vert^2 \nonumber \\
&\geq \frac{1}{4} \inf_{\widehat{\theta}} \max_{\theta \in \Omega_{m_1^*}} \mathbb{E} \Vert S^{\widehat{\theta}} (X_{n+1}) - S^{\theta} (X_{n+1}) \Vert^2,
\end{align}
where the infimum in the last line is taken over the all the estimators with values in $\Omega_{m_1^*}$. From here, we aim to apply Assouad's Lemma with Kullback version, see \citet[pages 117 to 119]{tsybakov2008introduction}. For the sake of completeness, we recall below this lemma. 
\begin{lemma}[Assouad version Kullback-Leibler]
Let $\rho$ denotes the Hamming distance, defined for all $\theta = (\omega_1, \ldots, \omega_{m})$ and $\theta' = (\omega'_1, \ldots, \omega'_{m})$ in $\Omega_{m} = \{0,1\}^m$ by
$$\rho(\theta, \theta') = \sum_{j = 1}^{m} \mathds{1}_{\{\omega_j \neq \widehat{\omega}_j \}}.$$
We also denote by $\mathbb{P}^{\otimes n}_{\theta}$ the distribution of $(X_i, Y_i)_{1 \leq i \leq n}$ under $S^{\theta}$. Assume that for all $\theta$, $\theta'$ such that $\rho (\theta, \theta') = 1$, we have the following upper bound of the Kullback-Leibler divergence between $\mathbb{P}_{\theta}^{\otimes n}$ and $\mathbb{P}_{\theta'}^{\otimes n}$,
$\KL (\mathbb{P}_{\theta}^{\otimes n}, \mathbb{P}_{\theta'}^{\otimes n}) \leq \alpha < \infty$. Then, 
$$\inf_{\widehat{\theta}} \max_{\theta \in \Omega_{m}} \mathbb{E} [ \rho (\theta, \widehat{\theta})] \geq \frac{m}{2} \max \left( \frac{1}{2} \exp(-\alpha), 1 - \sqrt{\alpha/2} \right),$$
where the infimum is taken over all the estimators $\widehat{\theta}$ with values in $\Omega_m$.
\end{lemma} 

\smallskip
Now, for all estimator $\widehat{\theta} = (\widehat{\omega}_1, \ldots, \widehat{\omega}_{m_1^*})$ and $\theta = (\omega_1, \ldots, \omega_{m_1^*})$ both with values in $\Omega_{m_1^*}$, we write
\begin{align*}
\mathbb{E} \Vert S^{\widehat{\theta}} (X_{n+1}) - S^{\theta} (X_{n+1}) \Vert^2 &= \mathbb{E} \Vert \sum_{j = 1}^{m_1^*} \mu_j (\widehat{\omega}_j - \omega_j ) \varphi_1 \otimes \varphi_j (X_{n+1}) \Vert^2 \\
&= \mathbb{E} \Vert \sum_{j = 1}^{m_1^*} \mu_j (\widehat{\omega}_j - \omega_j ) \varphi_1 \otimes \varphi_j \Gamma^{1/2} \Vert_{\HS}^2,
\end{align*} 
where the last equation stems from Lemma \ref{lemma-ps-Xn+1-HS}. Thus,
\begin{align*}
\mathbb{E} \Vert S^{\widehat{\theta}} (X_{n+1}) - S^{\theta} (X_{n+1}) \Vert^2 &= \sum_{j = 1}^{m_1^*} \lambda_j \mu_j^2 \mathbb{E} \left[ (\widehat{\omega}_j - \omega_j )^2 \right] \\
&= R^2/e \times [m_1^* \eta_{\alpha} (m_1^*)]^{-1} \mathbb{E} [ \rho (\theta, \widehat{\theta}) ].
\end{align*}
We then get from Equation \eqref{eq-minimax-MSPE-theta-hat} that
\begin{equation}
\label{eq-inf-Hamming}
\inf_{\widehat S_{n}} \sup_{S \Gamma^{1/2} \in \mathcal{W}_{\alpha,\beta}^R} \mathbb{E} \Vert \widehat S_{n} (X_{n+1}) - S (X_{n+1}) \Vert^2 \geq \frac{R^2}{4e} [m_1^* \eta_{\alpha} (m_1^*)]^{-1} \inf_{\widehat{\theta}} \max_{\theta \in \Omega_{m_1^*}} \mathbb{E} [ \rho (\theta, \widehat{\theta}) ].
\end{equation}
In order to apply Assouad's Lemma, we now upper bound $\KL (\mathbb{P}^{\otimes n}_{\theta}, \mathbb{P}^{\otimes n}_{\theta'})$ for $\theta$, $\theta'$ such that $\rho (\theta, \theta') = 1$. We are aware that the definition of the last Kullback-Leibler divergence requires that $\mathbb{P}^{\otimes n}_{\theta}$ is absolutely continuous with respect to $\mathbb{P}^{\otimes n}_{\theta'}$. This point will be clarified later. Furthermore, it is well known that 
$$\KL (\mathbb{P}^{\otimes n}_{\theta}, \mathbb{P}^{\otimes n}_{\theta'}) = n \KL (\mathbb{P}_{\theta}, \mathbb{P}_{\theta'}),$$
where $\mathbb{P}_{\theta}$ is the distribution of $(X, Y^{\theta})$, with $Y^{\theta} = S^{\theta} X  + \varepsilon$. Then, 
\begin{align*}
\KL (\mathbb{P}^{\otimes n}_{\theta}, \mathbb{P}^{\otimes n}_{\theta'}) &= n \int \ln \left( \frac{\mathrm{d} \mathbb{P}_{\theta}}{\mathrm{d} \mathbb{P}_{\theta'}} \right) \mathrm{d} \mathbb{P}_{\theta} \\
&= n \int \left[ \int \ln \left( \frac{\mathrm{d} \mathbb{P}_{Y^{\theta} | X}}{\mathrm{d} \mathbb{P}_{Y^{\theta'} | X}} \right) \mathrm{d} \mathbb{P}_{Y^{\theta} | X} \right] \mathrm{d} \mathbb{P}_{X}.
\end{align*}
The absolute continuity of $\mathbb{P}_{Y^{\theta} | X}$ with respect to $\mathbb{P}_{Y^{\theta'} | X}$ and the upper bound of the Kullback-Leibler divergence are provided by Cameron–Martin Theorem \citep{lifshits2012lectures}, which is stated below. 
\begin{theorem}[Cameron–Martin]
Let $Z$ be a centered Gaussian random variable in a Hilbert space $( \mathcal{X}, \langle \cdot , \cdot \rangle, \Vert \cdot \Vert )$, with a distribution measure $P$ and a covariance operator $\Gamma_Z$. We consider the subset $H_P \subset \mathcal{X}$ defined as $$H_P = \left\{ h \in \mathcal{X} \mbox{ such that } \Vert \Gamma_Z^{-1/2} (h) \Vert^2 < + \infty \right\}.$$
For all $h$ in $H_P$, we denote by $P_h$ the distribution mesure of $Z+h$. Then, $P_h$ is absolutely continuous with respect to $P$ and the density $\mathrm{d} P_h/\mathrm{d} P$ is given by
$$\mathrm{d} P_h/\mathrm{d} P : x \mapsto \exp \left\{ \langle x , \Gamma^{-1}_Z (h) \rangle - \frac{1}{2}\Vert \Gamma_Z^{-1/2} (h) \Vert^2 \right\}.$$
\end{theorem}
The next step is to apply Cameron-Martin Theorem to upper bound $\KL (\mathbb{P}^{\otimes n}_{\theta}, \mathbb{P}^{\otimes n}_{\theta'})$. Before this, there are some points to specify. We denote for the remaining of the proof $\Gamma_{\varepsilon}$, the covariance operator of $\varepsilon$. It is quite easy to see that conditionally to $X$, the distribution of $Y_{\theta}$ is Gaussian with mean $S^{\theta} X$ and covariance operator $\Gamma_{Y_{\theta}} = \Gamma_{\varepsilon}$. Let us set $h_{\theta} = S^{\theta} X$, for all $\theta$ in $\Omega_{m_1^*}$. Then, conditionally to $X$, the random variable $Y^{\theta}$ is a shift of $Y^{\theta'}$ with a shift equal to $h^{\theta} - h^{\theta'}$. Therefore, according to Cameron-Martin Theorem, we have 
\begin{align*}
\KL (\mathbb{P}^{\otimes n}_{\theta}, \mathbb{P}^{\otimes n}_{\theta'}) &= n \mathbb{E}_X \mathbb{E}_{\varepsilon} \langle Y^{\theta} , \Gamma^{-1}_{\varepsilon} (h^{\theta} - h^{\theta'}) \rangle - \frac{n}{2} \mathbb{E}_X \Vert \Gamma_{\varepsilon}^{-1/2} (h^{\theta} - h^{\theta'}) \Vert^2 \\
&= n \mathbb{E}_X \langle h^{\theta} , \Gamma^{-1}_{\varepsilon} (h^{\theta} - h^{\theta'}) \rangle - \frac{n}{2} \mathbb{E}_X \Vert \Gamma_{\varepsilon}^{-1/2} (h^{\theta} - h^{\theta'}) \Vert^2,
\end{align*}
where $\mathbb{E}_X$ and $\mathbb{E}_{\varepsilon}$ respectively designate the expectations with respect to the distributions of $X$ and $\varepsilon$. Then, according to Lemma \ref{lemma-ps-Xn+1-HS}, we get
\begin{equation*}
\KL (\mathbb{P}^{\otimes n}_{\theta}, \mathbb{P}^{\otimes n}_{\theta'}) = n  \langle S^{\theta} \Gamma^{1/2} , \Gamma^{-1}_{\varepsilon} (S^{\theta} - S^{\theta'}) \Gamma^{1/2} \rangle_{\HS} - \frac{n}{2} \Vert \Gamma_{\varepsilon}^{-1/2} (S^{\theta} - S^{\theta'}) \Gamma^{1/2} \Vert_{\HS}^2.
\end{equation*}
Now, from $\rho(\theta, \theta') = 1$, we know that there exists $j$ in $\{1, \ldots, m_1^* \}$ such that $\omega_j \neq \omega'_j$ and for all $l \neq j$, $\omega_l = \omega'_l$. Thus, 
\begin{align*}
\KL (\mathbb{P}^{\otimes n}_{\theta}, \mathbb{P}^{\otimes n}_{\theta'}) &= n \lambda_j \mu_j^2 \omega_j (\omega_j - \omega'_j) \langle \varphi_1, \Gamma^{-1}_{\varepsilon} (\varphi_1) \rangle - \frac{n}{2} \lambda_j \mu_j^2 \omega_j^2 \Vert \Gamma_{\varepsilon}^{-1/2} (\varphi_1) \Vert^2 \\
&= \frac{R^2}{2e \sigma_1} \omega_j (\omega_j - 2 \omega'_j) \times \frac{n}{m_1^* \eta_{\alpha} (m_1^*)} \\
&\leq \frac{R^2}{2e \sigma_1} \times \frac{n}{m_1^* \eta_{\alpha} (m_1^*)},
\end{align*}
where $1/\sigma_1 = \Vert \Gamma_{\varepsilon}^{-1/2} (\varphi_1) \Vert^2$. Bearing in mind that $m_1^* \eta_{\alpha} (m_1^*) \geq n/c_0$, for some positive constant $c_0$. By fitting $c_0$ value with respect to $R$, $\sigma_1$ and $\alpha$, we set the condition $ R^2/(2e \sigma_1) \times 1/c_0 (R, \sigma_1, \alpha) \leq 1/2$. Consequently, from Assouad's Lemma and Equation \eqref{eq-inf-Hamming}, we get
\begin{equation*}
\inf_{\widehat S_{n}} \sup_{S \Gamma^{1/2} \in \mathcal{W}_{\alpha,\beta}^R} \mathbb{E} \Vert \widehat S_{n} (X_{n+1}) - S (X_{n+1}) \Vert^2 \geq \frac{C(R,\sigma_1, \alpha)}{\eta_{\alpha} (m_1^*)},
\end{equation*}
where $C(R,\sigma_1, \alpha)$ is some positive constant only depending on its arguments. In other words, 
\begin{equation*}
\inf_{\widehat S_{n}} \sup_{S \Gamma^{1/2} \in \mathcal{W}_{\alpha,\beta}^R} \mathbb{E} \Vert \widehat S_{n} (X_{n+1}) - S (X_{n+1}) \Vert^2 \geq C(R,\sigma_1, \alpha) \inf_{m_1 \in \mathbb{N}\backslash\{0\}} \left\{ \sigma_{\varepsilon}^2 \frac{m_1}{n} +  \frac{3}{\eta_{\alpha} (m_1)} \right\}.
\end{equation*}

\subsubsection{Proof of Corollary \ref{coro:minimax_rates}}

In both cases (polynomial and exponential) we have to compute $\inf_{m_1\in\mathbb N\backslash\{0\}}\{\sigma_{\varepsilon}^2m_1/n+3/{\eta_{\alpha}(m_1)}\}$. 
We distinguish the two cases. 

\begin{enumerate}
\item If $\eta_{\alpha}(x)\asymp x^{\alpha}$, the dimension $m_1$ minimizing the minimax risk is up to a multiplicative constant, the solution $x^*_n=\arg\min_{x\geq 1}R_\alpha(x)$ with $R_\alpha(x)=x/n+x^{-\alpha}$.  we easily derive the result, $x_n^{*}=cn^{1/{\alpha+1}}$, for a constant $c$ which only depends on $\alpha$, and $R_{\alpha}(x^*_n)=cn^{-\alpha/(\alpha+1)}$, for another constant $c$.

\item If $\eta_{\alpha}(x)\asymp \exp(x^\alpha)$, we reason in a similar way with $R_\alpha(x)=x/n+\exp(-x^\alpha)$. The solution $x_n^*$ of the minimization problem of $R_\alpha$ verifies the equation $\alpha(x^*_n)^{\alpha-1}e^{-(x^*_n)^\alpha}=n^{-1}$ which cannot be written explicitly except in the case $\alpha= 1$ where $x^*_n= \ln(n)$ leading to a minimax rate of order $\ln(n)/n$. In the case $\alpha\neq 1$, we can upper-bound the risk as follows, for $n\geq e^1$  
\[
R_{\alpha}(x^*_n) =\frac{x^*_n}{n}+e^{-(x^*_n)^\alpha}\leq R_{\alpha}\left((\ln(n))^{1/\alpha}\right)=\frac{\ln(n)^{1/\alpha}}{n}+ \frac1n\leq 2\frac{(\ln(n))^{1/\alpha}}n
\]
implying that there exists $c >0$ such that, for $n$ large enough, $e^{-(x^*_n)^\alpha}\leq 2(\ln(n))^{1/\alpha}/n$, which implies
\begin{align*}x^*_n&\geq(\ln(n)-\ln((\ln(n))^{1/\alpha}-\ln(2))_+^{1/\alpha},\\
&= (\ln(n))^{1/\alpha}\left(1-\frac{1}{\alpha}\frac{\ln(\ln(n))}{n}-\frac{\ln(2)}{\ln(n)}\right)^{1/\alpha},\\
&\geq c(\ln(n))^{1/\alpha}.
\end{align*}
Moreover, 
$$R_{\alpha}(x^*_n)  =\frac{x^*_n}{n}+e^{-(x^*_n)^\alpha}\geq\frac{x^*_n}n\geq c\frac{(\ln(n))^{1/\alpha}}n.$$ 
Then the minimax rate is larger than $(\ln(n))^{1/\alpha}/n$.

\end{enumerate}

\subsection{Proof of the results of Section \ref{sec:adaptation}}\label{sec:proof_sec_adapt}

\subsubsection{Proof of Theorem \ref{thm-UB-Emp-Risk}}

Starting from the definition of $\widehat{m}_1$, we have for all $m_1$ in $\mathcal{M}_n$ that
\begin{equation}
\label{eq-Contrast-Pen}
\gamma_n (\widehat{S}_{\widehat{m}_1, \infty})  \leq   \gamma_n (\widehat{S}_{m_1, \infty}) +\pen(m_1) - \pen(\widehat{m}_1).
\end{equation}
Moreover, $\widehat{S}_{m_1, \infty}$ is the limit, in the sense of the Hilbert-Schmidt norm of $\widehat{S}_{m_1,m_2}$ when $m_2$ goes to infinity. Since the contrast function $\gamma_n : \mathcal{L}_2(\mathbb H)\rightarrow \mathbb R$ is continuous (where  we recall that $ \mathcal{L}_2(\mathbb H)$ is the space of Hilbert-Schmidt operators of $\mathbb H=L^2([0,1])$), then 
$$\gamma_n (\widehat{S}_{m_1, \infty}) =\lim_{m_2\to\infty}\gamma_n(\widehat{S}_{m_1,m_2})\leq \lim_{m_2\to\infty}\gamma_{n}(\widehat\Pi_{m_1,m_2}S)=\gamma_n(\widehat{\Pi}_{m_1,\infty}S),$$
where the inequality is a consequence of the definition of the estimators $(\widehat S_{m_1,m_2})_{m_1,m_2}$, see \eqref{eq:estim_MC}, and where $\widehat\Pi_{m_1,m_2}$ (resp. $\widehat\Pi_{m_1,\infty}$) stands for the projection operator onto the space $\Span\{\widehat\varphi_k\otimes\widehat\varphi_j,\;j=1,\ldots,m_1,\, k=1,\ldots m_2\}$  (resp. onto $\Span\{\widehat\varphi_k\otimes\widehat\varphi_j,\;j=1,\ldots,m_1,\, k\in\mathbb N\backslash\{0\}\}$).
Thus, since $\widehat{\Pi}_{m_1,\infty}S=S\widehat\Pi_{m_1}$ (where we recall that $\widehat{\Pi}_{m_1}$ is the projection operator of functions onto $\Span\{\widehat\varphi_j,\,j=1\ldots,m_1\}$, see \eqref{eq:egalite_proj}), Inequality \eqref{eq-Contrast-Pen} becomes
\begin{equation}
\label{eq-Contrast-Pen2}
\gamma_n (\widehat{S}_{\widehat{m}_1, \infty})  \leq   \gamma_n (S\widehat\Pi_{m_1}) +\pen(m_1) - \pen(\widehat{m}_1).
\end{equation}

We now seek for the link between the contrast function and the empirical norm. Let $T_1$ and $T_2$ be two linear operators. Then, 
\begin{align*}
\gamma_n (T_1) - \gamma_n (T_2) &= \frac{1}{n} \sum_{i = 1}^n \Vert Y_i - T_1 (X_i) \Vert^2 - \frac{1}{n} \sum_{i = 1}^n \Vert Y_i - T_2 (X_i) \Vert^2 \\
&= \frac{1}{n} \sum_{i = 1}^n \Vert [S - T_1] (X_i) + \varepsilon_i \Vert^2 - \frac{1}{n} \sum_{i = 1}^n \Vert [S - T_2] (X_i) +\varepsilon_i \Vert^2 \\
&= \Vert S - T_1 \Vert_n^2 - \Vert S - T_2 \Vert_n^2 + 2 \nu_n (T_2 - T_1), 
\end{align*} 
where the empirical process $\nu_n$ is defined for any linear operator $T$ by $\nu_n (T) = n^{-1} \sum_{i = 1}^n \langle T(X_i) , \varepsilon_i \rangle$. Back to Equation \eqref{eq-Contrast-Pen}, for all $m_1$ in $\mathcal{M}_n$ we write
\begin{align*}
\Vert S - \widehat{S}_{\widehat{m}_1, \infty} \Vert_n^2  &\leq  \Vert S -  S\widehat\Pi_{m_1} \Vert_n^2 + 2 \nu_n (\widehat{S}_{\widehat{m}_1, \infty} -  S\widehat\Pi_{m_1})  + \pen(m_1) - \pen(\widehat{m}_1) \\
&\leq \Vert S -  S\widehat\Pi_{m_1} \Vert_n^2 + 2 \Vert \widehat{S}_{\widehat{m}_1, \infty} -  S\widehat\Pi_{m_1} \Vert_n  \sup_{\substack{T \in V_{m \vee \widehat{m}_1, \infty} \\ \Vert T \Vert_{n} = 1}} \nu_n (T) + \pen(m_1) - \pen(\widehat{m}_1),
\end{align*}
where for all $m$ in $\mathbb{N}\backslash\{0\}$, $V_{m,\infty} = \Span \{ \widehat{\varphi}_k \otimes \widehat{\varphi}_j, j = 1 , \ldots m \mbox{ and } k \in \mathbb{N}\backslash\{0\} \}$. 
Let, $\zeta>0$ and $\tilde\zeta=4\zeta^{-1}+2$. Since $\tilde\zeta>0$, we remark that for all $x,y$ in $\mathbb{R}$, $2xy \leq \tilde\zeta x^2 + \tilde\zeta^{-1}y^2$. Then,
\begin{equation*}
\Vert S - \widehat{S}_{\widehat{m}_1, \infty} \Vert_n^2 \leq \Vert S - S\widehat\Pi_{m_1} \Vert_n^2 +\tilde \zeta^{-1}\| \widehat{S}_{\widehat{m}_1, \infty} -  S\widehat\Pi_{m_1} \Vert_n^2 + \tilde\zeta \sup_{\substack{T \in V_{m_1 \vee \widehat{m}_1, \infty} \\ \Vert T \Vert_{n} = 1}} \nu_n (T)^2 + \pen(m_1) - \pen(\widehat{m}_1).
\end{equation*}
We define for all $m,m'$ in $\mathcal{M}_n$, $p(m,m') = 8\tilde\zeta^{-1}(1+ \delta) \sigma^2_{\varepsilon} (m \vee m')/n$ such that $\pen(m) + \pen(m') \geq \tilde\zeta p(m,m')$. Therefore,
\begin{equation*}
\Vert S - \widehat{S}_{\widehat{m}_1, \infty} \Vert_n^2 \leq \Vert S -  S\widehat\Pi_{m_1} \Vert_n^2 + \tilde\zeta^{-1}\Vert \widehat{S}_{\widehat{m}_1, \infty} - S\widehat\Pi_{m_1} \Vert_n^2 + 2 \pen(m_1) + \tilde\zeta \left( \sup_{\substack{T \in V_{m_1 \vee \widehat{m}_1, \infty} \\ \Vert T \Vert_{n} = 1}} \nu_n (T)^2 - p(m_1,\widehat{m}_1) \right)_+,
\end{equation*}
where $(\cdot)_+$ is the positive part, defined for all $x$ in $\mathbb{R}$ as $x_+ = x \vee 0$. Next, since 
\[
\Vert \widehat{S}_{\widehat{m}_1, \infty} - S\widehat\Pi_{m_1} \Vert_n^2\leq 2 \Vert \widehat{S}_{\widehat{m}_1, \infty} - S \Vert_n^2+2\Vert S - S\widehat\Pi_{m_1} \Vert_n^2
\]
and $\tilde\zeta>2$
\begin{equation}
\label{eq-pre-Talagrand}
(1-2\tilde\zeta^{-1})\Vert S - \widehat{S}_{\widehat{m}_1, \infty} \Vert_n^2 \leq (1+2\tilde\zeta^{-1}) \Vert S -  S\widehat\Pi_{m_1}\Vert_n^2 + 2 \pen(m_1) + \tilde\zeta\left( \sup_{\substack{T \in V_{m_1 \vee \widehat{m}_1, \infty} \\ \Vert T \Vert_{n} = 1}} \nu_n (T)^2 - p(m_1,\widehat{m}_1) \right)_+.
\end{equation}
The upper-bound of the expectation of the last right-sided term in Equation \eqref{eq-pre-Talagrand} stems from Lemma \ref{lemma-Talagrand}. This upper-bound does not depend on the selected dimension $m$ in $\mathcal{M}_n$. Hence,  denoting $c(\zeta)=2(1+2\tilde\zeta^{-1})^{-1}=(2+\zeta)/(1+\zeta)$ we have,
\begin{equation*}
\mathbb{E} \Vert S - \widehat{S}_{\widehat{m}_1, \infty} \Vert_n^2 \leq  \frac{1+2\tilde\zeta^{-1}}{1-2\tilde\zeta^{-1}}\inf_{m_1 \in \mathcal{M}_n} \left\{ \mathbb{E} \Vert S -  S\widehat\Pi_{m_1} \Vert_n^2 + c(\zeta)\pen(m_1) \right\} + \frac{C(\delta,p)}{n}
\end{equation*}
and the result comes from the fact that $\frac{1+2\tilde\zeta^{-1}}{1-2\tilde\zeta^{-1}}=1+\zeta$.

\begin{lemma}
\label{lemma-Talagrand}
Under Assumption $\mathcal{A}_{7}$, the following inequality holds for all $m$ in $\mathcal{M}_n$,
\begin{equation*}
\sum_{m' \in \mathcal{M}_n} \mathbb{E} \left[ \left( \sup_{\substack{T \in V_{m \vee m', \infty} \\ \Vert T \Vert_{n} = 1}} \nu_n^2 (T) - p(m,m') \right)_+ \right] \leq \frac{C(\delta,p)}{n}, 	
\end{equation*}
where $p(m,m') = 2(1+ \delta) \frac{m \vee m'}{n} \sigma^2$ and for all $x$ in $\mathbb{R}$, $x_+ = x \vee 0$. 
\end{lemma}

\medskip
\noindent\textit{Proof of Lemma \ref{lemma-Talagrand}.} The main ideas of this proof are inspired from \citet[Lemma 5]{brunel2016non}, with an extension of the demonstrated results in the infinite dimensional frame of this paper. Denote by $\bar{T}_X = \left( T(X_1), \ldots, T(X_n) \right) \in \mathbb{H}^{n}$ and $\bar{\varepsilon} = \left( \varepsilon_1, \ldots, \varepsilon_n \right) \in \mathbb{H}^n$. Then, $$\nu_n(T) = \frac{1}{n} \langle \bar{T}_X, \bar{\varepsilon} \rangle^{\otimes n},$$ where $\langle \cdot , \cdot \rangle^{\otimes n}$ is the scalar product on $\mathbb{H}^{n}$ induced by $\langle \cdot , \cdot \rangle$, which is defined for all $x = (x_1, \ldots, x_n)$ and $y = (y_1, \ldots, y_n)$ in $\mathbb{H}^{n}$ as $\langle x , y \rangle^{\otimes n} = \sum_{i = 1}^{n} \langle x_i , y_i \rangle$. The associated norm is denoted $\Vert \cdot \Vert^{\otimes n}$. Subsequently, 
\begin{equation*}
\sup_{\substack{T \in V_{m, \infty} \\ \Vert T \Vert_{n} = 1 }} \nu_n (T) = \sup_{\substack{z \in \widetilde{V}_{m, \infty} \\ \Vert z \Vert^{\otimes n} = \sqrt{n}}} \frac{1}{n} \langle z , \bar{\varepsilon} \rangle^{\otimes n},
\end{equation*}
where $\widetilde{V}_{m, \infty}$ is the linear subspace of $\mathbb{H}^{n}$ defined as
$$\widetilde{V}_{m, \infty} = \left\{ z \in \mathbb{H}^{n}, \exists T \in V_{m, \infty} \mbox{ such that } z = \bar{T}_X \right\}.$$ Therefore, denoting by $\Pi^{\otimes n}_{\tilde{V}_{m,\infty}}$ the orthogonal projection onto $\widetilde{V}_{m, \infty}$, we get
\begin{align*}
\sup_{\substack{T \in V_{m, \infty} \\ \Vert T \Vert_{n} = 1 }} \nu_n (T) &= \frac{1}{\sqrt{n}} \sup_{\substack{z \in \widetilde{V}_{m,\infty} \\ \Vert z \Vert^{\otimes n} = 1}} \langle z , \bar{\varepsilon} \rangle^{\otimes n} = \frac{1}{\sqrt{n}} \Vert \Pi^{\otimes n}_{\tilde{V}_{m,\infty}} (\bar{\varepsilon}) \Vert^{\otimes n}.
\end{align*}
From here, we use the Talagrand-type inequality stated in Lemma \ref{lemma-Adapt-Baraud}, which is an adaptation of \citep[Corollary 5.1]{baraud2000model} in our setting. For all $x > 0$, we write
\begin{equation*}
\mathbb{P} \left( n \sup_{\substack{T \in V_{m, \infty} \\ \Vert T \Vert_{n} = 1 }} \nu_n^2 (T) \geq m \sigma_{\varepsilon}^2 + 2 \sigma_{\varepsilon}^2 \sqrt{mx} +  \sigma_{\varepsilon}^2 x \right) \leq C(p) \sigma_{\varepsilon}^p \tau_p \frac{m}{x^{p/2}},
\end{equation*}
where  $\tau_p = \mathbb{E} \Vert \varepsilon \Vert^{p}$. Now, for all $\theta > 0$, notice that $2\sqrt{mx} \leq \theta m + \theta^{-1} x$, this means that
\begin{equation*}
\mathbb{P} \left(\sup_{\substack{T \in V_{m, \infty} \\ \Vert T \Vert_{n} = 1 }} \nu_n^2 (T) \geq (1 + \theta) \frac{m \sigma_{\varepsilon}^2}{n} + (1 + \theta^{-1}) \frac{\sigma_{\varepsilon}^2 x}{n} \right) \leq C(p) \frac{m}{x^{p/2}},
\end{equation*}
We set, 
$$Q_{m \vee m'} = \left( \sup_{\substack{T \in V_{m \vee m', \infty} \\ \Vert T \Vert_{n} = 1 }} \nu_n^2 (T) - p(m,m') \right)_+.$$
We write for all $m,m'$ in $\mathcal{M}_n$, 
\begin{align*}
	\mathbb{E} \left[ Q_{m \vee m'} \right]	&= \int_{0}^{+ \infty} \mathbb{P} \left( Q_{m \vee m'} \geq t \right) \mathrm{d} t \\
	&\leq C(\delta, p) \frac{m \vee m'}{n^{p/2}} \int_{0}^{+ \infty} \frac{\mathrm{d} t}{\left( t + \sigma_{\varepsilon}^2 m \vee m' /n (1 + \delta) \right)^{p/2}} \leq C(\delta,p) (m \vee m')^{2 - p/2}/n. 
\end{align*}	
Knowing that $p > 4$, we have that $(m \vee m')^{2 - p/2} \leq 1$ which gives the desired result.


\begin{lemma}
\label{lemma-Adapt-Baraud}
Let $L$ be a linear subspace of $\mathbb{H}^n$ and $m \in \mathbb{N}\backslash\{0\}$ such that $dim(L) \leq m$. Denote by $\zeta$ the map defined on $\mathbb{H}^{n}$ as $$\zeta : s \mapsto \Vert \Pi^{\otimes n}_{L} (s) \Vert^{\otimes n},$$ where $\Pi^{\otimes n}_{L}$ designates the orthogonal projection onto $L$ with respect to $\langle \cdot , \cdot \rangle^{\otimes n}$. Then, for all $p \geq 2$ such that $\tau_p = \mathbb{E} \Vert \varepsilon \Vert^p < + \infty$ and for all $x > 0$, we have
$$\mathbb{P} \left( \zeta^2 (\bar{\varepsilon}) \geq m \sigma_{\varepsilon}^2 + 2 \sigma_{\varepsilon}^2 \sqrt{mx} + \sigma_{\varepsilon}^2 x \right) \leq C(p)\frac{m}{x^{p/2}}.$$ 
\end{lemma}

\medskip
\noindent\textit{Proof of Lemma \ref{lemma-Adapt-Baraud}.} Let $(\psi_1, \ldots, \psi_{\ubar{m}})$ be an orthonormal basis of $L$, with $\ubar{m} \leq m$. Denote by $\mathcal{B}_n$ the unit ball of $\mathbb{H}^{n}$ with respect to the norm $\Vert \cdot \Vert^{\otimes n}$, we have 
\begin{align*}
\zeta^2 (\bar{\varepsilon}) &= \left[ \sup_{u \in \mathcal{B}_n} \langle \Pi^{\otimes n}_{L} (\bar{\varepsilon}) , u \rangle^{\otimes n} \right]^2 = \left[ \sup_{u \in \mathcal{B}_n} \langle  \bar{\varepsilon} , \Pi^{\otimes n}_{L} (u) \rangle^{\otimes n} \right]^2 \\
&= \left[ \sup_{u \in \mathcal{B}_n} \sum_{j = 1}^{\ubar{m}} \langle \bar{\varepsilon} , \psi_j \rangle^{\otimes n} \langle u , \psi_j \rangle^{\otimes n} \right]^2. 
\end{align*}
Let $\upsilon_i$ be the $n$-tuple of $\mathbb{H}^n$ composed by $\varepsilon_i$ in the $i^{\text{th}}$ position and $0_{\mathbb{H}}$ elsewhere. By construction, $(\upsilon_i)_i$ are independent and $\bar{\varepsilon} = \sum_{i = 1}^n \upsilon_i$ and if we denote for $u$ in $\mathcal{B}_n$, $g_u : x \mapsto \sum_{j = 1}^{\ubar{m}} \langle x , \psi_j \rangle^{\otimes n} \langle u , \psi_j \rangle^{\otimes n}$, we get 
$$\zeta^2 (\bar{\varepsilon}) = \left[ \sup_{u \in \mathcal{B}_n} \sum_{i = 1}^n g_u (\upsilon_i) \right]^2.$$
Now, we rely on \citet[Theorem 5.2]{baraud2000model} and we write for all $t > 0$, 
\begin{align}
\label{eq-interm-concentr}
\mathbb{P} \left( \zeta (\bar{\varepsilon}) \geq \mathbb{E} \left[ \zeta (\bar{\varepsilon}) \right] + t \right) &\leq t^{-p} \mathbb{E} \vert \zeta (\bar{\varepsilon}) - \mathbb{E} \left[ \zeta (\bar{\varepsilon}) \right] \vert^{p} \nonumber \\
&\leq C(p) t^{-p} \left( \mathbb{E}_1 +  \mathbb{E}_2^{p/2} \right),
\end{align}
where $\mathbb{E}_1$ and $\mathbb{E}_2$ are defined as
$$\mathbb{E}_1 = \mathbb{E} \left[ \max_{i = 1, \ldots, n} \sup_{u \in \mathcal{B}_n} \vert g_u (\upsilon_i) \vert^p \right] \quad \mbox{and} \quad \mathbb{E}_2 = \mathbb{E} \left[ \sup_{u \in \mathcal{B}_n} \sum_{i = 1}^n g_u^2 (\upsilon_i) \right].$$
By Cauchy-Schwarz, it is quite straightforward to see that
\begin{align*}
\mathbb{E}_1 &= \mathbb{E} \left[ \max_{i = 1, \ldots, n} \sup_{u \in \mathcal{B}_n} \vert \sum_{j = 1}^{\ubar{m}} \langle \varepsilon_i , \psi_{j,i} \rangle \langle u , \psi_j \rangle^{\otimes n} \vert^p \right] \\
&\leq \mathbb{E} \left[ \sum_{i = 1}^n \Vert \varepsilon_i \Vert^p \sup_{u \in \mathcal{B}_n} \Vert \sum_{j = 1}^{\ubar{m}} \langle u , \psi_j \rangle^{\otimes n} \psi_{j,i} \Vert^p \right],
\end{align*}
where $\psi_j = (\psi_{j,1}, \ldots, \psi_{j,n})$ for all $j$ in $\{ 1, \ldots, \ubar{m} \}$. Notice that, $\Vert \sum_{j = 1}^{\ubar{m}} \langle u , \psi_j \rangle^{\otimes n} \psi_{j,i} \Vert^2 \leq 1$ for all $u$ in $\mathcal{B}_n$ and $i$ in $\{1, \ldots , n \}$. This leads to   
\begin{align}
\label{eq-upper-bound-E1}
\mathbb{E}_1 &\leq \tau_p \times \sum_{i = 1}^n \sup_{u \in \mathcal{B}_n} \Vert \sum_{j = 1}^{\ubar{m}} \langle u , \psi_j \rangle^{\otimes n} \psi_{j,i} \Vert^2 \nonumber \\
&\leq \tau_p \times \sum_{i = 1}^n \sup_{u \in \mathcal{B}_n} \left\{ \sum_{l = 1}^{\ubar{m}} \vert \langle u , \psi_l \rangle^{\otimes n} \vert^2 \sum_{j = 1}^{\ubar{m}} \Vert \psi_{j,i} \Vert^2 \right\} \nonumber \\
&\leq \tau_p \times \sum_{i = 1}^n \sum_{j = 1}^{\ubar{m}} \Vert \psi_{j,i} \Vert^2 \nonumber \\
&\leq m \tau_p.
\end{align}
To upper bound $\mathbb{E}_2$, we start with 
\begin{align}
\label{eq-interm-upper-bound-E2}
\mathbb{E}_2 &= \mathbb{E} \left[ \sup_{u \in \mathcal{B}_n} \sum_{i = 1}^n \vert \sum_{j = 1}^{\ubar{m}} \langle \varepsilon_i , \psi_{j,i} \rangle^{\otimes n} \langle u , \psi_j \rangle^{\otimes n} \vert^2 \right] \nonumber \\
&\leq \mathbb{E} \left[ \sup_{u \in \mathcal{B}_n} \sum_{i = 1}^n \Vert \varepsilon_i \Vert^2 \Vert \sum_{j = 1}^{\ubar{m}} \langle u , \psi_j \rangle^{\otimes n} \psi_{j,i} \Vert^2  \right] \nonumber \\
&\leq \mathbb{E}_{2,1} + \mathbb{E}_{2,2},
\end{align}
where $\mathbb{E}_{2,1}$ and $\mathbb{E}_{2,2}$ are defined below
\begin{align*}
\mathbb{E}_{2,1} &= \mathbb{E} \left[ \sup_{u \in \mathcal{B}_n} \sum_{i = 1}^n \Vert \varepsilon_i \Vert^2 \Vert \sum_{j = 1}^{\ubar{m}} \langle u , \psi_j \rangle^{\otimes n} \psi_{j,i} \Vert^2 \bm{1}_{\Vert \varepsilon_i \Vert \leq c} \right], \\ 
\mathbb{E}_{2,2} &= \mathbb{E} \left[ \sup_{u \in \mathcal{B}_n} \sum_{i = 1}^n \Vert \varepsilon_i \Vert^2 \Vert \sum_{j = 1}^{\ubar{m}} \langle u , \psi_j \rangle^{\otimes n} \psi_{j,i} \Vert^2 \bm{1}_{\Vert \varepsilon_i \Vert > c}  \right],
\end{align*}
with $c > 0$ will be set latter. It is clear that
\begin{equation}
\label{eq-upper-bound-E21}
\mathbb{E}_{2,1} \leq c^2. 
\end{equation}
On the flip side, 
\begin{align}
\label{eq-upper-bound-E22}
\mathbb{E}_{2,2} &\leq c^{2-p} \times \mathbb{E} \left[ \sup_{u \in \mathcal{B}_n} \sum_{i = 1}^n \Vert \varepsilon_i \Vert^p \Vert \sum_{j = 1}^{\ubar{m}} \langle u , \psi_j \rangle^{\otimes n} \psi_{j,i} \Vert^2 \right] \nonumber \\
&\leq c^{2-p} m \tau_p, 
\end{align}
where the last inequality holds with the same arguments when upper bounding $\mathbb{E}_1$. Taking $c^p = m \tau_p$ and combining Equations \eqref{eq-interm-upper-bound-E2}, \eqref{eq-upper-bound-E21} and \eqref{eq-upper-bound-E22} yields
\begin{equation}
\label{eq-upper-bound-E2}
\mathbb{E}_2^{p/2} \leq m \tau_p.
\end{equation}
From Equations \eqref{eq-interm-concentr}, \eqref{eq-upper-bound-E1} and \eqref{eq-upper-bound-E2}, we get for all $t > 0$ that
\begin{equation*}
\mathbb{P} \left( \zeta (\bar{\varepsilon}) \geq \mathbb{E} [ \zeta (\bar{\varepsilon}) ] + t \right) \leq C(p)m \tau_p/t^p.
\end{equation*}
Now, the well-known inequality $\mathbb{E} [ \zeta (\bar{\varepsilon}) ]^2 \leq  \mathbb{E} [ \zeta^2 (\bar{\varepsilon}) ]$ gives
\begin{equation*}
\mathbb{P} \left( \zeta^2 (\bar{\varepsilon}) \geq \mathbb{E} [ \zeta^2 (\bar{\varepsilon}) ] + 2 \sqrt{\mathbb{E} [\zeta^2 (\bar{\varepsilon})] t^2} + t^2 \right) \leq C(p)m \tau_p/t^p.
\end{equation*}
Starting from the last inequality, showing that $\mathbb{E} [\zeta^2 (\bar{\varepsilon})] \leq m \sigma_{\varepsilon}^2$ and taking $t^2 = x$ ends the proof. Indeed, knowing that the $(\varepsilon_i)_i$ are independent and centered, we have 
\begin{align*}
\mathbb{E} [\zeta^2 (\bar{\varepsilon})] &= \sum_{j = 1}^{\ubar{m}} \mathbb{E} \left[ \left( \sum_{i = 1}^n \langle \varepsilon_i , \psi_{j,i} \rangle \right)^2 \right] \\
&\leq \sum_{j = 1}^{\ubar{m}} \sum_{i = 1}^n \mathbb{E} \langle \varepsilon_i , \psi_{j,i} \rangle^2  \\
&\leq m \sigma_{\varepsilon}^2, 
\end{align*}
where the last upper-bound stems from Cauchy-Schwarz inequality. 

\section*{Acknowledgments}
The authors want to thank André Mas for his careful reading of our work.

\section*{Funding}
The research of the authors is partly supported by the french Agence Nationale de la Recherche (ANR-18-CE40-0014, projet SMILES).

\bibliographystyle{apalike}
\bibliography{Ref}

\end{document}